\documentclass[10pt,reqno]{amsart}

\usepackage[utf8]{inputenc}
\usepackage[T1]{fontenc}
\usepackage{amsmath,amssymb,amsthm}
\usepackage{mathrsfs}
\usepackage{color}
\usepackage[margin=1.0in]{geometry}
\usepackage[colorlinks=false, pdfborder={0 0 0}] {hyperref}  

\allowdisplaybreaks

\numberwithin{equation}{section}

\theoremstyle{plain}
\newtheorem{theorem}{Theorem}[section]
\newtheorem{lemma}[theorem]{Lemma}
\newtheorem{proposition}[theorem]{Proposition}
\newtheorem{corollary}[theorem]{Corollary}

\theoremstyle{definition}
\newtheorem{remark}{Remark}[section]

\newcommand{\R}{\mathbb{R}}
\newcommand{\T}{\mathbb{T}}
\newcommand{\dv}{\operatorname{div}}
\newcommand{\eps}{\varepsilon}
\newcommand{\Mloc}{\mathcal{M}}
\newcommand{\Mstar}{\mathcal{M}_{*}}
\newcommand{\Lop}{\mathcal{L}}
\newcommand{\Nspace}{\mathfrak{N}}
\newcommand{\intv}{\int_{\R^{3}}}

\newcommand{\tp}{\widetilde{p}}
\newcommand{\tE}{\widetilde{E}}

\newcommand{\ba}{\begin{aligned}}
\newcommand{\ea}{\end{aligned}}
\newcommand{\be}{\begin{equation}}
\newcommand{\ee}{\end{equation}}
\newcommand{\nm}{\nonumber}

\def\di{\displaystyle}
\def\r{\rho}
\def\wt{\widetilde}

\def\v{\varepsilon}
\def\mb{\mathbf}
\def\f{\frac}
\def\t{\theta}

\def\d{\delta}
\def\s{\sigma}

\def\z{\zeta}
\def\na{\nabla}
\def\pa{\partial}

\begin{document}

\title[Large-time behavior for 3D fluid-particle system]
      {Large-time behavior and stability of planar rarefaction waves for the three-dimensional full fluid-particle system}

\author{Hailiang Li}
\address{School of Mathematical Sciences, Capital Normal University,
         Beijing 100048, P.~R.~China}
\email{hailiang-li@mail.cnu.edu.cn}

\author{Dehua Wang}
\address{Department of Mathematics, University of Pittsburgh,
         Pittsburgh, PA 15260, USA}
\email{dhwang@pitt.edu}

\author{Teng Wang}
\address{Department of Mathematics, School of Mathematics, Statistics and
         Mechanics, Beijing University of Technology, Beijing 100124,
         P.~R.~China}
\email{tengwang@amss.ac.cn}

\date{\today}

\begin{abstract}
The large-time behavior of the three-dimensional full compressible Euler-Vlasov-Fokker-Planck system is studied. In this system, an inviscid polytropic gas and a dispersed particle phase exchange both momentum and energy. We prove the nonlinear time-asymptotic stability of a planar rarefaction wave in a domain that is unbounded in one direction and periodic in the two transverse directions, for initial data that are a small third-order Sobolev perturbation of a spatially varying local Maxwellian determined by the wave profile. For the non-isentropic case, we are not aware of any earlier stability result for a nonlinear wave pattern: the global results available concern perturbations of a global Maxwellian around a constant state, and the wave-pattern analysis carried out so far is for the isentropic model. The proofs rely on three key features of the energy exchange. A microscopic H-theorem on the five-dimensional macroscopic space controls the heat-flux moments of the microscopic component. An effective nonconservative viscous system for the particle density, velocity, and temperature is derived via a Chapman-Enskog expansion about the local Maxwellian, with viscosity and heat conductivity that depend on both the density and the temperature. A double relaxation mechanism dissipates the relative velocity and the temperature difference between the two phases simultaneously; thus, the momentum and energy estimates must be closed together. In particular, the kinetic coupling by itself provides enough dissipation to stabilize a nonlinear wave of the inviscid fluid.
\end{abstract}

\maketitle

%%%%%%%%%%%%%%%%%%%%%%%%%%%%%%%%%%%%%%%%%%%%%%%%%%%%%%%%%%%%%%%%%%%%%%%%%%%%%%
\section{Introduction}
 
Compressible fluid-particle models describe clouds of small particles (droplets, dust grains, or bubbles) dispersed in a carrier gas, and are used in combustion theory \cite{Williams1958,Williams1985}, spray propulsion \cite{MassotVilledieu}, biomedical sprays \cite{Baranger2005}, and pollutant-dispersion modeling \cite{Sportisse}. In the thin-spray regime, the particle volume fraction is negligible, but the mass fraction is not. The dispersed phase is represented statistically by a distribution function $f(t,x,v)$ satisfying a Vlasov-Fokker-Planck equation and exerting drag on the gas.  When particles and gas exchange momentum and thermal energy, their temperatures relax toward each other, yielding the full compressible Euler-Vlasov-Fokker-Planck (Euler-VFP) system:
\begin{equation}\label{EVFP}
\begin{cases}
\rho_{t}+\dv_{x}(\rho u)=0,\\[2pt]
(\rho u)_{t}+\dv_{x}(\rho u\otimes u)+\nabla_{x}p
   =\rho\displaystyle\intv (v-u)f\,dv,\\[6pt]
(\rho E)_{t}+\dv_{x}(\rho uE+pu)
   =\rho\displaystyle\intv\big[(v-u)\cdot v-3\theta\big]f\,dv,\\[6pt]
f_{t}+v\cdot\nabla_{x}f
   =\rho\,\dv_{v}\big((v-u)f+\theta\nabla_{v}f\big).
\end{cases}
\end{equation}
The fluid state is described by the density $\rho=\rho(t,x)\ge 0$, the velocity $u=u(t,x)=(u_{1},u_{2},u_{3})(t,x)$, the total energy $E(t,x)\ge0$ and the
absolute temperature $\theta(t,x)\ge0$, for $t\ge0$ and $x\in\Omega:=\R\times\T^{2}$, where $\T^{2}$ is the two-dimensional torus. The dispersed phase is described by $f=f(t,x,v)\ge0$ with $v\in\R^{3}$ the particle velocity. The pressure $p$, the total energy $E$ and the internal energy $e$ obey the constitutive relations:
\begin{equation}\label{constitutive}
p=R\rho\theta=A\rho^{\gamma}\exp\Big(\frac{\gamma-1}{R}S\Big),\qquad
E=e+\frac{|u|^{2}}{2},\qquad
e=\frac{R}{\gamma-1}\theta+\mathrm{const.},
\end{equation}
where $S$ is the entropy, $\gamma>1$ is the adiabatic exponent and $A$, $R$ are positive fluid constants. Throughout the paper we consider the polytropic
monatomic case $\gamma=5/3$ and normalize $R=1$. The model \eqref{EVFP} with temperature transport was derived in \cite{Boudin2009}.

We consider \eqref{EVFP} on $\Omega=\R\times\T^{2}$, with periodic boundary conditions in the transverse directions and initial data
\begin{equation}\label{initial}
\begin{cases}
\big(\rho,u,\theta,f\big)(t,x,v)\big|_{t=0}
   =\big(\rho_{0}(x),u_{0}(x),\theta_{0}(x),f_{0}(x,v)\big),\quad x\in\Omega,\\[4pt]
(\rho_{0},u_{0},\theta_{0})(x)\to(\rho_{\pm},u_{\pm},\theta_{\pm}),
   \qquad\text{as }x_{1}\to\pm\infty,\\[4pt]
f_{0}(x,v)\to \Mloc_{[n_{\pm},u_{\pm},\theta_{\pm}]}(v)
  :=\dfrac{n_{\pm}}{\sqrt{(2\pi\theta_{\pm})^{3}}}\,
    e^{-\frac{|v-u_{\pm}|^{2}}{2\theta_{\pm}}},
   \qquad\text{as }x_{1}\to\pm\infty,
\end{cases}
\end{equation}
where $\rho_{\pm}>0$, $n_{\pm}>0$, $u_{\pm}=(u_{1\pm},0,0)$ and $\theta_{\pm}>0$ are prescribed constants.

We note that the Fokker-Planck operator
\[
\Lop_{[\rho,u,\theta]}f=\rho\,\dv_{v}\big((v-u)f+\theta\nabla_{v}f\big)
\]
has density-dependent friction and temperature-dependent diffusion. Its equilibria are Maxwellians with the local fluid velocity and temperature, and its relaxation rate is proportional to the fluid density. Therefore, a global Maxwellian is not a natural reference profile when the fluid state is non-constant.
Set $n=\intv f\,dv$ and $nw=\intv vf\,dv$ as the particle density and momentum. Integrating $\eqref{EVFP}_{4}$ against $1,\,v,\,\tfrac12|v|^{2}$ and using $\intv\dv_{v}\big((v-u)f+\theta\nabla_{v}f\big)dv=0$, one finds that the coupling terms
\begin{equation}\label{exchange}
\rho\intv(v-u)f\,dv=\rho n(w-u),
\quad
\rho\intv\big[(v-u)\cdot v-3\theta\big]f\,dv
   =3\rho n(T-\theta)+\rho n\,w\cdot(w-u)
\end{equation}
appear in the particle equations with the opposite sign, so that the total momentum $\rho u+nw$ and the total energy of the mixture are conserved. Here $T$ is the particle temperature, defined in \eqref{npart} below. Thus \eqref{exchange} contains two coupled relaxation mechanisms, one acting on the relative velocity $w-u$ and the other on the temperature difference $T-\theta$. In the barotropic fluid-particle models, only the first of these is present.

%When the carrier fluid is viscous and heat-conducting, one considers the full compressible Navier-Stokes-Vlasov-Fokker-Planck (NS-VFP) system instead:
%\begin{equation}\label{NSVFP}
%\begin{cases}
%\rho_{t}+\dv_{x}(\rho u)=0,\\[2pt]
%(\rho u)_{t}+\dv_{x}(\rho u\otimes u)+\nabla_{x}p
%   =\dv_{x}\big(2\mu D(u)+\lambda\,\dv_{x}\!u\,\mathbb{I}\big)
%    +\rho\displaystyle\intv (v-u)f\,dv,\\[6pt]
%(\rho E)_{t}+\dv_{x}(\rho uE+pu)
%   =\kappa\Delta_{x}\theta
%    +\dv_{x}\big(u(2\mu D(u)+\lambda\,\dv_{x}\!u\,\mathbb{I})\big)\\[2pt]
%\hspace{5.35cm}
%    +\rho\displaystyle\intv\big[(v-u)\cdot v-3\theta\big]f\,dv,\\[6pt]
%f_{t}+v\cdot\nabla_{x}f
%   =\rho\,\dv_{v}\big((v-u)f+\theta\nabla_{v}f\big),
%\end{cases}
%\end{equation}
%where $D(u)=(\nabla_{x}u+(\nabla_{x}u)^{\top})/2$ is the deformation tensor, $\mathbb{I}$ is the $3\times3$ identity matrix, the shear and bulk viscosities
%satisfy $\mu>0$, $3\mu+2\lambda\ge0$, and $\kappa>0$ is the heat conductivity.

When the carrier fluid is viscous and heat-conducting, one considers the full compressible Navier-Stokes-Vlasov-Fokker-Planck (NS-VFP) system, where the fluid equations already provide the spatial dissipation needed to stabilize a nonlinear wave.
In this paper, we focus on the inviscid system \eqref{EVFP},  %In \eqref{NSVFP}, the fluid equations already provide the spatial dissipation needed to stabilize a nonlinear wave. 
where this dissipation is absent and must instead be generated by the kinetic equation and interphase relaxation. We study how dissipation is transferred from the kinetic phase to the hyperbolic phase and show that this transfer is sufficient to stabilize a planar rarefaction wave.

\subsection{Related results}
\label{SSprevious}

 For the kinetic-fluid systems without energy exchange, one considers the barotropic NS-VFP system with density-independent friction.
%\begin{equation}\label{NSVFPbaro}
%\begin{cases}
%\rho_{t}+\dv_{x}(\rho u)=0,\\[2pt]
%(\rho u)_{t}+\dv_{x}(\rho u\otimes u)+\nabla_{x}p
%  =\dv_{x}\big(2\mu D(u)+\lambda\,\dv_{x}\!u\,\mathbb{I}\big)
%   +\displaystyle\intv (v-u)f\,dv,\\[6pt]
%f_{t}+v\cdot\nabla_{x}f=\dv_{v}\big((v-u)f+\nabla_{v}f\big).
%\end{cases}
%\end{equation}
Mellet-Vasseur \cite{MelletVasseur2007,MelletVasseur2008} proved the existence and the asymptotic behavior of weak solutions in a three-dimensional bounded
domain; Chae-Kang-Lee \cite{ChaeKangLee} obtained global classical solutions and their exponential decay on $\T^{3}\times\R^{3}$; Wang \cite{WangWJ}
constructed global strong solutions in $H^{2}$ together with optimal decay rates for the spatial derivatives. Dropping the viscosity in
%\eqref{NSVFPbaro} 
the barotropic NS-VFP system leads to the barotropic Euler-VFP system, for which Carrillo-Goudon \cite{CarrilloGoudon} analyzed the asymptotic regimes and the
hydrodynamic expansions, and Duan-Liu \cite{DuanLiu} established global well-posedness and convergence to equilibrium under small perturbations. When
the friction force is proportional to the fluid density $\rho$, the coupling generates additional nonlinear terms and, as observed already in
\cite{MelletVasseur2007}, the global theory becomes substantially more delicate;
see \cite{BarangerDesvillettes,CaoJiang,ChoiKwon,LiMuWang,LiNiWang,LiNiWu2,WangYu} and the references therein. In particular, Li-Mu-Wang \cite{LiMuWang}
and Li-Ni-Wu \cite{LiNiWu2} obtained global classical solutions with optimal decay rates near a global Maxwellian, and Li-Ni-Wang \cite{LiNiWang} proved
global well-posedness together with the inviscid limit.

For the kinetic-fluid systems with energy exchange, there are fewer rigorous results in the non-isentropic case. Numerical studies of fluid-particle models
with energy exchange were carried out in \cite{Boudin2009,GoudonJinYan}. Mu-Wang \cite{MuWang} proved global existence and the optimal large-time
behavior of strong solutions to the three-dimensional non-isentropic NS-VFP system  %\eqref{NSVFP}
 in $H^{4}$, both in the whole space and in periodic
domains; Li-Ni-Wu \cite{LiNiWu1} lowered the initial regularity from $H^{4}(\R^{3})$ to $H^{2}(\R^{3})$ and derived optimal decay rates; very
recently Li-Ni-Xin \cite{LiNiXin} constructed global classical solutions of the non-isentropic NS-VFP system %\eqref{NSVFP} 
and analyzed the vanishing dissipation limit, thereby
obtaining global classical solutions of the non-isentropic Euler-VFP system. % \eqref{EVFP}. 
See also \cite{CaoJiang,LiLiuYang,LiShou}.
The results above concern perturbations of a global Maxwellian around a constant fluid state. They do not address the stability of spatially varying local Maxwellians or nonlinear wave patterns such as rarefaction waves, shock profiles, and contact discontinuities. For fluid-particle systems, the only result on wave patterns known to us is the isentropic analysis of Li-Wang-Wang \cite{LiWangWang2022}. However, the argument in \cite{LiWangWang2022} does not extend to the non-isentropic case, as explained in Section \ref{SSnovelty} below.

Regarding the planar rarefaction waves in several space dimensions, for multi-D viscous conservation laws, Xin \cite{Xin} first proved the stability of planar
rarefaction waves in the scalar case; decay rates were obtained by Ito \cite{Ito} and Nishikawa-Nishihara \cite{NishikawaNishihara}. Hokari-Matsumura
\cite{HokariMatsumura} treated a $2\times2$ artificial system with a positive-definite viscosity matrix in two dimensions. For the compressible
Navier-Stokes equations with periodic transverse directions, the time-asymptotic stability of planar rarefaction waves was established by
Li-Wang \cite{LiWang2018} and Li-Wang-Wang \cite{LiWangWang2018}, and Huang-Xu-Yuan \cite{HuangXuYuan} obtained the stability of planar rarefaction
waves for the three-dimensional Navier-Stokes equations under periodic perturbations. At the kinetic level, planar rarefaction waves were studied for
the Boltzmann equation by Wang-Wang \cite{WangWang2019} and for the bipolarVlasov-Poisson-Boltzmann system in \cite{LiWangWang2018bipolar,
WangWang2020}. 
In all these works, the macroscopic equations are viscous, unlike those considered in this paper.

\subsection{Main result}
\label{SSmain}

We are concerned with the nonlinear time-asymptotic stability of a planar rarefaction wave for the inviscid Euler-VFP system \eqref{EVFP}. To this end,
for a solution $f$ of $\eqref{EVFP}_{4}$ the particle density $n$, mean
velocity $w$ and temperature $T$ are defined by
\begin{equation}\label{npart}
\begin{split}
n(t,x)=\intv f\,dv,\quad
n w_{i}(t,x)=\intv v_{i}f\,dv\  (i=1,2,3),\quad
n\Big(\frac{3}{2}T+\frac{|w|^{2}}{2}\Big)(t,x)=\intv\frac{|v|^{2}}{2}f\,dv,
\end{split}
\end{equation}
and the associated local Maxwellian is
\begin{equation}\label{localM}
\Mloc=\Mloc_{[n,w,T]}(t,x,v)
=\frac{n(t,x)}{\sqrt{\big(2\pi T(t,x)\big)^{3}}}
 \exp\Big(-\frac{|v-w(t,x)|^{2}}{2T(t,x)}\Big).
\end{equation}
Following the macro-micro decomposition of \cite{LiuYangYu,LiuYu} and applying a Chapman-Enskog expansion to $\eqref{EVFP}_{4}$ around $\Mloc$, we obtain
at the formal fluid level the \emph{nonconservative two-phase system}:
\begin{equation}\label{twophase}
\begin{cases}
\rho_{t}+\dv_{x}(\rho u)=0,\\[2pt]
(\rho u)_{t}+\dv_{x}(\rho u\otimes u)+\nabla_{x}p=\rho n(w-u),\\[2pt]
(\rho E)_{t}+\dv_{x}(\rho uE+pu)=3\rho n(T-\theta)+\rho n\,w\cdot(w-u),\\[4pt]
n_{t}+\dv_{x}(nw)=0,\\[2pt]
(nw)_{t}+\dv_{x}(nw\otimes w)+\nabla_{x}\tp
   =\dv_{x}\Big(\dfrac{nT}{2\rho}\,\mathbb{S}(w)\Big)-\rho n(w-u),\\[6pt]
(n\tE)_{t}+\dv_{x}(nw\tE+\tp w)
   =\dv_{x}\Big(\dfrac{5nT}{6\rho}\nabla_{x}T\Big)
   +\dv_{x}\Big(\dfrac{nT}{2\rho}\,w\cdot\mathbb{S}(w)\Big)\\[6pt]
\hspace{4.6cm}-\,3\rho n(T-\theta)-\rho n\,w\cdot(w-u),
\end{cases}
\end{equation}
where $\tE=\frac32T+\frac{|w|^{2}}{2}$, the pressures are $p=\rho\theta$ and
$\tp=nT$, and the effective viscous stress tensor is
\begin{equation}\label{stress}
\mathbb{S}(w):=2D(w)-\frac{2}{3}\dv_{x}\!w\,\mathbb{I},
\qquad
\big(\mathbb{S}(w)\big)_{ij}
 =\partial_{x_{j}}w_{i}+\partial_{x_{i}}w_{j}
  -\frac{2}{3}\,\dv_{x}\!w\,\delta_{ij}.
\end{equation}
Section \ref{SSdecomposition} provides the details. The analysis is based on the relaxation terms $\rho n(w-u)$ and $3\rho n(T-\theta)$ in \eqref{twophase}, along with the effective stress and heat flux.
\begin{equation}\label{effcoeff}
\dv_{x}\Big(\frac{nT}{2\rho}\,\mathbb{S}(w)\Big),
\qquad
\dv_{x}\Big(\frac{5nT}{6\rho}\nabla_{x}T\Big).
\end{equation}

The end states $(\rho_{\pm},u_{\pm},\theta_{\pm})$ are connected by a rarefaction wave of the Riemann problem for the fluid Euler system
$\eqref{twophase}_{1,2,3}$,
\begin{equation}\label{Euler-fluid}
\begin{cases}
\rho_{t}+\dv_{x}(\rho u)=0,\\[2pt]
(\rho u)_{t}+\dv_{x}(\rho u\otimes u)+\nabla_{x}(\rho\theta)=0,\\[2pt]
\Big[\rho\Big(\frac32\theta+\frac{|u|^{2}}{2}\Big)\Big]_{t}
 +\dv_{x}\Big[\rho u\Big(\frac32\theta+\frac{|u|^{2}}{2}\Big)+\rho\theta u\Big]=0,
\end{cases}
\qquad x\in\R\times\T^{2},\ t>0,
\end{equation}
with Riemann data
\begin{equation}\label{Riemann-fluid}
(\rho,u,\theta)(0,x)=
\begin{cases}
(\rho_{-},u_{-},\theta_{-}),& x_{1}<0,\\
(\rho_{+},u_{+},\theta_{+}),& x_{1}>0,
\end{cases}
\qquad (x_{2},x_{3})\in\T^{2},
\end{equation}
and, likewise, $(n_{\pm},u_{\pm},\theta_{\pm})$ are connected by a rarefaction
wave of the particle Euler system $\eqref{twophase}_{4,5,6}$,
\begin{equation}\label{Euler-part}
\begin{cases}
n_{t}+\dv_{x}(nw)=0,\\[2pt]
(nw)_{t}+\dv_{x}(nw\otimes w)+\nabla_{x}(nT)=0,\\[2pt]
\Big[n\Big(\frac32T+\frac{|w|^{2}}{2}\Big)\Big]_{t}
 +\dv_{x}\Big[nw\Big(\frac32T+\frac{|w|^{2}}{2}\Big)+nTw\Big]=0,
\end{cases}
\qquad x\in\R\times\T^{2},\ t>0,
\end{equation}
with Riemann data
\begin{equation}\label{Riemann-part}
(n,w,T)(0,x)=
\begin{cases}
(n_{-},u_{-},\theta_{-}),& x_{1}<0,\\
(n_{+},u_{+},\theta_{+}),& x_{1}>0,
\end{cases}
\qquad (x_{2},x_{3})\in\T^{2}.
\end{equation}
The construction of the background profile is determined by the coupling. The two Riemann problems are solved consistently, in the sense that the two $3$-rarefaction fans are selected so that their velocities and temperatures coincide. This causes the relaxation terms in \eqref{twophase} to vanish on the background profile. The densities then satisfy $\rho^{r}/\rho_{+}=n^{r}/n_{+}$ and may differ if $\rho_{+}\neq n_{+}$ (see Section \ref{SScons-rare}). Since the self-similar rarefaction wave is only Lipschitz continuous, we use the smooth approximations $(\bar\rho,\bar u,\bar\theta)(t,x_{1})$ and $(\bar n,\bar u,\bar\theta)(t,x_{1})$, obtained as in \cite{MatsumuraNishihara1992} from a smoothed inviscid Burgers profile (see \eqref{dbur} and \eqref{au}). These profiles solve the one-dimensional Euler system:
\begin{equation}\label{1DEuler-fluid}
\begin{cases}
\bar\rho_{t}+(\bar\rho\bar u_{1})_{x_{1}}=0,\\[2pt]
(\bar\rho\bar u_{1})_{t}+(\bar\rho\bar u_{1}^{2}+\bar\rho\bar\theta)_{x_{1}}=0,
\qquad
(\bar\rho\bar u_{i})_{t}+(\bar\rho\bar u_{1}\bar u_{i})_{x_{1}}=0,\ i=2,3,\\[2pt]
\Big[\bar\rho\Big(\frac32\bar\theta+\frac{|\bar u|^{2}}{2}\Big)\Big]_{t}
 +\Big[\bar\rho\bar u_{1}\Big(\frac32\bar\theta+\frac{|\bar u|^{2}}{2}\Big)
   +\bar\rho\bar\theta\bar u_{1}\Big]_{x_{1}}=0,
\end{cases}
\end{equation}
and
\begin{equation}\label{1DEuler-part}
\begin{cases}
\bar n_{t}+(\bar n\bar u_{1})_{x_{1}}=0,\\[2pt]
(\bar n\bar u_{1})_{t}+(\bar n\bar u_{1}^{2}+\bar n\bar\theta)_{x_{1}}=0,
\qquad
(\bar n\bar u_{i})_{t}+(\bar n\bar u_{1}\bar u_{i})_{x_{1}}=0,\ i=2,3,\\[2pt]
\Big[\bar n\Big(\frac32\bar\theta+\frac{|\bar u|^{2}}{2}\Big)\Big]_{t}
 +\Big[\bar n\bar u_{1}\Big(\frac32\bar\theta+\frac{|\bar u|^{2}}{2}\Big)
   +\bar n\bar\theta\bar u_{1}\Big]_{x_{1}}=0,
\end{cases}
\end{equation}
respectively, and hence also the three-dimensional systems \eqref{Euler-fluid} and \eqref{Euler-part}.

We use the following notation for function spaces. For $s\in\mathbb{N}$, let $H^{s}_{x}(\Omega)$ be the usual Sobolev space on $\Omega=\R\times\T^{2}$ with
norm $\|f\|_{H^{s}_{x}}=\sum_{j=0}^{s}\|\nabla^{j}_{x}f\|$, where $\|\cdot\|:=\|\cdot\|_{L^{2}(\Omega)}$, and write
$\|(f,g)\|_{H^{s}}=\|f\|_{H^{s}}+\|g\|_{H^{s}}$. 
Let $\Mstar=\Mloc_{[n_{*},w_{*},T_{*}]}(v)$ be a global Maxwellian with constants$n_{*}>0$, $w_{*}\in\R^{3}$, $T_{*}>0$; the weighted space
$H^{s}_{x}\big(L^{2}_{v}(1/\sqrt{\Mstar})\big)$ on $\Omega\times\R^{3}$ carries
the norm
\begin{equation}\label{weightednorm}
\|f\|_{H^{s}_{x}(L^{2}_{v}(1/\sqrt{\Mstar}))}
=\sum_{j=0}^{s}\Big(\intv\frac{1}{\Mstar}
 \big\|\nabla^{j}_{x}f\big\|^{2}\,dv\Big)^{\frac12}.
\end{equation}

Our main result reads as follows.

\begin{theorem}\label{thm:main}
  Let $(\rho_\pm,n_\pm,u_\pm,\theta_\pm)$ be the end state constants satisfying \eqref{endstate},  $(\bar\rho,\bar u,\bar\theta)(t,x_{1})$ and
$(\bar n,\bar u,\bar\theta)(t,x_{1})$ be the smooth planar $3$-rarefaction waves of \eqref{1DEuler-fluid}-\eqref{1DEuler-part} given by \eqref{au}.
Set the wave strength $\varepsilon:=|(\r_+-\r_-,n_+-n_-,u_+-u_-, \t_+-\t_-)|$, and   the initial perturbation
\begin{equation}\label{E0}
\begin{split}
\mathcal{E}(0):=\ &\big\|\big(\rho_{0}-\bar\rho(0,\cdot),\,
   u_{0}-\bar u(0,\cdot),\,\theta_{0}-\bar\theta(0,\cdot)\big)\big\|^{2}_{H^{3}_{x}}\\
&+\big\|f_{0}(x,v)
   -\Mloc_{[\bar n(0,x_{1}),\bar u(0,x_{1}),\bar\theta(0,x_{1})]}
   \big\|^{2}_{H^{3}_{x}(L^{2}_{v}(1/\sqrt{\Mstar}))},
\end{split}
\end{equation}
where $(\rho_{0},u_{0},\theta_{0},f_{0})$ is the initial data \eqref{initial}  with $f_{0}\geq 0$, and $\Mstar=\Mloc_{[n_{*},w_{*},T_{*}]}(v)$ is a global
Maxwellian of the form \eqref{localM} with constants $n_{*}>0$, $w_{*}\in\R^{3}$, $T_{*}>0$ determined in \eqref{Mstar-choice}. 
There exists $\eps_{0}\in(0,1)$ such that if
\begin{equation}\label{smallness}
\mathcal{E}(0)+\eps\le\eps_{0},
\end{equation}
  the initial-value problem \eqref{EVFP} and \eqref{initial} for the three-dimensional full compressible Euler-VFP system admits a unique global strong solution $(\rho,u,\theta,f)$ such that
$$
\rho-\bar\rho,\ u-\bar u,\ \theta-\bar\theta\in C\big([0,\infty);H^3_x(\Omega)\big),
\qquad
f-\Mloc_{[\bar n,\bar u,\bar\theta]}\in C\big([0,\infty);H^3_x(L^2_v(1/\sqrt{\Mstar}))\big),
$$
and $f(t,x,v)\geq0$ for all $t\geq0$, satisfying
\begin{equation}\label{uniform}
\begin{split}
\sup_{t\ge0}\big\|\big(\rho-\bar\rho,\,u-\bar u,\,\theta-\bar\theta\big)(t)
  \big\|^{2}_{H^{3}_{x}}
&\le C_{0}\big(\mathcal{E}(0)+\eps^{\frac18}\big),\\
\sup_{t\ge0}\big\|f(t,x,v)
  -\Mloc_{[\bar n(t,x_{1}),\bar u(t,x_{1}),\bar\theta(t,x_{1})]}
  \big\|^{2}_{H^{3}_{x}(L^{2}_{v}(1/\sqrt{\Mstar}))}
&\le C_{0}\big(\mathcal{E}(0)+\eps^{\frac18}\big),
\end{split}
\end{equation}
for a constant $C_{0}>0$ independent of $t$. Moreover, the planar
$3$-rarefaction wave is time-asymptotically stable:
\begin{equation}\label{asymptotic}
\begin{split}
\lim_{t\to\infty}\big\|\big(\rho-\bar\rho,\,u-\bar u,\,
  \theta-\bar\theta\big)(t)\big\|_{L^{\infty}_{x}}&=0,\\
\lim_{t\to\infty}\big\|f(t,x,v)
  -\Mloc_{[\bar n(t,x_{1}),\bar u(t,x_{1}),\bar\theta(t,x_{1})]}
  \big\|_{L^{\infty}_{x}(L^{2}_{v}(1/\sqrt{\Mstar}))}&=0.
\end{split}
\end{equation}
\end{theorem}

\begin{remark}
Several remarks are as follows.

\noindent
 (1) The reference state in \eqref{E0}-\eqref{asymptotic} is a local Maxwellian associated with a rarefaction profile, not a global equilibrium. Consequently, the basis functions, the projections $\mathbf{P}_{0}$ and $\mathbf{P}_{1}$, and the linearized operator $\Lop_{[w,T]}$ depend on $(t,x)$ through the solution itself, so the spectral information available in the constant-coefficient case cannot be used.
%\end{remark}
%
%\begin{remark}

\noindent
(2) Both phases converge to the same velocity and temperature profile $(\bar u,\bar\theta)$, but their densities converge to two distinct profiles, $\bar\rho$ and $\bar n$. The first results from double relaxation, and the second from the absence of relaxation in the mass equations.
%\end{remark}
%
%\begin{remark}

\noindent
(3) The parameter $\eps$ is the strength of the rarefaction wave and the smoothing scale of the approximate profile \eqref{dbur} is normalized to one, so that all the smallness in the background is carried by $\eps$ alone (see Remark \ref{rk:eps}). The rate $\eps^{1/8}$ in
\eqref{uniform} comes from the non-integrable decay
$\|\bar u_{1x_{1}}\|^{2}_{L^{2}_{x_{1}}}\sim(1+t)^{-1}$ and from the correction
function $\bar{\mathcal{G}}$ introduced in \eqref{Gbar-intro} below to absorb
it.
%\end{remark}
%
%\begin{remark}

\noindent
(4)  We expect the framework developed herein to extend to other nonlinear waves in non-isentropic kinetic-fluid systems, particularly viscous contact waves and shock profiles. These cases, however, require additional arguments. For example, in the case of a contact wave, the background lacks a monotone profile that can be exploited.  %We leave them to future work.
\end{remark}

\subsection{Novelty and main difficulties}
\label{SSnovelty}

 We summarize the main contributions of the paper and then explain why the isentropic argument in \cite{LiWangWang2022} does not carry over.

 Theorem \ref{thm:main} provides, to our knowledge, the first stability result for a nonlinear wave in a fluid-particle model with energy exchange, and the first in which the reference state of such a model is a spatially varying local Maxwellian rather than a global equilibrium. The global results quoted in Section \ref{SSprevious} for the NS-VFP system
  %\eqref{EVFP} and \eqref{NSVFP} 
 all concern perturbations of a constant state.

The Chapman-Enskog expansion of $\eqref{EVFP}_{4}$ about the local Maxwellian \eqref{localM} yields the nonconservative viscous system $\eqref{twophase}_{4,5,6}$, which appears to be new for the full Euler-VFP system. Its viscosity $nT/2\rho$ and heat conductivity $5nT/6\rho$ depend on both the density and the temperature, whereas the classical Chapman-Enskog expansion of the Boltzmann equation yields coefficients that depend on the temperature alone. These coefficients are also not intrinsic to the particle phase, since they involve the fluid density $\rho$: the two phases are coupled through the dissipation, not only through the relaxation terms, and a loss of control over $\rho$ degrades the dissipation available for $(n,w,T)$. Transport coefficients that depend on additional state variables arise in the physics of high-temperature reacting mixtures \cite{ZeldovichRaizer} and in viscous radiative or reactive gases with internal energy production \cite{Feireisl}; see also \cite{Ducomet,LiaoZhao,UmeharaTani}.

 The null space of the Fokker-Planck operator $\dv_{v}\big((v-u)(\cdot)+\theta\nabla_{v}(\cdot)\big)$ is one-dimensional and spanned by the Maxwellian with fluid parameters $(\rho,u,\theta)$. It is therefore too small to capture relaxation through momentum and energy exchanges. We instead work in the five-dimensional space $\Nspace=\mathrm{span}\{\chi_{0},\chi_{1},\chi_{2},\chi_{3},\chi_{4}\}$ generated by the local Maxwellian \eqref{localM},
\begin{equation}\label{basis}
\chi_{0}=\frac{1}{\sqrt{n}}\,\Mloc,\qquad
\chi_{i}=\frac{v_{i}-w_{i}}{\sqrt{nT}}\,\Mloc\ \ (i=1,2,3),\qquad
\chi_{4}=\frac{1}{\sqrt{6n}}\Big(\frac{|v-w|^{2}}{T}-3\Big)\Mloc,
\end{equation}
 and establish a microscopic $H$-theorem that provides uniform coercivity of the linearized operator $\Lop_{[w,T]}$ on $\Nspace^{\perp}$ (Lemmas \ref{Borel}-\ref{co-1} and Corollary \ref{co-2}). The five macroscopic modes are removed simultaneously, and the resulting estimate controls the weighted velocity gradient and the velocity moments of $\mathbf{P}_{1}G$, including the heat-flux moments, which are absent in the isentropic theory.

 The fluid component of \eqref{EVFP} is purely hyperbolic. Our estimates show that the relative-velocity relaxation, the temperature relaxation, and the effective viscous and thermal dissipation \eqref{effcoeff} generated by the kinetic equation, together with the microscopic Fokker-Planck coercivity, provide sufficient regularization to stabilize a nonlinear wave. The resulting behavior is wave-like and differs from the diffusive behavior known for the compressible Euler equations with damping \cite{HsiaoLiu} and for the Fokker-Planck equation \cite{LinWangWu}.

\smallskip
 We next outline three differences between the present problem and the isentropic model \cite{LiWangWang2022}.

\smallskip
\noindent\emph{(1)  Enlarged macroscopic space.}  
%The energy exchange makes temperature an independent macroscopic variable, so the macroscopic space becomes five-dimensional (mass, momentum, energy) instead of four and includes the energy mode $\chi_4$. The enlargement itself causes no difficulty: by \eqref{eigen}, the operator $\Lop_{[w,T]}$ is strictly negative on $\mathrm{span}\{\chi_1,\dots,\chi_4\}$, so coercivity on $\mathfrak{N}^\bot$ is no harder than on $\{\chi_0\}^\bot$. What the enlargement requires is that the five moments, which define the hydrodynamic variables $(n,w,T)$ through \eqref{npart}, be removed simultaneously, so that the decomposition and the moment system \eqref{ns-1} are consistent. The difficulties lie elsewhere: by \eqref{basis}, the energy basis function depends on both $w$ and $T$, so that $\mathfrak{N}$, $\mb{P}_0$, $\mb{P}_1$ and $\Lop_{[w,T]}$ vary with the solution, and coercivity must then be transferred to the fixed weight $\Mstar$. This is what the microscopic $H$-theorem of \S\ref{SSH-Theorem} provides.
Energy exchange makes temperature an independent macroscopic variable, expanding the macroscopic space from four to five dimensions (mass, momentum, energy) and adding the energy mode $\chi_4$. This creates no difficulty: by \eqref{eigen}, $\Lop_{[w,T]}$ is strictly negative on $\mathrm{span}\{\chi_1,\dots,\chi_4\}$, so coercivity on $\mathfrak{N}^\bot$ is no harder than on $\{\chi_0\}^\bot$. However, all five moments defining hydrodynamic variables $(n,w,T)$ via \eqref{npart} must be removed simultaneously to keep the decomposition and moment system \eqref{ns-1} consistent. The main challenge is that, by \eqref{basis}, the energy basis function depends on $w$ and $T$, making $\mathfrak{N}$, $\mb{P}_0$, $\mb{P}_1$, and $\Lop_{[w,T]}$ solution-dependent. Coercivity must therefore be transferred to the fixed weight $\Mstar$, as provided by the microscopic $H$-theorem in Section \ref{SSH-Theorem}.

\smallskip
\noindent\emph{(2) Double relaxation.}  
%The isentropic model relaxes only the relative velocity $w-u$, whereas here the thermal relaxation $T-\theta$ is present as well and the two are coupled: by \eqref{exchange} the energy exchange contains the work term $\rho n\,w\cdot(w-u)$, which has no definite sign, so that the momentum and energy equations cannot be treated separately. Adding the two momentum equations of \eqref{twophase} cancels the drag force, while a suitable difference produces a coercive term for $w-u$; adding the two energy equations preserves the total energy, while a suitable combination of the temperature equations produces dissipation for $T-\theta$. The work terms are controlled by cancellations between the momentum and energy identities and by the smallness of the rarefaction profile. This makes the basic energy estimate (Proposition \ref{Prop3.1}) more involved than its isentropic counterpart.
%
Unlike the isentropic model, which relaxes only the relative velocity $w-u$, this model also includes coupled thermal relaxation $T-\theta$. By \eqref{exchange}, energy exchange includes the sign-indefinite work term $\rho n\,w\cdot(w-u)$, preventing separate treatment of the momentum and energy equations. Adding the momentum equations in \eqref{twophase} cancels drag, while a suitable difference yields coercivity in $w-u$. Adding the energy equations preserves total energy, and a suitable combination of temperature equations dissipates $T-\theta$. Cancellations between momentum and energy identities, together with the smallness of the rarefaction profile, control the work terms. Thus, the basic energy estimate (Proposition \ref{Prop3.1}) is more involved than its isentropic counterpart.

\smallskip
\noindent\emph{(3)  Effective heat flux and density-temperature-dependent coefficients.} 
%The Chapman-Enskog expansion yields both an effective viscous stress, as in the isentropic case, and an effective heat flux proportional to $\nabla_{x}T$. Since the transport coefficients in \eqref{effcoeff} depend on $n$, $T$ and $\rho$, differentiating them generates nonlinear terms coupling the density and temperature derivatives of both phases, and these terms have to be absorbed by the dissipation that they perturb.
%
The Chapman-Enskog expansion yields effective viscous stress, as in the isentropic case, and heat flux proportional to $\nabla_{x}T$. Because the transport coefficients in \eqref{effcoeff} depend on $n$, $T$, and $\rho$, differentiating them produces nonlinear terms coupling density and temperature derivatives across both phases. These terms must be absorbed by the dissipation they perturb.

\subsection{Strategy of the proof}
\label{SSstrategy}

The proof of Theorem \ref{thm:main} combines local existence results for strong solutions with uniform-in-time a priori estimates. We outline the main steps.

\smallskip
\noindent\emph{(1) The macro-micro system.} 
%Let $\mathbf{P}_{0}$ be the orthogonal projection onto the macroscopic space $\Nspace$ generated by \eqref{basis}, and let $\mathbf{P}_{1}=\mathbf{I}-\mathbf{P}_{0}$. Writing $f=\Mloc+G$ with $\mathbf{P}_{0}f=\Mloc$ and $\mathbf{P}_{1}f=G$, the microscopic component satisfies
Let $\mathbf{P}_{0}$ project orthogonally onto the macroscopic space $\Nspace$ generated by \eqref{basis}, and set $\mathbf{P}_{1}=\mathbf{I}-\mathbf{P}_{0}$. Write $f=\Mloc+G$, where $\mathbf{P}_{0}f=\Mloc$ and $\mathbf{P}_{1}f=G$. The microscopic component satisfies
\begin{equation}\label{micro}
G_{t}+\mathbf{P}_{1}(v\cdot\nabla_{x}G)-\rho\Lop_{[w,T]}G
=-\mathbf{P}_{1}(v\cdot\nabla_{x}\Mloc)
 +\rho(w-u)\cdot\nabla_{v}G-\rho(T-\theta)\Delta_{v}G .
\end{equation}
%Applying $\Lop_{[w,T]}^{-1}$ to $\mathbf{P}_{1}(v\cdot\nabla_{x}\Mloc)$ gives the leading Chapman-Enskog term, whose momentum moment produces the effective stress $\mathbb{S}(w)$ and whose energy moment produces the heat flux proportional to $\nabla_{x}T$; this is how \eqref{twophase} arises. The remainder is kept in the microscopic equation and absorbed by the Fokker-Planck dissipation (Section \ref{SSdecomposition}). Since $\|\bar u_{1x_{1}}\|^{2}_{L^{2}_{x_{1}}}\sim(1+t)^{-1}$ is not integrable in time, the microscopic perturbation is measured relative to the correction
%
Applying $\Lop_{[w,T]}^{-1}$ to $\mathbf{P}_{1}(v\cdot\nabla_{x}\Mloc)$ yields the leading Chapman-Enskog term: its momentum moment gives the effective stress $\mathbb{S}(w)$, and its energy moment gives heat flux proportional to $\nabla_{x}T$, producing \eqref{twophase}. The remainder stays in the microscopic equation and is absorbed by Fokker-Planck dissipation (Section \ref{SSdecomposition}). Because $\|\bar u_{1x_{1}}\|^{2}_{L^{2}_{x_{1}}}\sim(1+t)^{-1}$ is not time-integrable, the microscopic perturbation is measured relative to the correction
\begin{equation}\label{Gbar-intro}
\bar{\mathcal{G}}(t,x,v)=\frac{1}{\rho T}\,\Lop_{[w,T]}^{-1}
\Big\{\mathbf{P}_{1}\Big[{v_{1}}\Big(\frac{|v-w|^{2}}{2T}\bar\theta_{x_{1}}
 +{v_{1}}\bar u_{1x_{1}}\Big)\Mloc_{[n,w,T]}\Big]\Big\}.
\end{equation}

\smallskip
\noindent\emph{(2) The microscopic $H$-theorem.} 
%{Both the operator and the macroscopic projection vary with the solution, and the energy basis function depends on $w$ and $T$. We show that the quadratic form of $\Lop_{[w,T]}$ controls the weighted velocity gradient and the velocity moments of $\mathbf{P}_{1}G$, orthogonality to all five macroscopic variables being used here, and we then transfer the estimate to the fixed weighted space $L^{2}_{v}(1/\sqrt{\Mstar})$ used in the nonlinear energy method (Lemma \ref{co-1} and Corollary \ref{co-2}).}
Both the operator and macroscopic projection depend on the solution, while the energy basis function depends on $w$ and $T$. We show that the quadratic form of $\Lop_{[w,T]}$ controls the weighted velocity gradient and $\mathbf{P}_{1}G$ velocity moments, assuming orthogonality to all five macroscopic variables. We then transfer the estimate to the fixed weighted space 
$L^{2}_{v}(1/\sqrt{\Mstar})$ used in the nonlinear energy method (Lemma \ref{co-1} and Corollary \ref{co-2}).

\smallskip
\noindent\emph{(3) The double relaxation.} 
%As described in \S\ref{SSnovelty}, sums and differences of the two momentum equations of \eqref{twophase} yield a coercive term for $w-u$, and the corresponding combination of the energy and temperature equations yields dissipation for $T-\theta$; the indefinite work terms $\rho n\,w\cdot(w-u)$ are controlled by cancellation between the two identities. This is carried out in Proposition \ref{Prop3.1}.
%
As described in Section \ref{SSnovelty}, combining the momentum equations in \eqref{twophase} yields coercivity for $w-u$, and combining the energy and temperature equations yields dissipation for $T-\theta$. Cancellation between the identities controls the indefinite work terms $\rho n\,w\cdot(w-u)$ (Proposition \ref{Prop3.1}).

\smallskip
\noindent\emph{(4) Higher-order spatial dissipation.} 
%The fluid component has neither viscosity nor heat conduction. The effective stress and heat flux in $\eqref{twophase}_{5,6}$ nevertheless control $\nabla_{x}w$ and $\nabla_{x}T$, and the double relaxation estimates control $w-u$ and $T-\theta$ together with their derivatives, hence transfer that control to $u$ and $\theta$. The densities are handled by means of the hyperbolic structure, the differentiated continuity equations and the monotonicity of the $3$-rarefaction profile. The kinetic phase thus supplies the regularizing mechanism that the inviscid fluid phase lacks (Propositions \ref{Prop3.2} and \ref{Prop3.3}).
%
Although the fluid has no viscosity or heat conduction, the effective stress and heat flux in $\eqref{twophase}_{5,6}$ control $\nabla_{x}w$ and $\nabla_{x}T$. Double relaxation estimates control $w-u$ and $T-\theta$ and their derivatives, thereby extending this control to $u$ and $\theta$. Hyperbolic structure, differentiated continuity equations, and monotonicity of the $3$-rarefaction profile control the densities. Thus, the kinetic phase regularizes the inviscid fluid phase (Propositions \ref{Prop3.2} and \ref{Prop3.3}).

\smallskip
\noindent\emph{(5) Closing the estimates.} 
%The lower-order fluid estimates, the higher-order macroscopic estimates and the microscopic estimates are combined into a single closed a priori inequality. The resulting uniform bound continues the local solution globally, and the time-integrability of the dissipative quantities, together with the bounds on their time derivatives, gives the convergence \eqref{asymptotic} to the planar rarefaction wave.
%
Combining the lower-order fluid, higher-order macroscopic, and microscopic estimates yields a closed a priori inequality. Its uniform bound extends the local solution globally, and the dissipative quantities' time-integrability and bounded time derivatives imply convergence \eqref{asymptotic} to the planar rarefaction wave.

%\subsection{Notation and outline}
%\label{SSnotation}
\medskip

%\noindent 
{\bf Notation.} Positive generic constants that are uniform in time are all denoted by $C$;
$A\sim B$ means that $C^{-1}A\le B\le CA$ for some universal constant $C$. We
write $\partial^{\alpha}$ for either $\partial_{t}$ or $\nabla_{x}$, with
$\nabla_{x}=(\partial_{x_{1}},\partial_{x_{2}},\partial_{x_{3}})$, and set
$\partial_{i}:=\partial_{x_{i}}$ $(i=1,2,3)$, $\nabla:=\nabla_{x}$,
$\dv:=\dv_{x}$, $\|\cdot\|_{H^{s}}:=\|\cdot\|_{H^{s}_{x}}$, where $\dv$ and
$\nabla$ act on the spatial variable $x\in\Omega$. Finally
\[
\int f\,dv:=\int_{\R^{3}}f\,dv,
\qquad
\iint f\,dvdx:=\int_{\R^{3}}\!\!\int_{\Omega}f\,dvdx .
\]

The paper is organized as follows. In Section \ref{SPreliminaries} we set up the macro-micro
decomposition around the local Maxwellian, derive the effective nonconservative
system \eqref{twophase} through the Chapman-Enskog expansion, prove the
microscopic $H$-theorem, and collect the properties of the smooth planar
$3$-rarefaction wave together with a Gagliardo-Nirenberg inequality on
$\Omega=\R\times\T^{2}$. Section 3 reduces Theorem \ref{thm:main} to the
uniform a priori estimates of Theorem \ref{priori}, which are stated there as
Propositions \ref{Prop3.1}-\ref{Prop3.3} and proved in Section \ref{Slow-es}. Two technical proofs are deferred to Appendices \ref{Appen-1} and \ref{Appen-2}.%%%%%%%%%%%%%%%%%%%%%%%%%%%%%%%%%%%%%%%%%%%%%%%%%%%%%%%%%%%%%%%%%%%%%%%%%%%%%%

\section{Preliminaries}
\label{SPreliminaries}

In this section, we introduce the macro-micro decomposition around a local Maxwellian for the Vlasov--Fokker--Planck equation. This provides a unified framework for studying the time-asymptotic stability of basic wave patterns for the multi-dimensional full compressible Euler--VFP system \eqref{EVFP}. In particular, a new nonconservative viscous compressible fluid model is derived via the Chapman--Enskog expansion of the Vlasov--Fokker--Planck equation, owing to a suitable choice of the local Maxwellian. Its transport coefficients depend on both temperature and density. We then establish a microscopic $H$-theorem for the Vlasov--Fokker--Planck equation around the local Maxwellian, and recall the existence and time-decay properties of the planar rarefaction wave. Finally, we state a Gagliardo--Nirenberg inequality on the domain  $\Omega:=\mathbb{R}\times \mathbb{T}^2$.

\subsection{Macro-micro decomposition}
\label{SSdecomposition}
We consider the macro--micro decomposition of solutions to the Vlasov--Fokker--Planck equation \eqref{EVFP}$_4$ into its macroscopic (fluid dynamical) and microscopic (kinetic) components. The null space of the original Fokker--Planck operator
$$\Lop(\cdot)=\dv_v\big((v-u)(\cdot)+\t\na_v(\cdot) \big)
$$
is the one dimensional space
$$
{\rm span}\Big\{ \f{\r(t,x)}{\sqrt{ (2\pi \t(t,x))^3}} \exp{\Big(-\f{|v-u(t,x)|^2}{2\t(t,x)}\Big)} \Big\},
$$
which is insufficient to capture the desired relaxation friction damping in the momentum and energy exchanges. Therefore, motivated by the recent work \cite{LiWangWang2022}, we define a new space $\mathfrak{N}$ using the specific local Maxwellian \eqref{localM}, and introduce macroscopic and microscopic projections onto this space, as described below.

Let $\mathcal{M}$ be the local Maxwellian defined in \eqref{localM}, and define the linear space $\mathfrak{N}$ by
$$
\mathfrak{N}={\rm span}\{\chi_0,\chi_1,\chi_2,\chi_3,\chi_4\},
$$
where $\chi_i~(i=0,\cdots,4)$ are pairwise orthogonal and given by
\begin{equation*}
 \chi_0(v) = {\f1{\sqrt{n}}\Mloc}, \quad
 \chi_i(v) ={\f{v_i-w_i}{\sqrt{nT}}\Mloc}~~ (i=1,2,3), \quad
 \chi_4(v) = {\f{1}{\sqrt{6n}}\Big(\f{|v-w|^2}{ T}-3\Big)\Mloc}
\end{equation*}
with
$$
\langle \chi_i(v),\chi_j(v)\rangle:=\int_{\mathbb{R}^3}\f{1}{\Mloc}\chi_i(v)\chi_j(v) dv=\delta_{ij}, \qquad i,j=0,\cdots,4.
$$
In terms of this orthogonal basis, the macroscopic projection $\mb{P}_0$ onto
$\mathfrak{N}$ and the microscopic projection $\mb{P}_1$ onto $\mathfrak{N}^\bot$ are defined by
\begin{equation}\label{p0}
 \mb{P}_0g = {\di\sum_{j=0}^4\langle g,\chi_j\rangle\chi_j},\qquad
 \mb{P}_1g= g-\mb{P}_0g.
\end{equation}

We decompose the solution $f(t,x,v)$ of the Vlasov--Fokker--Planck equation \eqref{EVFP}$_4$ with respect to the local Maxwellian \eqref{localM} as
$$
f(t,x,v)=\Mloc(t,x,v)+G(t,x,v),
$$
where $\mb{P}_0f=\Mloc$ is called the macroscopic component of $f$, $\mb{P}_1f=G$ is called the microscopic component of $f$.

For given fluid velocity $u(x,t)$ and temperature $\t(t,x)$, equation $\eqref{EVFP}_4$ can be rewritten as
\begin{equation}\label{VFP-L}
f_t+v\cdot\na_x f=\r \mathcal{L}_{[w,T]}f+\r \dv_v\big((w-u)f-(T-\t)\na_v f\big),
\end{equation}
where $w=(w_1,w_2,w_3)$ is the velocity and $T$ is the temperature associated with $f(t,x,v)$, defined in \eqref{npart}. The linearized Fokker--Planck operator $\mathcal{L}_{[w,T]}$ is given by
\begin{equation}\label{L}
\mathcal{L}_{[w,T]}f=\dv_v\big((v-w)f+T\na_vf\big).
\end{equation}
It is straightforward to verify that $\Lop_{[w,T]}\Mloc=0$.  With the decomposition $f=\Mloc+G$, equation $\eqref{VFP-L}$ becomes
\begin{equation}\label{M+G}
(\Mloc+G)_t+v\cdot\nabla_x(\Mloc+G)=\r \Lop_{[w,T]}G+\r \dv_v\big((w-u)(\Mloc+G)-(T-\t)\na_v(\Mloc+G)\big).
\end{equation}
Multiplying the equation \eqref{M+G} by $1$, $v_i$ $(i=1,2,3)$, $\f12|v|^2$ respectively, and integrating by parts with respect to $v$ over ${\mathbb R}^3$, we obtain the following fluid-dynamical system for $(n, w, T)$:
\begin{equation}\label{ns-1}
\begin{cases}
 n_{t}+ \dv_x (nw )=0, \\
\displaystyle (n w)_t+\dv_x (n w\otimes w) +\nabla_x (nT) =-\r n(w-u)-\int_{\mathbb{R}^3} v\otimes v\cdot\nabla_xG\,dv,\\
\displaystyle \Big[n\Big(\f32T+\f{|w|^2}{2}\Big)\Big]_t+\dv_x\Big[nw\Big(\f32T+\f{|w|^2}{2}\Big)+wnT\Big]\\
=-3\r n(T-\t)-\r nw\cdot(w-u)-\int\f{1}{2}|v|^2v\cdot\na_xG\,dv.
\end{cases}
\end{equation}
The fluid--dynamical system \eqref{ns-1} for $(n,w,T)$ is not closed,  and an additional equation for the microscopic component $G$ is needed. This is obtained
by applying the projection operator $\mb{P}_1$ to \eqref{M+G}:
\begin{equation}\label{non-fluid}
G_t+\mb{P}_1(v\cdot\nabla_x\Mloc)+\mb{P}_1(v\cdot\nabla_x G)
=\r \Lop_{[w,T]}G+\r (w-u)\cdot\na_v G-\r (T-\t)\Delta_v G,
\end{equation}
where we have used the identities
$$
\langle \dv_v\big((w-u)f-(T-\t)\na_v f\big),\chi_0\rangle=0
$$
and
\begin{equation*}
\begin{split}
&\quad\,\,\mb{P}_1\big[\dv_v\big((w-u)f-(T-\t)\nabla_v f\big)\big]\\
&=\dv_v\big((w-u)f-(T-\t)\na_v f\big)-\mb{P}_0\big[\dv_v\big((w-u)f-(T-\t)\na_v f\big)\big]\\
&=(w-u)\cdot\na_vf-(T-\t)\Delta_v f-(w-u)\cdot\na_v\Mloc+(T-\t)\Delta_v\Mloc\\
&=(w-u)\cdot\na_vG-(T-\t)\Delta_vG.
\end{split}
\end{equation*}

The linearized operator $\Lop_{[w,T]}$ is coercive and invertible
on $\mathfrak{N}^\bot$ (see Lemmas \ref{le-diss} and \ref{co-1}).
This property, together with \eqref{non-fluid}, implies the following expression for $G$:
\begin{equation}\label{Pi}
\begin{aligned}
G &=\r^{-1}\Lop_{[w,T]}^{-1}[\mb{P}_1(v\cdot\nabla_x\Mloc)]+\Lop_{[w,T]}^{-1}[\r^{-1}G_t+\r^{-1}\mb {P}_1(v\cdot\nabla_xG)
-(w-u)\cdot\na_vG+(T-\t)\Delta_vG]\\
 :&= \r^{-1}\Lop_{[w,T]}^{-1}[\mb{P}_1(v\cdot\nabla_x\Mloc)] +\Pi.
\end{aligned}
\end{equation}
Substituting \eqref{Pi} into \eqref{ns-1}, we obtain the following full compressible fluid--dynamic system for $(n,w,T)$:
\begin{equation}\label{ns-2}
\left\{
\begin{array}{l}
\di n_{t}+{\rm div}_x (n w)=0,\\
\di (n w)_t+{\rm div}_x  (nw\otimes w) +\nabla_x(nT)=-\rho n(w-u)\\
\di \qquad\qquad -\int v\otimes v\cdot\nabla_x \big(\r^{-1}\Lop_{[w,T]}^{-1}[\mb{P}_1(v\cdot\nabla_x\Mloc)]\big) dv
-\int v\otimes v\cdot\nabla_x\Pi dv,\\
\di \Big[n\Big(\f32T+\f{|w|^2}{2}\Big)\Big]_t+\dv_x\Big[nw\Big(\f32T+\f{|w|^2}{2}\Big)+nTw\Big]=-3\rho n(T-\t)-\rho nw\cdot(w-u)\\
\di\qquad\qquad-\int \f{1}{2}|v|^2 v\cdot\nabla_x \big(\r^{-1}\Lop_{[w,T]}^{-1}[\mb{P}_1(v\cdot\nabla_x\Mloc)]\big) dv
-\int\f{1}{2}|v|^2v\cdot\na_x\Pi dv.
\end{array}
\right.
\end{equation}
Note that the two integral terms in \eqref{ns-2},
$$
-\int v\otimes v\cdot\nabla_x \big(\r^{-1}\Lop_{[w,T]}^{-1}[\mb{P}_1(v\cdot\nabla_x\Mloc)]\big) dv,\quad
-\int \f{1}{2}|v|^2 v\cdot\nabla_x \big(\r^{-1}\Lop_{[w,T]}^{-1}[\mb{P}_1(v\cdot\nabla_x\Mloc)]\big) dv
$$
correspond to the viscosity and heat conductivity, respectively.

The following lemma gives the explicit expressions for the viscosity and heat conductivity terms in \eqref{ns-2}.

\begin{lemma}\label{coefficient}
Let $\Mloc$ be the Maxwellian defined in \eqref{localM}. Then we have
\be\label{viscosity}
\di -\int v\otimes v\cdot\nabla_x \big(\r^{-1}\Lop_{[w,T]}^{-1}\big[\mb{P}_1(v\cdot\nabla_x\Mloc)\big]\big) dv
=\dv_x\big(\f{nT}{2\r}\, \mathbb{S}(w)\big),
\ee
\be\label{heat-conduct}
\di -\int \f{1}{2}|v|^2 v\cdot\na_x\big(\r^{-1}\Lop^{-1}_{[w,T]}\big[\mb{P}_1(v\cdot\nabla_x\Mloc)\big]\big) dv
=\dv_x\big(\f{5nT}{6\r}\,\na_x T\big)+\dv_x\big(\f{nT}{2\r} w\cdot\mathbb{S}(w)\big),
\ee
where
$$
\di \mathbb{S}(w):=2\mathbb{D}(w)-\f23\dv_x w\mathbb{I}, \quad (\mathbb{S}(w))_{ij}=\pa_{x_j}w_i+\pa_{x_i}w_j-\f23\dv_x w\,\d_{ij},~~i,j=1,2,3.
$$
\end{lemma}
The proof of this lemma is given in Appendix \ref{Appen-1}.
\smallskip
Substituting \eqref{viscosity} and \eqref{heat-conduct} into \eqref{ns-2} yields following system for $(n,w,T)$:
\be\label{mix-e}
\left\{
\begin{array}{l}
\di n_{t}+ \dv_x (n w)=0,\\
\di (n w)_t+ \dv_x  (nw\otimes w) +\nabla_x (nT)=\dv_x\Big(\f{nT}{2\r}\, \mathbb{S}(w)\Big)-\r n(w-u)-\int v\otimes v\cdot\nabla_x\Pi dv,\\
\di \Big[n\Big(\f32 T+\f{|w|^2}{2}\Big)\Big]_t+\dv_x\Big[n w \Big(\f32T+\f{|w|^2}{2}\Big)+n T w\Big]=\dv_x\Big(\f{5nT}{6\r}\,\na_x T\Big)\\
\di+\dv_x\Big(\f{nT}{2\r} w\cdot\mathbb{S}(w)\Big)-3\r n(T-\t)-\r nw\cdot(w-u)-\int\f12|v|^2v\cdot\na_x\Pi dv,
\end{array}
\right.
\ee
where $\Pi$ is defined in \eqref{Pi}, and $(u,\theta)$ are  the fluid velocity and temperature from the Euler--VFP system \eqref{EVFP}.

\smallskip
\begin{remark}
The fluid-dynamical model \eqref{mix-e} is derived here for the first time for the full compressible Euler--VFP system \eqref{EVFP}. A key feature of this model is that its transport coefficients depend on both density and temperature. This is notable in view of the classical Chapman--Enskog expansion for the Boltzmann equation, which yields the compressible Navier--Stokes equations with temperature-dependent transport coefficients. Indeed, while the viscous stress tensor $\mathbb{S}(w)$ in both cases takes the same form, the coefficients obtained here depend on the density as well as the temperature.
 The corresponding Prandtl number follows directly from $
\mu_{\rm eff}=nT/2\r$, $\kappa_{\rm eff}=5nT/6\r$, and $c_p=\f{\gamma}{\gamma-1}R=\f52$:
$$
{\rm Pr}=\f{c_p\mu_{\rm eff}}{\kappa_{\rm eff}}
=\f52\cdot\f{nT/2\r}{5nT/6\r}=\f32,
$$
compared with ${\rm Pr}\approx\f23$ for a monatomic gas governed by the Boltzmann equation. Thus, relative to viscous dissipation, the Fokker--Planck coupling produces weaker thermal dissipation than classical kinetic theory does.
\end{remark}

\begin{remark}
We note that the physics of high-temperature gases in reacting mixtures, as described, for instance, in the monograph \cite{ZeldovichRaizer}, indicates that both viscosity and thermal conductivity may depend on other state variables, such as density $n$ and temperature $T$. This observation is consistent with the discussion in \cite{Feireisl}, where internal energy production is considered, for instance, in viscous, radiative, or reactive gases. Further related studies can be found in \cite{Ducomet,LiaoZhao,UmeharaTani} and references therein.
\end{remark}

\subsection{Microscopic $H$-theorem}
\label{SSH-Theorem}

We list some dissipative properties of the linearized Fokker--Planck operator $\Lop_{[w,T]}$ defined in \eqref{L} in the weighted space $L^2_v$. To this end, we first present the spectral structure of $\Lop_{[w,T]}$, which gives rise to its dissipative properties.

\begin{lemma}\label{Borel}
Let $\Mloc$ be the local Maxwellian defined in \eqref{localM} and $\chi_0,\dots,\chi_4$ be the basis functions \eqref{basis}. Then
\be\label{eigen}
\Lop_{[w,T]}\chi_0=0,\qquad \Lop_{[w,T]}\chi_i=-\chi_i\ \ (i=1,2,3),\qquad
\Lop_{[w,T]}\chi_4=-2\chi_4,
\ee
so that $\ker\Lop_{[w,T]}=\mathrm{span}\{\Mloc\}$ is one dimensional, and moreover
\be\label{spec}
-\int\f{g\Lop_{[w,T]} g}{\Mloc}dv\geq 2\int\f{g^2}{\Mloc}dv
\qquad\text{for every }g\in \mathfrak{N}^\bot.
\ee
\end{lemma}

%\noindent {\it Proof.} 
\noindent {\it Proof.} 
For any polynomial $\varphi$, a direct computation from \eqref{L} gives the conjugation identity
\be\label{conj}
\Lop_{[w,T]}(\varphi\Mloc)=\big(T\Delta_v\varphi-(v-w)\cdot\na_v\varphi\big)\Mloc,
\ee
the contributions $3\varphi\Mloc$ and $\pm\varphi|v-w|^2\Mloc/T$ cancelling. Hence $\Lop_{[w,T]}$ is, up to the normalization of the variance, the
Ornstein--Uhlenbeck operator: if $\varphi$ is a Hermite polynomial of degree $k$ associated with $\Mloc$, then $T\Delta_v\varphi-(v-w)\cdot\na_v\varphi=-k\varphi$, so that $\varphi\Mloc$ is an eigenfunction of $\Lop_{[w,T]}$ with eigenvalue $-k$. Taking $\varphi=1$,
$\varphi=v_i-w_i$ and $\varphi=|v-w|^2/T-3$ in \eqref{conj} yields \eqref{eigen}; for the last one, $T\cdot\f{6}{T}-\f{2|v-w|^2}{T}=-2\big(\f{|v-w|^2}{T}-3\big)$.
Since $\Lop_{[w,T]}$ is self-adjoint and nonpositive for $\langle\cdot,\cdot\rangle$, the spectrum of $-\Lop_{[w,T]}$ on $L^2_v(1/\sqrt{\Mloc})$ is $\{0,1,2,\dots\}$.
By \eqref{eigen} the eigenspace of $-\Lop_{[w,T]}$ for the eigenvalue $1$ is exactly $\mathrm{span}\{\chi_1,\chi_2,\chi_3\}\subset\mathfrak{N}$, while $\mathfrak{N}$ also contains the eigenvalue-$0$ eigenfunction $\chi_0$ and the radial element $\chi_4$ of the eigenvalue-$2$ eigenspace. Consequently $\mathfrak{N}^\bot$ is spanned by eigenfunctions of $-\Lop_{[w,T]}$ with eigenvalues $k\geq 2$, which is \eqref{spec}. \hfill $\Box$

\medskip
We establish the dissipative property of the linearized Fokker--Planck operator $\Lop_{[w,T]}$ in \eqref{L} around the local Maxwellian $\Mloc$ defined in \eqref{localM} as follows.

\begin{lemma}\label{le-diss}
Let $\Mloc$ be the local Maxwellian defined in \eqref{localM}, and suppose $T_\star\leq T\leq T^\star$ for some positive constants $T_\star$, $T^\star$, then there exists a positive constant $\sigma(w,T_\star,T^{\star})$ depending on $w$, $T_\star$ and $T^{\star}$, such that for any $f$ for which  $\int\f{f\Lop_wf}{\Mloc}dv$ is well-defined, the following inequality holds:
\be\label{diss-1}
-\int\f{f\Lop_{[w,T]} f}{\Mloc} dv\geq \sigma(w,T_\star,T^\star)\int\f{|\na_v(\mb{P}_1f)|^2
+(1+|v|^2)|\mb{P}_1f|^2}{\Mloc}dv.
\ee
Consequently, for any $g\in \mathfrak{N}^\bot$,
\be\label{diss-21}
-\int\f{g\Lop_{[w,T]} g}{\Mloc} dv\geq \sigma(w,T_\star,T^\star)\int\f{|\na_v g|^2+(1+|v|^2)g^2}{\Mloc}dv.
\ee
\end{lemma}

The proof of this lemma is given in Appendix \ref{Appen-2}.

As a direct consequence of Lemma \ref{le-diss}, the following microscopic H-theorem is crucial for controlling the non-fluid  component $G$ in \eqref{non-fluid} in the stability analysis.
\begin{lemma}[$H$-theorem]\label{co-1}
There exist two positive constants $\bar \sigma =\bar \sigma (w_*,T_*)$ and $\eta_0=\eta_0(w_*,T_*)$ such that if
$|w-w_*|+|T-T_*|<\eta_0$, then for any $g\in \mathfrak{N}^\bot$
\be\label{diss-3}
-\int\f{g\Lop_{[w,T]} g}{\Mloc_*} dv\geq \bar\sigma \int\f{|\na_v g|^2+(1+|v|^2)g^2}{\Mloc_*}dv,
\ee
where $\Mloc_*=\f{n_*}{\sqrt{ (2\pi T_*)^3}} \exp{\big(-\f{|v-w_*|^2}{2T_*}\big)}$ is a global Maxwellian with constants $n_*>0$, $w_*\in\mathbb{R}^3$ and $T_*>0$.
\end{lemma}

%\noindent {\it Proof.}  
\noindent {\it Proof.} 
From the proof of Lemma \ref{le-diss}, it follows that \eqref{diss-1} and \eqref{diss-21} hold when $w$ and $T$ are constants. More precisely, there exists a positive constant $\sigma(w_*,T_*)$ depending only on $w_*$ and $T_*$, such that  for any $g\in \mathfrak{N}^\bot$, it holds
\be\label{diss-212}
-\int\f{g\Lop_{[w_*,T_*]} g}{\Mloc_*} dv\geq \sigma(w_*,T_*)\int\f{|\na_v g|^2+(1+|v|^2)g^2}{\Mloc_*}dv.
\ee
Using \eqref{diss-212} and Cauchy's inequality yields
$$
\ba
-\int\f{g\Lop_{[w,T]} g}{\Mloc_*} dv&=-\int\f{g\Lop_{[w_*,T_*]} g}{\Mloc_*} dv-\int\f{g(\Lop_{[w,T]} g
-\Lop_{[w_*,T_*]} g)}{\Mloc_*} dv\\
&=-\int\f{g\Lop_{[w_*,T_*]} g}{\Mloc_*} dv-\int\f{g}{\Mloc_*}(w_*-w)\cdot\na_v g dv
+\int (T-T_*)\na_v g\cdot\na_v\big(\f{g}{\Mloc_*}\big)dv\\
&=-\int\f{g\Lop_{[w_*,T_*]} g}{\Mloc_*} dv-\int\f{g}{\Mloc_*}(w_*-w)\cdot\na_v g dv\\
&\quad+(T-T_*)\int \Big(\f{|\na_v g|^2}{\Mloc_*}+\f{g}{\Mloc_*}\f{v-w_*}{T}\cdot\na_v g\Big) dv\\
& \geq \sigma(w_*,T_*)\int\f{|\na_v g|^2+(1+|v|^2)g^2}{\Mloc_*}dv\\
&\quad-C(w_*,T_*)(|w-w_*|+|T-T_*|)\int\f{|\na_v g|^2+(1+|v|^2)g^2}{\Mloc_*}dv,
\ea
$$
which leads to \eqref{diss-3} with the choice  $\eta_0=\f{\sigma(w_*,T_*)}{2C(w_*,T_*)}$. The proof of Lemma \ref{co-1} is complete.

\hfill $\Box$

\begin{remark}\label{rk:Mstar}
 The hypothesis $|w-w_*|+|T-T_*|<\eta_0$ in Lemma \ref{co-1} is not an additional restriction but rather a condition on the choice of the global Maxwellian $\Mstar$. Throughout, we fix $(n_*,w_*,T_*)$ so that
\be\label{Mstar-choice}
|w_*-u_\pm|+|T_*-\t_\pm|\leq C\eps,
\qquad \t_\pm<2T_*,
\ee
which is possible because the end states satisfy \eqref{endstate} and
$\eps$ is the wave strength. With this choice, \eqref{rho-bd}-\eqref{assum-1}
give
$$
|w-w_*|+|T-T_*|\leq |w-\bar u|+|T-\bar\t|+|\bar u-w_*|+|\bar\t-T_*|\leq C(\d+\eps)
$$
uniformly in $(t,x)$, so that the hypothesis of Lemma \ref{co-1} holds as soon as
$\eps_0$ is small. The second condition in \eqref{Mstar-choice} guarantees
$\int|\Mloc_{[n,w,T]}|^2\Mstar^{-1}dv<\infty$, so that the weighted norms
\eqref{weightednorm} used in \eqref{E0} and \eqref{uniform} are finite.
In the a priori estimates of Section \ref{Slow-es}, we therefore never choose
$\eta_0$, which is determined by $(w_*,T_*)$ through Lemma \ref{co-1}, but
only make $\d+\eps$ small. Every occurrence of $\eta_0$ below is to be read
as the quantity $C(\d+\eps)$ furnished by this remark.
\end{remark}

As a direct corollary of Lemma \ref{co-1} and H\"older inequality, we have
\begin{corollary}\label{co-2}
Under the assumptions of Lemma \ref{co-1}, for any $g\in \mathfrak{N}^\bot$
$$
\left\{
\ba
&\int\f{|\na_v(\Lop^{-1}_{[w,T]} g)|^2+(1+|v|^2)|\Lop^{-1}_{[w,T]} g|^2}{\Mloc} dv\leq \bar\sigma^{-2}
\int\f{(1+|v|^2)^{-1}}{\Mloc}|g|^2dv,\\
&\int\f{|\na_v(\Lop^{-1}_{[w,T]} g)|^2+(1+|v|^2)|\Lop^{-1}_{[w,T]} g|^2}{\Mloc_*} dv\leq \bar\sigma^{-2}
\int\f{(1+|v|^2)^{-1}}{\Mloc_*}|g|^2dv.
\ea
\right.
$$
\end{corollary}
\begin{remark}\label{rk:Nperp}
% Although the space $\mathfrak{N}$ depends on $(n,w,T)$, and hence on
%$(t,x)$ through the solution, its orthogonal complement does not. Indeed, since
%$\chi_0/\Mloc=n^{-1/2}$, $\chi_i/\Mloc=(v_i-w_i)(nT)^{-1/2}$ and
%$\chi_4/\Mloc=(|v-w|^2/T-3)(6n)^{-1/2}$, the conditions
%$\langle g,\chi_j\rangle=0$, $j=0,\dots,4$, are equivalent to
%\be\label{Nperp}
%\int g\,dv=0,\qquad \int g\,v\,dv=0,\qquad \int g\,|v|^2dv=0,
%\ee
%that is, to the vanishing of the collision invariants of $g$, which involves
%no parameter. This is what makes the varying-projection framework consistent;
%in particular $\mathfrak{N}^\bot$ coincides with the orthogonal complement of the
%macroscopic space associated with the global Maxwellian $\Mstar$, which is used
%tacitly in the proof of Lemma \ref{co-1}. A further consequence of \eqref{Nperp},
%used repeatedly below, is that $\mathfrak{N}^\bot$ is invariant under $\na_v$ and
%$\Delta_v$: for $g\in\mathfrak{N}^\bot$ one has
%$\langle\na_vg,\chi_0\rangle=0$,
%$\langle\pa_{v_k}g,\chi_i\rangle=-\d_{ik}(nT)^{-1/2}\int g\,dv=0$ and
%$\langle\pa_{v_k}g,\chi_4\rangle=-\f{2}{T\sqrt{6n}}\int g(v_k-w_k)dv=0$,
%and similarly for $\Delta_vg$. Hence $\Lop^{-1}_{[w,T]}$ may be applied to
%$\na_v(\Lop^{-1}_{[w,T]}h)$ and $\Delta_v(\Lop^{-1}_{[w,T]}h)$ in Lemma
%\ref{ex} below.
%
Although the space $\mathfrak{N}$ depends on $(n,w,T)$, and thus on $(t,x)$ through the solution, its orthogonal complement is parameter-independent. Since $\chi_0/\Mloc=n^{-1/2}$, $\chi_i/\Mloc=(v_i-w_i)(nT)^{-1/2}$, and $\chi_4/\Mloc=(|v-w|^2/T-3)(6n)^{-1/2}$, the conditions $\langle g,\chi_j\rangle=0$, $j=0,\dots,4$, are equivalent to
\be\label{Nperp}
\int g\,dv=0,\qquad \int g\,v\,dv=0,\qquad \int g\,|v|^2dv=0,
\ee
so the collision invariants of $g$ vanish independently of parameters. Thus the varying-projection framework is consistent: $\mathfrak{N}^\bot$ equals the orthogonal complement of the macroscopic space for the global Maxwellian $\Mstar$, used implicitly in Lemma \ref{co-1}.
Equation \eqref{Nperp} also implies that $\mathfrak{N}^\bot$ is invariant under $\na_v$ and $\Delta_v$. For $g\in\mathfrak{N}^\bot$, $\langle\na_vg,\chi_0\rangle=0$, $\langle\pa_{v_k}g,\chi_i\rangle=-\d_{ik}(nT)^{-1/2}\int g\,dv=0$, and $\langle\pa_{v_k}g,\chi_4\rangle=-\f{2}{T\sqrt{6n}}\int g(v_k-w_k)dv=0$; similarly for $\Delta_vg$. Hence in Lemma \ref{ex}, $\Lop^{-1}_{[w,T]}$ can be applied to $\na_v(\Lop^{-1}_{[w,T]}h)$ and $\Delta_v(\Lop^{-1}_{[w,T]}h)$.

\end{remark}

Moreover, we have the following commutator-type identity
\begin{lemma}\label{ex}
For any $h(v)\in \mathfrak{N}^\bot$ and for $z=t$ or $x_i~(i=1,2,3)$, the following identity holds:
\be\label{L-ex}
\begin{array}{l}
\pa_z(\Lop^{-1}_{[w,T]}h)=\Lop^{-1}_{[w,T]}(\pa_z h)+\pa_z w\cdot\Lop^{-1}_{[w,T]}[\na_v(\Lop_{[w,T]}^{-1}h)]\\[2mm]
\di\qquad\qquad\quad
-\f{\pa_z T}{T}\Lop^{-1}_{[w,T]}[h-3\Lop^{-1}_{[w,T]}h-(v-w)\cdot\na_v(\Lop^{-1}_{[w,T]}h)]
\end{array}
\ee

\end{lemma}

\noindent {\it Proof.} 
On one hand, we have for any $f(v)$,
\be\label{l1}
\pa_z(\Lop_{[w,T]} f)=\Lop_{[w,T]}(\pa_z f)-\pa_z w\cdot\na_v f+\pa_z T\Delta_v f.
\ee
On the other hand,
\be\label{l2}
\Lop_{[w,T]}f=3f+(v-w)\cdot\na_v f+T\Delta_v f,
\ee
which together with \eqref{l1} gives
$$
\pa_z(\Lop_{[w,T]} f)=\Lop_{[w,T]}(\pa_z f)-\pa_z w\cdot\na_v f
+\f{\pa_z T}{T}[\Lop_{[w,T]}f-3f-(v-w)\cdot\na_v f].
$$
Setting $h=\Lop_{[w,T]} f$, so that $f=\Lop_{[w,T]}^{-1}h$, and combining the above identity leads to
$$
\begin{array}{l}
\Lop_{[w,T]}^{-1}(\pa_z h)=\pa_z(\Lop_{[w,T]}^{-1}h)-\pa_z w\cdot\Lop_{[w,T]}^{-1}[\na_v(\Lop_{[w,T]}^{-1}h)]\\
\di\qquad\qquad\quad+\f{\pa_z T}{T}\Lop^{-1}_{[w,T]}[h-3\Lop^{-1}_{[w,T]}h-(v-w)\cdot\na_v(\Lop^{-1}_{[w,T]}h)].
\end{array}
$$
This yields the desired identity \eqref{L-ex}, and completes the proof of Lemma \ref{ex}.  \hfill $\Box$

\subsection{Planar rarefaction wave}
\label{SScons-rare}

In this subsection, we construct the smooth rarefaction wave and summarize its properties. We consider the rarefaction wave solutions to the Euler systems \eqref{Euler-fluid}--\eqref{Riemann-fluid} and \eqref{Euler-part}--\eqref{Riemann-part} for a polytropic gas for $\gamma=5/3$ and $R=1$. Recalling \eqref{constitutive}, we have
$$
S=-\ln \r+\f32\ln \t +\f32\ln \f{1}{A},\qquad p=\r\t=A\r^{\f53}e^{\f23S},
$$
for some positive constant $A$.  A direct calculation shows that each Euler system \eqref{Euler-fluid} for $(\r,u_1, S)$ has three distinct eigenvalues
$$
\ba
\lambda_i(\r,u_1,S)&=u_1+(-1)^{\f{i+1}{2}}\sqrt{\pa_{\r}p(\r,S)}
=u_1+(-1)^{\f{i+1}{2}}\sqrt{A\f53\r^{\f23}e^{\f23S}}\\
&=u_1+(-1)^{\f{i+1}{2}}\sqrt{\f53\t}, \quad i=1,3,\\
\lambda_2(\r,u_1,S)&=u_1,
\ea
$$
with the corresponding right eigenvectors
\be\label{eigen-vec}
r_i(\r,u_1,S)=((-1)^{\f{i+1}{2}}\r,\sqrt{\pa_{\r}p(\r,S)},0)^t,~i=1,3, \qquad\quad
r_2(\r,u_1,S)=(\pa_{\scriptscriptstyle S}p,0,-\pa_{\r}p)^t,
\ee
satisfying
$$
r_i(\r,u_1,S)\cdot \nabla_{(\r,u_1,S)}\lambda_i(\r,u_1,S)\neq 0,~ i=1,3,\quad
{\rm
and}\quad
r_2(\r,u_1,S)\cdot \nabla_{(\r,u_1,S)}\lambda_2(\r,u_1,S)\equiv 0
$$
in the non-vacuum region. Thus the two $i$-Riemann invariants $\Lambda_i^{(j)}(\r,u_1,S)$ can be defined as (cf. \cite{Smoller})
\begin{equation}\label{RI}
\Lambda_i^{(1)}(\r,u_1,S)=u_1+(-1)^{\f{i-1}{2}}\int^{\r}\f{\sqrt{\pa_{z}p(z,S)}}{z}dz,\qquad
\Lambda_i^{(2)}(\r,u_1,S)=S,\qquad i=1,3,
\end{equation}
such that
$$
\nabla_{(\r,u_1,S)} \Lambda_i^{(j)}(\r,u_1,S)\cdot r_i(\r,u_1,S)\equiv0,\quad i=1,3,~ j=1,2.
$$
To compute the eigenvalues and eigenvectors of the Euler system \eqref{Euler-part}, we define the macroscopic entropy by
\begin{equation*}
-\frac{3}{2}n \wt S:=\int\Mloc\ln \Mloc d v=n\Big(\ln n-\f32\ln(2\pi T)-\f32\Big),
\end{equation*}
which yields
\be\label{p-entropy}
 \wt S:=-\f23\ln n+\ln(2\pi T)+1,\quad  \tilde p=nT=\f{1}{2\pi e} n^{\f53}e^{\wt S}.
\ee
A direct computation shows that the Euler system \eqref{Euler-part} has three distinct eigenvalues in the variables $(n,w_1, \wt S)$:
$$
\ba
\lambda_i(n,w_1,\wt S)&=w_1+(-1)^{\f{i+1}{2}}\sqrt{\pa_n \tilde p(n, \wt S)}
=w_1+(-1)^{\f{i+1}{2}}\sqrt{\f53\f{1}{2\pi e}n^{\f23}e^{\wt S}}\\
&=w_1+(-1)^{\f{i+1}{2}}\sqrt{\f53 T}, \quad i=1,3,\\
\lambda_2(\r,u_1,\wt S)&=w_1.
\ea
$$
The corresponding right eigenvectors $r_i(n,w_1, \wt S)$ and $i$-Riemann invariants $\Lambda_i^{(j)}(n,w_1, \wt S)$ $(i=1,3,~ j=1,2)$  can be defined in a manner analogous to \eqref{eigen-vec} and \eqref{RI}.

Without loss of generality, we consider the time-asymptotic stability of the planar $3$-rarefaction wave for the Euler systems \eqref{Euler-fluid}--\eqref{Riemann-fluid}
and \eqref{Euler-part}--\eqref{Riemann-part}. The corresponding results for the planar $1$-rarefaction wave and for the superposition of planar $1$-rarefaction wave and $3$-rarefaction wave follow by similar argument. The planar $3$-rarefaction wave for the Euler systems \eqref{Euler-fluid}--\eqref{Riemann-fluid} and \eqref{Euler-part}--\eqref{Riemann-part} can be expressed explicitly in terms of the Riemann problem for the following one-dimensional inviscid Burgers equation:
\begin{equation}\label{bur}
\left\{
\begin{array}{ll}
B_t+BB_{x_1}=0,\\
B(0,x_1)=\left\{\begin{array}{ll}
B_-,&x_1<0,\\
B_+,&x_1>0.
\end{array}
\right.
\end{array}
\right.
\end{equation}
When $B_-<B_+$, the Riemann problem $(\ref {bur})$ admits a unique self-similar rarefaction wave fan solution $B^r(t, x_1) = B^r(x_1/t)$ given by
$$
B^r\Big(\f {x_1}{t}\Big)=\left\{\begin{array}{lr}
B_-,&\f {x_1}{t}\leq B_-,\\
\f {x_1}{t},&B_-\leq \f {x_1}{t}\leq B_+,\\ B_+,&\f {x_1}{t}\geq B_+.
\end{array}
\right.
$$
Then the planar $3$-rarefaction wave solutions $(\r^r,u^r,\t^r)(x_1/t)$ to \eqref{Euler-fluid}--\eqref{Riemann-fluid}
 and  $(n^r,u^r,\t^r)(x_1/t)$ to \eqref{Euler-part}--\eqref{Riemann-part}
 are defined explicitly as follows.
\begin{eqnarray}   \label{3-rw}
\left\{
\begin{array}{l}
\di B_\pm=\lambda_3(\r_\pm,u_{1\pm},\t_\pm)=\lambda_3(n_\pm, u_{1\pm},\t_\pm)=u_{1\pm}+\sqrt{\f53\t_\pm}, \\
\di B^r(\f {x_1}{t})= \lambda_3(\r^r,u_1^r,\t^r)(\f {x_1}{t})=\lambda_3(n^r,u_1^r,\t^r)(\f {x_1}{t})=u_1^{r}(\f{x_1}{t})+\sqrt{\f53\t^r(\f {x_1}{t})},\\
\di  u_1^r(\f{x_1}{t})=3\sqrt{\f53A(\r^r(\f{x_1}{t}))^{\f23}e^{\f23S_+}}
+u_{1+}-3\sqrt{\f53A\r^{\f23}_+e^{\f23S_+}}\\
\di\qquad\quad\,\,=3\sqrt{\f53\f{1}{2\pi e}(n^r(\f{x_1}{t}))^{\f23}e^{\wt S_+}}
+u_{1+}-3\sqrt{\f53\f{1}{2\pi e} n^{\f23}_+e^{\wt S_+}},\\
\di u_i^r(\f {x_1}{t})=0,~i=2,3.
\end{array} \right.
\end{eqnarray} 
The end states $(\rho_\pm,n_\pm,u_\pm,\theta_\pm)$ satisfy
\begin{align}\label{endstate}
S_-=S_+,\quad \widetilde S_-=\widetilde S_+ \quad
\Big(\mathrm{equivalently}~\frac{\rho_-}{\rho_+}=\frac{n_-}{n_+}
=\Big(\frac{\theta_-}{\theta_+}\Big)^{\frac{3}{2}}
\Big),\\
\nm u_{1-}-3\sqrt{\frac{5\theta_-}{3}}=
u_{1+}-3\sqrt{\frac{5\theta_+}{3}},\quad
\theta_-<\theta_+.
\end{align}
Recall that a planar rarefaction wave is one-dimensional in the three-dimensional space and independent of the transverse variables $x_2$ and $x_3$.

In order to fully exploit the friction relaxation effect arising from the momentum and energy exchanges between fluid and particles in the stability analysis of the planar rarefaction wave, we require in \eqref{3-rw} that the velocity $u_1^r$ and temperature $\t^r$ of the two rarefaction coincide, while the densities may differ. Moreover, from $\eqref{3-rw}_3$, it follows that
$$
\sqrt{A e^{\f23 S_+}}\big((\r^r)^{\f13}-\r_+^{\f13}\big)=\sqrt{\f{e^{\wt S_+}}{2\pi e}}\big((n^r)^{\f13}-n_+^{\f13}\big).
$$
Using \eqref{constitutive} and \eqref{p-entropy}, we have
$$
Ae^{\f23S_+}=\r_+^{-\f23}\t_+ \quad{\rm and}\quad
\f{e^{\wt S_+}}{2\pi e}=n_+^{-\f23}\t_+.
$$
Consequently,
$$
\f{\r^{r}}{\r_+}=\f{n^r}{n_+},
$$
which implies that the densities $\rho^{r}(x_1/t)$ and $n^{r}(x_1/t)$ of the two rarefaction waves may differ, provided $\rho_{+}\neq n_{+}$.

Note that the planar rarefaction wave solution constructed in \eqref{3-rw} is only Lipschitz continuous at the edge of the rarefaction wave, and has a singularity at $t=0$, $x=(0,x_2,x_3)$ with any $(x_2,x_3)\in \mathbb{R}^2$. To obtain the existence of strong solutions of Euler/NS--VFP system, we need more regularity for the planar rarefaction wave. We therefore construct a smooth approximation of the $3$-rarefaction wave defined in \eqref{3-rw}, following the approach in \cite{MatsumuraNishihara1986}. The smooth approximate rarefaction wave is constructed as follows.
\begin{equation}\label{dbur}
\left\{
\begin{array}{l}
\di \bar B_{t}+\bar B\bar B_{x_1}=0,\\
\di \bar B( 0,x_1)=\bar B_0(x_1)=\f{B_++B_-}{2}+\f{B_+-B_-}{2}\tanh(x_1).
\end{array}
\right.
\end{equation}
If $B_-<B_+$, then the problem \eqref{dbur} admits a unique classical solution $\bar B(t,x_1)$ given by
\begin{equation}\label{b-s}
\bar B(t,x_1)=\bar B_0(x_0(t,x_1)),\qquad x=x_0(t,x_1)+\bar B_0(x_0(t,x_1))t.
\end{equation}

Correspondingly, the smooth $3$-rarefaction wave profiles $(\bar{\r},\bar{u},\bar\t)(t,x_1)$ for the full compressible Euler system \eqref{Euler-fluid}-\eqref{Riemann-fluid} and $(\bar{n},\bar{u},\bar\t)(t,x_1)$ for \eqref{Euler-part}-\eqref{Riemann-part} are defined by
\begin{eqnarray}\label{au}
\left\{
\begin{array}{l}
\di B_\pm=\lambda_3(\r_\pm,u_{1\pm},\t_\pm)=\lambda_3(n_\pm,u_{1\pm},\t_\pm)=u_{1\pm}+\sqrt{\f53\t_\pm}, \\
\di \bar B(1+t,x_1)=\lambda_3(\bar{\r},\bar{u}_1,\bar\t)(t,x_1)=\lambda_3(\bar{n},\bar{u}_1,\bar\t)(t,x_1)
=\bar u_{1}(t,x_1)+\sqrt{\f53\bar\t(t,x_1)},\\
\di  \bar u_1(t,x_1)=3\sqrt{\f53A \bar \r^{\f23}(t,x_1)e^{\f23S_+}}
+u_{1+}-3\sqrt{\f53A\r^{\f23}_+e^{\f23S_+}}\\
\di\qquad\qquad=3\sqrt{\f53\f{1}{2\pi e} \bar n^{\f23}(t,x_1)e^{\wt S_+}}
+u_{1+}-3\sqrt{\f53\f{1}{2\pi e} n^{\f23}_+e^{\wt S_+}},\\
\di \bar u_i(t,x_1)=0,~i=2,3,
\end{array} \right.
\end{eqnarray}
where $\bar B(t, x_1)$ is the solution of Burgers equation $(\ref{dbur})$ given in \eqref{b-s}.
As in the non-smooth case,
the velocity $\bar u_1$ and temperature $\bar\t$ of the two smooth $3$-rarefaction wave profiles in \eqref{au} are taken to be equal in order to exploit the relaxation friction damping of fluid--particle model, while the densities $\bar\r$ and $\bar n$ may differ if $\rho_+\neq n_+$.

We list some properties for the approximate rarefaction wave below.
\begin{lemma}[\cite{MatsumuraNishihara1986}]\label{appr}
The Burgers equation $(\ref{dbur})$ has a unique smooth global solution $\bar B(t,x_1)$ satisfying the following properties:
\begin{itemize}
\item[(1)] $B_-<\bar B(t,x_1)<B_+, \ \bar B_{x_1}(t,x_1)>0,$ \ for  $x_1\in\mathbb{R}, \ t\geq 0.$
\item[(2)] For any $\ t> 0$ and p $\in[1,\infty]$, there exists a constant $C_{p}$ such that
\begin{equation*}
\begin{array}{ll}
%\di \|\bar B(t,\cdot)-B^r(\frac\cdot t)\|_{L^p}\leq  C_p\v^{-\frac1p}(B_+-B_-),\\
\di
\| \bar B_{x_1}(t,\cdot)\|_{L^p}\leq  C_{p}\min\{(B_+-B_-),
(B_+-B_-)^{\frac1p}t^{-1+\frac1p}\}, \\
\di
  \| \bar B_{x_1x_1}(t,\cdot)\|_{L^p}\leq
C_{p}\min\{(B_+-B_-),t^{-1}\},\\
\di \| \bar B_{x_1x_1x_1}(t,\cdot)\|_{L^p}\leq
C_{p}\min\{(B_+-B_-),t^{-1}\},\\
\di
|\bar B_{x_1x_1}(t,x_1)|\leq C\bar B_{x_1}(t, x_1).
\end{array}
\end{equation*}
\item[(3)] The smooth rarefaction wave $\bar B(t,x_1)$ and the self-similar rarefaction wave fan $B^r(x_1/t)$ are time asymptotically equivalent, i.e.,
$$
\lim_{t\rightarrow+\infty}\sup_{x_1\in \mathbb{R}} |\bar B(t,x_1)-B^r(x_1/t)|=0.
$$
\end{itemize}
\end{lemma}
Correspondingly, the $3$-rarefaction wave profiles $(\bar\rho, \bar u, \bar \t)(t,x_1)$ and $(\bar n, \bar u, \bar \t)(t,x_1)$ for
the Euler systems \eqref{1DEuler-fluid} and \eqref{1DEuler-part} have the following properties.
\begin{lemma}[\cite{MatsumuraNishihara1986}]\label{appu}
Let $\varepsilon=|(\r_+-\r_-,n_+-n_-,u_+-u_-, \t_+-\t_-)|$ denote the strength of the $3$-rarefaction waves defined in (\ref{au}). Then the following properties hold:
\begin{itemize}
\item[(1)]  $\bar u_{1x_1}=\f34 \bar B_{x_1}>0$, \, for  $x_1\in\mathbb{R}, \ t\geq 0.$ \\
 $\bar u_{1x_1}=\sqrt{\f53Ae^{\f23S_+}}\bar\r^{-\f23}\bar\r_{x_1}
=\sqrt{\f53\f{1}{2\pi e}e^{\wt S_+}}\bar n^{-\f23}\bar n_{x_1}$, \,
$\bar u_{1x_1}=\f{\sqrt{15}}{2}\f{\bar\t_{x_1}}{\sqrt{\bar\t}}$,
\item[(2)] For any $\ t> 0$ and p $\in[1,\infty]$, there exists a constant $C_{p}$ such that
\begin{align} \nm
& \| (\bar\r,\bar n,\bar u_1)_{x_1}(t,\cdot)\|_{L^p}\leq
C_{p} \min\{\v,\v^{\frac1p}(1+t)^{-1+\frac1p}\},\\ \nm
\nm
&  \|(\bar\r,\bar n, \bar u_1)_{x_1x_1}(t,\cdot)\|_{L^p}\leq
C_{p} \min\{\v,(1+t)^{-1}\},\\ \nm
& \|(\bar\r,\bar n, \bar u_1)_{x_1x_1x_1}(t,\cdot)\|_{L^p}\leq
C_{p} \min\{\v,(1+t)^{-1}\},\\ \nm
& |\bar u_{1x_1x_1}(t,x_1)|\leq C\bar u_{1x_1}(t,x_1).
\end{align}
\item[(3)] The smooth $3$-rarefaction wave is time-asymptotically equivalent to the self-similar $3$-rarefaction wave:
\begin{equation*}
\lim_{t\rightarrow+\infty}\sup_{x_1\in\mathbb{R}}| (\bar\r, \bar n, \bar
u, \bar \t)(t,x_1)-(\r^r, n^r, u^r, \t^r)(x_1/t)|=0.
\end{equation*}
\end{itemize}
\end{lemma}

\begin{remark}\label{rk:eps}
%We use a single small parameter: the smoothing scale in \eqref{dbur} is normalized to one and $\eps$ is the wave strength, so that
%$\|\bar B_0^{(k)}\|_{L^p}\leq C\eps$ for every $k\geq1$ and every $p$. The powers of the smoothing parameter that appear in the two-parameter construction
%of \cite{MatsumuraNishihara1986} are therefore not available here, but they are
%not needed: all the smallness consumed in Section \ref{Slow-es} is produced by the $\min$ structure of Lemma \ref{appu}(2) itself. Indeed, for the second and
%third derivatives,
%$$
%\int_0^\infty\|(\bar\r,\bar n,\bar u_1,\bar\t)_{x_1x_1}(\tau,\cdot)\|^2d\tau
%\leq C\int_0^\infty\min\{\eps^2,(1+\tau)^{-2}\}d\tau
%\leq C\Big(\eps^2\cdot\f1\eps+\eps\Big)\leq C\eps,
%$$
%the crossover occurring at $1+\tau=\eps^{-1}$, and likewise for $(\bar\r,\bar n,\bar u_1,\bar\t)_{x_1x_1x_1}$; while for the first derivative
%$\|\bar u_{1x_1}\|^2_{L^2_{x_1}}\leq C\min\{\eps^2,\eps(1+t)^{-1}\}$ is not integrable in time, which is precisely why the correction $\bar{\mathcal{G}}$ of
%\eqref{Gbar-intro} is introduced. Finally $\min\{\eps,(1+t)^{-1}\}\leq\eps^{\f18}(1+t)^{-\f78}$, which is the origin of the
%rate $\eps^{\f18}$ in \eqref{uniform}.
We use one small parameter: the smoothing scale in \eqref{dbur} is normalized to one, and $\eps$ is the wave strength, giving $\|\bar B_0^{(k)}\|_{L^p}\leq C\eps$ for every $k\geq1$ and $p$. Thus, unlike the two-parameter construction in \cite{MatsumuraNishihara1986}, smoothing-parameter powers are unavailable but unnecessary: all smallness in Section \ref{Slow-es} comes from the min structure of Lemma \ref{appu}(2). For the second and third derivatives,
$$
\int_0^\infty\| (\bar\r,\bar n,\bar u_1,\bar\t)_{x_1x_1}(\tau,\cdot)\|^2d\tau
\leq C\int_0^\infty\min\{\eps^2,(1+\tau)^{-2}\}d\tau
\leq C(\eps^2\cdot\f1\eps+\eps)\leq C\eps,
$$
with crossover at $1+\tau=\eps^{-1}$; the same holds for $(\bar\r,\bar n,\bar u_1,\bar\t)_{x_1x_1x_1}$. The first-derivative bound $\|\bar u_{1x_1}\|^2_{L^2_{x_1}}\leq C\min\{\eps^2,\eps(1+t)^{-1}\}$ is not time-integrable, motivating the correction $\bar{\mathcal{G}}$ in \eqref{Gbar-intro}. Finally, $\min\{\eps,(1+t)^{-1}\}\leq\eps^{\f18}(1+t)^{-\f78}$, yielding the rate $\eps^{\f18}$ in \eqref{uniform}.

\end{remark}

We conclude with the Gagliardo--Nirenberg inequality on  $\Omega:=\mathbb{R}\times \mathbb{T}^2$, which has been applied to the stability analysis of planar rarefaction wave for the three-dimensional Boltzmann equation \cite{WangWang2019} and related coupled models \cite{WangWang2020}, as well as for the isentropic Euler/NS--VFP system \cite{LiWangWang2022}. For the proof, we refer to \cite{LiWangWang2022}.

\begin{lemma}[\cite{LiWangWang2022,WangWang2019}]\label{sobolev}
For $g\in H^2(\Omega)$, $x=(x_1,x_2,x_3)\in \Omega$,  the following Gagliardo-Nirenberg inequality holds:
\begin{equation}\label{GN}
\|g\|^2_{L^{\infty}(\Omega)}\leq 2\|g\|_{L^2(\Omega)}\|\pa_{x_1}g\|_{L^2(\Omega)}
+ C\|\na_x g\|_{L^2(\Omega)}\|\na_x^2 g\|_{L^2(\Omega)},
\end{equation}
where $C$ is a positive constant.
\end{lemma}

%%%%%%%%%%%%%%%%%%%%%%%%%%%%%%%%%%%%%%%%%%%%%%%%%%%%%%%%%%%%%%%%%%%%%%%%%%%%%%%%%%%

\section{Nonlinear Stability}
\label{Sstability}

In this section, we prove Theorem \ref{thm:main} on the time-asymptotic stability of the planar rarefaction wave for the 3D full compressible Euler--VFP system \eqref{EVFP}. The proof combines the local existence of classical solutions with uniform a priori estimates. The local well-posedness can be established by adapting the iteration argument developed in \cite[Lemma 2.1]{WangWang2019}, especially for the kinetic part; see also \cite[Theorem 2.7]{LiNiWu1}.
For completeness, we state the local-in-time existence of a solution for the three-dimensional compressible Euler--VFP system  \eqref{EVFP} and \eqref{initial}.
\begin{lemma}\label{le-local}
{\rm (Local-in-time existence)} For any suitable small constant $\Xi>0$, there exists a positive constant $\mathcal{T}^*(\Xi)>0$, such that if the initial data $f_0(x,v)\geq0$ and
$$
\ba
&\|(\r_0(x)-\bar\r(0,x_1),u_0(x)-\bar u(0,x_1),\t_0(x)-\bar \t(0,x_1))\|_{H^3_x}\\
&+\big\|f_0(x,v)-\Mloc_{[\bar n(0,x_1),\bar u(0,x_1),\bar \t(0,x_1)]}\big\|_{H_x^3(L_v^2(1/\sqrt{\Mloc_*}))}
\leq \frac{\Xi}{2\sqrt{C_0}},
\ea
$$
where $C_0:=1/\min\{1,\bar\sigma\}\geq 1$ and the positive constant $\bar\sigma$ is defined in Lemma \ref{co-1}. Then the initial value problem \eqref{EVFP} and \eqref{initial} admits a unique solution $(\r,u,\t,f)$ satisfying
\begin{equation}
\begin{array}{l}
\di \sup_{0\leq t\leq \mathcal{T}^*(\Xi)}\|(\r(t,x)-\bar\r(t,x_1),u(t,x)-\bar u(t,x_1),\t(t,x)-\bar \t(t,x_1))\|^2_{H^3_x}
\leq \Xi,\\
\di \sup_{0\leq t\leq \mathcal{T}^*(\Xi)}\big\|f(t,x,v)-\Mloc_{[\bar n(t,x_1),\bar u(t,x_1),\bar \t(t,x_1)]}\big\|^2_{H_x^3(L_v^2(1/\sqrt{\Mloc_*}))}
\leq \Xi.
\end{array}
\end{equation}
\end{lemma}

It remains to establish the time-independent estimates \eqref{uniform}, after which global existence follows by a standard continuity argument.
For definiteness, we fix the ordering of the small parameters here, which is used without further comment in Sections \ref{Sstability} and
\ref{Slow-es}. Let $C_1$ denote the constant on the right-hand side of
\eqref{full-es}. Given $\mathcal{E}(0)$ and $\eps$ as in Theorem \ref{thm:main},
we set
$$
\d^2:=2C_1\big(\mathcal{N}(0)+\eps^{\f18}\big),
$$
so that, by \eqref{initial-relation}, $\d^2\leq C(\mathcal{E}(0)+\eps^{\f18})$ and
hence $\d\to0$ as $\eps_0\to0$; with this choice, \eqref{full-es} improves the a
priori assumption \eqref{assumption} by a factor $2$ on any interval on which the
latter holds, which is what the continuation argument requires. The parameter
$\eps_0$ is then fixed small enough that, first, all the absorptions performed in
Section \ref{Slow-es} are admissible, i.e., each quantity
$C(\d+\eps^{\f12})$ occurring there is smaller than the corresponding coercive
constant, and second, $C(\d+\eps)\leq\eta_0$ with $\eta_0$ as in Lemma
\ref{co-1}, so that Remark \ref{rk:Mstar} applies. Since $\d$ and $\eps$ are then
both controlled by $\eps_0$, no circularity is involved.
The wave phenomenon is distinct from the diffusive behavior of the damped Euler system or the pure Fokker--Planck equation, arising instead from fluid--particle interactions via the relaxation drag force in the fluid  equations and the Vlasov force in the kinetic equation.

First, we introduce the perturbation
\be\label{perturb}
\begin{cases}
(\phi,\psi,\z,\tilde \phi,\tilde \psi,\tilde \z)(t,x)=(\r-\bar\r, u-\bar u,\t-\bar\t,n-\bar n, w-\bar u,T-\bar\t)(t,x),\\
\widetilde{G}(t,x,v)=G(t,x,v)-\bar{\mathcal{G}}(t,x,v),
\end{cases}
\ee
where the correction function $\bar{\mathcal{G}}(t,x,v)$ is defined by
\begin{equation*}\label{bG}
\bar{\mathcal{G}}(t,x,v)=\f{1}{\r T}\Lop_{[w,T]}^{-1}\Big\{\mb{P}_1\Big[v_1\Big(\f{|v-w|^2}{2T}\bar\t_{x_1}+v_1\bar u_{1x_1}\Big)\Mloc_{[n,w,T]}\Big]\Big\}.
\end{equation*}
This correction is necessary since $\|\bar u_{1x_1}\|_{L^2_{x_1}}^2\sim (1+t)^{-1}$ is not integrable in $t$. Here $(n,w,T)$ are defined in \eqref{npart}.

From the system \eqref{twophase}$_{1,2,3}$, \eqref{ns-1} and \eqref{1DEuler-fluid}--\eqref{1DEuler-part}, we obtain the fluid perturbation $(\phi,\psi,\z,\tilde\phi,\tilde\psi,\tilde\z)$:
\begin{align}
\begin{cases}\label{perturb-0}
 & \phi_t+u\cdot\na\phi+\r\dv\psi+\psi\cdot\na\bar\r+\phi\dv\bar u=0,\\
 & \r\psi_t+\r u\cdot\na\psi+\t\na\phi+\r\na\z+\r\psi\cdot\na\bar u+\big(\z-\f{\bar\t}{\bar\r}\phi\big)\na\bar\r
=\r n(w-u),\\
&\f32\r\z_t+\f32\r u\cdot\na\z+\r\t\dv\psi+\f32\r\psi\cdot\na\bar\t+\r\z\dv\bar u=3\r n(T-\t)+\r n|w-u|^2,\\
&\tilde\phi_t+w\cdot\na\tilde\phi+ n\dv\tilde\psi+\tilde\psi\cdot\na\bar n+\tilde\phi\dv\bar u=0, \\
& n\tilde\psi_t+n w\cdot\na\tilde\psi+T\na\tilde\phi+n \na\tilde\z+n\tilde\psi\cdot\na\bar u
+\big(\tilde\z-\f{\bar\t}{\bar n}\tilde\phi\big)\na\bar n
=-\r n(w-u)- \int v\otimes v\cdot\na_xGdv,\\
&\f32n\tilde\z_t+\f32nw\cdot\na\tilde\z+nT\dv\tilde\psi+\f32n\tilde\psi\cdot\na\bar\t+n\tilde\z\dv\bar u=-3\r n(T-\t)\\
&\qquad\qquad\qquad\qquad
-\int \f12|v|^2v\cdot\na_xG dv+w\cdot\int v\otimes v\cdot\na_xG dv
\end{cases}
\end{align}
with initial data
\be\label{per-i}
\left\{
\begin{array}{l}
(\phi,\psi,\z)(t=0,x)=(\r_0(x)-\bar\r_0(x_1), u_0(x)-\bar u_0(x_1),\t_0(x)-\bar \t_0(x_1)),\\
(\tilde\phi,\tilde\psi,\tilde\z)(t=0,x)=(n_0(x)-\bar n_0(x_1), w_0(x)-\bar u_0(x_1),T_0(x)-\bar \t_0(x_1)),
\end{array}
\right.
\ee
where $(\r_0(x), u_0(x), \t_0(x))$ are given in $\eqref{initial}_1$, and
\be\label{n0}
\ba
 &n_0(x):=\int f_0(x,v)dv,\quad n_0(x)w_0(x):=\int vf_0(x,v)dv,\\
 &n_0(x)\Big(\f32T_0(x)+\f{|w_0(x)|^2}{2}\Big):=\int_{\mathbb{R}^3}\f12|v|^2f_0(x,v)dv
\ea
\ee
with $f_0(x,v)$ as in $\eqref{initial}_2$.
From \eqref{non-fluid} and \eqref{perturb}, the equation for $\widetilde{G}$ is
\begin{equation}\label{tilde-G}
\ba
\widetilde{G}_t-\r \Lop_{[w,T]}\widetilde{G}
&=-\f{1}{T}\mb{P}_1\big[v\cdot\big(\na\tilde\psi\cdot v+\f{|v-w|^2}{2T}\na\tilde\z\big)\Mloc\big]
-\mb{P}_1(v\cdot\na_xG)\\
&\quad\,-\bar{\mathcal{G}}_t+\r(w-u)\cdot\na_vG-\r(T-\t)\Delta_vG
\ea
\end{equation}
with initial data
\begin{equation}\label{Gi}
\widetilde{G}(t=0,x,v)=G(t=0,x,v)-\bar{\mathcal{G}}(t=0,x,v),
\end{equation}
where
\be\label{G-initial}
G(t=0,x,v)=f_0(x,v)-\Mloc_{[n_0(x), w_0(x),T_0(x)]}(v).
\ee
As mentioned above, to prove Theorem \ref{thm:main}, it suffices to establish  suitable a priori estimates and to verify the stability of planar rarefaction wave in \eqref{asymptotic} under the following a priori assumption:
\begin{equation}\label{assumption}
\ba
\mathcal{N}(\mathcal{T})=\sup_{0\leq t \leq\mathcal{T}}
\Big\{
\|(\phi,\psi,\z)(t)\|^2_{H^3}+\|(\tilde\phi,\tilde\psi,\tilde\z)(t)\|^2_{H^2}+\|\na^2(\phi_t,\psi_t,\z_t)(t)\|^2\\
\di+\int\int\frac{|\na_x^2\partial^{\alpha} f(t)|^2}{\Mloc_*}
+\f{|(\na_x\partial^{\alpha}G,\partial^{\alpha}G,\wt{G})(t)|^2}{\Mloc_*}dvdx
\Big\}\leq \d^2\ll1,
\ea
\end{equation}
where $\partial^{\alpha}=\partial_{t}$ or $\na_x$ and $\d>0$ is a small positive constant depending only on the initial perturbation and the wave strength $\varepsilon$.

From \eqref{smallness} and \eqref{assumption}, the initial data of the perturbation  $(\phi,\psi,\z,\tilde\phi,\tilde\psi,\tilde\z,\widetilde{G})$ satisfy
\be\label{initial-relation}
\mathcal{N}(0)\leq C(\mathcal{E}(0)+\varepsilon) \quad{\rm and} \quad
\mathcal{E}(0)\leq C(\mathcal{N}(0)+\varepsilon).
\ee

For $\d$ sufficiently small, the condition $\di\sup_{0\leq t\leq T}\|(\phi,\z,\tilde\phi,\tilde\z)(t)\|_{H^2}\leq\d$
together with the Sobolev embedding theorem implies $|\phi|\leq \f12\r_-$, $|\tilde\phi|\leq \f12 n_-$, $|(\z,\tilde\z)|\leq \f12\t_-$,
and $|(u,w)|\leq C$, where $C>0$ depends on $\r_-,n_-,\t_-, u_{\pm}$. Consequently, since $0<\rho_-<\bar\r<\rho_+$ and $0<n_-<\bar n<n_+$,  the densities
$\r=\bar \r+\phi$ and $n=\bar n+\tilde \phi$ satisfy
\be\label{rho-bd}
0<\f12\r_-< \r< \f12\r_-+\r_+, \quad 0<\f12 n_-<  n< \f12 n_-+ n_+.
\ee
Similarly, the temperature $\t=\bar\t+\z$ and $T=\bar\t+\tilde\z$ satisfy
\be\label{t-bd}
0<\f12\t_-< \t, \ T< \f12\t_-+\t_+,
\ee
since $0<\t_-< \bar \t< \t_+$. Moreover, by the three-dimensional Sobolev inequality, it holds
\begin{equation}\label{assum-1}
\ba
\sup_{t\in[0,\mathcal{T}]}\|\pa^{\alpha}(\r,u,\t,n,w,T)(t)\|_{L^{\infty}_x}
&\leq C\sup_{t\in[0,\mathcal{T}]}(\|\pa^{\alpha}(\phi,\psi,\z,\tilde\phi,\tilde\psi,\tilde\z)(t)\|_{L^{\infty}_x}+\|\pa^{\alpha}\bar u(t)\|_{L^{\infty}_x})\\
&
\leq C \sup_{t\in[0,\mathcal{T}]}\|\pa^{\alpha}(\phi,\psi,\z,\tilde\phi,\tilde\psi,\tilde\z)(t)\|_{H^{2}_x}+C\v\leq C(\d+\v),
\ea
\end{equation}
\begin{equation}\label{assum-2}
\ba &\sup_{t\in[0,\mathcal{T}]}\Big\|\int\f{|\pa^{\alpha}G(t)|^2}{\Mloc_*}dv\Big\|_{L^{\infty}_x}
=\sup_{t\in[0,\mathcal{T}]}\int\Big\|\f{\pa^{\alpha}G(t)}{\sqrt{\Mloc_*}}\Big\|^2_{L^{\infty}_x} dv
\leq C\sup_{t\in[0,\mathcal{T}]}\int\Big\|\f{\pa^{\alpha}G(t)}{\sqrt{\Mloc_*}}\Big\|^2_{H^{2}_x} dv\\
&\qquad \leq C\sup_{t\in[0,\mathcal{T}]}\int\int \f{|\pa^{\alpha}G(t)|^2+|\na_x\pa^{\alpha}G(t)|^2+|\na^2_x\pa^{\alpha}G(t)|^2}{\Mloc_*}dvdx
\leq C(\d+\v)^2,
\ea
\end{equation}
where $C>0$ is a generic constant. Since
$$ \sup_{t\in[0,\mathcal{T}]}\Big\|\na_x^2\pa^{\alpha}\Big(n,nw,n(\f32T+\f{|w|^2}{2})\Big)(t)\Big\|^2\leq C\int\int\f{|\na_x^2\pa^{\alpha}f(t)|^2}{\mathbf{M}_*} dvdx,
$$
which together with \eqref{assumption} yields
\begin{align}\label{assum-3}
&\sup_{t\in[0,\mathcal{T}]}\|\na_x^2\pa^{\alpha}(n,w,T)(t)\|^2\\
\nm &\leq C\sup_{t\in[0,\mathcal{T}]}\Big(\Big\|\na_x^2\pa^{\alpha}\Big(n,nw,n(\f32T+\f{|w|^2}{2})\Big)(t)\Big\|^2
+\||\na_{x}^2(n,w,T)|\,|\pa^{\alpha}(n,w,T)|(t)\|^2\Big)\\
\nm &+C\sup_{t\in[0,\mathcal{T}]}(\||\na_{x}\pa^{\alpha}(n,w,T)|\,|\na_x(n,w,T)|(t)\|^2
+\||\na_{x}(n,w,T)|^2|\pa^{\alpha}(n,w,T)|(t)\|^2)
\leq C(\d+\v)^2,
\end{align}
and
\be\label{assum-4}
\sup_{t\in[0,\mathcal{T}]}\|\na_x^2\pa^{\alpha}(\r,u,\t)(t)\|^2\leq
C\sup_{t\in[0,\mathcal{T}]}(\|\na_x^2\pa^{\alpha}(\phi,\psi,\z)(t)\|^2+\|\na_{x}^2\pa^{\alpha}\bar u\|^2)\leq C(\d+\v)^2,
\ee
where $C>0$ is a generic constant. Under the a priori assumption \eqref{assumption}, we have the following estimates.

\begin{theorem}\label{priori}
{\rm (A priori estimates)} Let $\mathcal{T}>0$, and suppose the assumptions of Theorem  \ref{thm:main} hold. Let $(\r,u,\t,f)$ be the strong solution of the 3D Euler--VFP system \eqref{EVFP}, \eqref{initial}, with the corresponding perturbations $(\phi,\psi,\z,\tilde\phi,\tilde\psi,\tilde\z)$ and $\wt{G}$ solving the initial value problems \eqref{perturb-0}-\eqref{per-i} and \eqref{tilde-G}-\eqref{Gi}, respectively. Then, under the a priori assumption \eqref{assumption}, the following estimate holds for  $t\in[0,\mathcal{T}]:$
\begin{equation}\label{full-es}
\begin{array}{ll}
\di \|(\phi,\psi,\z)(t)\|_{H^3}^2+\|(\tilde\phi,\tilde\psi,\tilde\z)(t)\|_{H^2}^2+\|\na^2(\phi_t,\psi_t,\z_t)(t)\|^2\\[1mm]
\di+\int_0^t(\|\pa^{\alpha}(\phi,\psi,\z,\tilde\phi,\tilde\psi,\tilde\z)\|^2_{H^2}
+\|(\tilde\psi-\psi,\tilde\z-\z)\|^2_{H^3}+\|\sqrt{\bar u_{1x_1}}(\phi,\psi_1,\z,\tilde\phi,\tilde\psi_1,\tilde\z)\|^2)d\tau\\[1mm]
\di+\int\int\f{1}{\Mloc_*}|(\na_x^2\pa^{\alpha}f,\na_x\pa^{\alpha} G,\pa^{\alpha}G,\wt{G})(t)|^2dxdv\\[1mm]
\di+\int_0^t\int\int\f{1}{\Mloc_*}|\na_v(\na_x^2\pa^{\alpha}G,\na_x\pa^{\alpha}G,\pa^{\alpha}G,\wt{G})|^2dvdxd\tau\\[1mm]
\di+\int_0^t\int\int\f{1+|v|^2}{\Mloc_*}|(\na_x^2\pa^{\alpha}G,\na_x\pa^{\alpha}G,\pa^{\alpha}G,\wt{G})|^2
 dvdxd\tau
\leq C(\mathcal{N}(0)+\v^{\frac{1}{8}}),
\end{array}
\end{equation}
where $C>0$ is a generic constant independent of time $\mathcal{T}>0$ and $\pa^{\alpha}$ denotes $\pa_t$ or $\na_x$.
\end{theorem}

The proof of Theorem \ref{priori} will be given in Section \ref{Slow-es}.

\begin{remark}
Since $\mathcal{N}(0)\leq C(\mathcal{E}(0)+\varepsilon)$ by
\eqref{initial-relation}, and $\varepsilon\leq C\varepsilon^{\f18}$ for $\varepsilon\leq 1$, the right hand side of \eqref{full-es} is bounded by $C(\mathcal{E}(0)+\varepsilon^{\f18})$, where $\mathcal{E}(0)$ is defined in \eqref{E0}.
\end{remark}

Once Theorem \ref{priori} is proved, the proof of Theorem \ref{thm:main} is complete.

\noindent{{\emph{Proof of Theorem \ref{thm:main}}}}: The global-in-time existence of solution follows immediately
from the local-in-time existence Lemma \ref{le-local} combined with the uniform-in-time a priori estimates in Theorem \ref{priori}. It remains to justify the
time-asymptotic stability of the planar rarefaction wave stated in \eqref{asymptotic}. Indeed, from \eqref{full-es}, we have
\begin{align}\label{3.7} &\int_0^{\infty}\|\na(\phi,\psi,\z)\|^2d\tau+\int_0^{\infty}\Big|\f{d}{dt}\|\na(\phi,\psi,\z)\|^2\Big|dt\\
\nm & \leq C\int_0^{\infty}(\|\na(\phi,\psi,\z)\|^2+\|\na(\phi_t,\psi_t,\z_t)\|^2)dt\leq C,
\end{align}
and
\be\label{3.8}
\begin{cases}
\di \int_0^{\infty}\int\int \f{|\na_xG|^2}{\Mloc_*}dv dxdt+
\int_0^{\infty}\int\int \f{|\na_x(\Mloc-\Mloc_{[\bar n,\bar u,\bar\t]})|^2}{\Mloc_*}dv dxdt\leq C,\\
\di\int_0^{\infty}\Big|\f{d}{dt}\int\int \f{|\na_xG|^2}{\Mloc_*}dvdx\Big|dt
\leq \int_0^{\infty}\int\int \f{|\na_xG|^2+|\na_xG_t|^2}{\Mloc_*}dv dx dt\leq C,\\
\di\int_0^{\infty}\Big|\f{d}{dt}\int\int \f{|\na_x(\Mloc-\Mloc_{[\bar n,\bar u,\bar\t]})|^2}{\Mloc_*}dvdx\Big|dt\\
\di\leq \int_0^{\infty}\int\int \f{|\na_x(\Mloc-\Mloc_{[\bar n,\bar u,\bar\t]})|^2}{\Mloc_*}dvdxdt
+ \int_0^{\infty}\int\int \f{|\pa_t\na_x(\Mloc-\Mloc_{(\bar n,\bar u,\bar\t)})|^2}{\Mloc_*}dv dx dt\\
\di\leq C\int_0^{\infty}(\|\na(\tilde\phi,\tilde\psi,\tilde\z)\|^2+\|\na(\tilde\phi_t,\tilde\psi_t,\tilde\z_t)\|^2)dt\leq C.
\end{cases}
\ee
Then \eqref{3.7} and \eqref{3.8} indicate
\be\label{4.2}
\di \lim_{t\rightarrow\infty}\|\na(\phi,\psi,\z)\|^2=0,\quad
\lim_{t\rightarrow\infty}\int_{\mathbb{R}^3}\int\f{|\na_xG|^2
+|\na_x(\Mloc-\Mloc_{[\bar n,\bar u,\bar\t]})|^2}{\Mloc_*}dx dv=0.
\ee
By the three-dimensional Gagliardo--Nirenberg inequality \eqref{GN}, we have
$$
\ba
&\|(\phi,\psi,\z)\|_{L_x^{\infty}}\leq C\|(\phi,\psi,\z)\|^{\f12}\|\pa_{x_1}(\phi,\psi,\z)\|^{\f12}+
C\|\na_x(\phi,\psi,\z)\|^{\f12}\|\na^2_x(\phi,\psi,\z)\|^{\f12}\leq C,\\[2mm]
& \Big\|\int \f{|G|^2}{\Mloc_*} dv\Big\|_{L^{\infty}_x}
 \leq C\Big(\int\int \f{|G|^2}{\Mloc_*} dv dx\Big)^{\f12}\Big(\int\int \f{|\pa_{x_1}G|^2}{\Mloc_*} dv dx\Big)^{\f12}\\
&~~\qquad\qquad\qquad\quad
+C\Big(\int\int \f{|\na_xG|^2}{\Mloc_*} dv dx\Big)^{\f12}\Big(\int\int \f{|\na^2_xG|^2}{\Mloc_*} dv dx\Big)^{\f12},\\
\ea
$$
$$
\ba
&\quad\Big\|\int \f{|\Mloc-\Mloc_{[\bar n,\bar u,\bar\t]}|^2}{\Mloc_*} dv\Big\|_{L^{\infty}_x}\\
&\leq C\Big(\int\int \f{|\Mloc-\Mloc_{[\bar n,\bar u,\bar\t]}|^2}{\Mloc_*} dv dx\Big)^{\f12}
\Big(\int\int \f{|\pa_{x_1}(\Mloc-\Mloc_{[\bar n,\bar u,\bar\t]})|^2}{\Mloc_*} dv dx\Big)^{\f12}\\
&+C\Big(\int\int \f{|\na_x(\Mloc-\Mloc_{[\bar n,\bar u,\bar\t]})|^2}{\Mloc_*} dv dx\Big)^{\f12}
\Big(\int\int \f{|\na^2_x(\Mloc-\Mloc_{[\bar n,\bar u,\bar\t]})|^2}{\Mloc_*} dv dx\Big)^{\f12},
\ea
$$
which together with \eqref{full-es} and \eqref{4.2} yields
$$
\begin{array}{l}
\di\lim_{t\rightarrow\infty}\sup_{x\in\Omega}
|(\r(t,x)-\bar\r(t,x_1),u(t,x)-\bar u(t,x_1),\t(t,x)-\bar\t(t,x_1))|=0,\\
\di\lim_{t\rightarrow\infty}\sup_{x\in\Omega}\int \f{|G|^2+|\Mloc-\Mloc_{[\bar n,\bar u,\bar\t]}|^2}{\Mloc_*}
(t,x,v)dv=0.
\end{array}
$$
Thus
\begin{equation*}
\di \lim_{t\rightarrow\infty}\sup_{x\in\Omega}\int_{\mathbb{R}^3}\f{|f-\Mloc_{[\bar n,\bar u,\bar\t]}|^2}{\Mloc_*}(t,x,v)dv
\leq 2\lim_{t\rightarrow\infty}\sup_{x\in\Omega}\int_{\mathbb{R}^3} \f{|G|^2+|\Mloc-\Mloc_{[\bar n,\bar u,\bar\t]}|^2}{\Mloc_*}
(t,x,v)dv=0,
\end{equation*}
which verifies \eqref{asymptotic}. Thus, the proof of Theorem \ref{thm:main} is complete. \hfill $\Box$

 %
%%%%%%%%%%%%%%%%%%%%%%%%%%%%%%%%%%%%%%%%%%%%%%%%%%%%%%%%%%%%%%%%%%%%%%%%%%%%%%%%%%%%%%%%%%%%%%%%%%%%%%%%

\section{A priori Estimates }
\label{Slow-es}

In this section, we prove the a priori estimates in Theorem \ref{priori} for the 3D full compressible Euler--VFP system \eqref{EVFP}, by combining the basic energy estimates in Proposition \ref{Prop3.1} with  the higher-order estimates in Proposition \ref{Prop3.2} and \ref{Prop3.3}. This approach is motivated by the energy method for the Boltzmann equation developed in \cite{LiuYangYu,LiuYu}. In particular, we fully exploit the relaxation friction damping $w-u$ and $T-\t$ arising from the fluid-particle interactions, together with the viscous dissipation from the terms $\dv (nT\mathbb{S}(w)/\rho)$ and $\dv(n T\na T/\rho)$ in the new viscous compressible fluid model \eqref{mix-e} for $(n,w,T)$, in order to obtain the desired higher-order regularity estimates for $(\rho, u,\t)$ and $(n,w,T)$, relative to those for the microscopic part $G$. This observation is essential for establishing the global dynamical stability of the planar rarefaction wave for the 3D Euler--VFP system \eqref{EVFP}.

We list the lower-order estimates in Proposition \ref{Prop3.1} and the higher-order estimates in Propositions \ref{Prop3.2} and \ref{Prop3.3} below.

\begin{proposition}\label{Prop3.1}
Let $(\phi,\psi,\z,\tilde \phi, \tilde\psi,\tilde\z)$ be the strong solution of the initial value problem \eqref{perturb-0}--\eqref{per-i},  and let $\wt{G}$ be the strong solution of the initial value problem \eqref{tilde-G}--\eqref{Gi}. Then, under the a priori assumption \eqref{assumption}, the following estimate holds for $t\in[0,T]:$
\begin{align}\label{lower}
& \|(\phi,\psi,\z,\tilde\phi,\tilde\psi,\tilde\z)(t)\|^2
+\int\int\f{|\widetilde{G}(t)|^2}{\Mloc_*} dxdv
+\int_0^t\int\int\f{|\na_v\wt{G}|^2+(1+|v|^2)|\widetilde{G}|^2}{\Mloc_*}dvdxd\tau\\
\nm &+\int_0^t(\|\sqrt{\bar u_{1x_1}}(\phi,\psi_1,\z,\tilde\phi,\tilde\psi_1,\tilde\z)\|^2
+\|\pa^{\alpha}(\phi,\tilde\phi,\tilde\psi,\tilde\z)\|^2+\|(\psi_t,\z_t)\|^2+\|(\tilde\psi-\psi,\tilde\z-\z)\|^2)d\tau\\
\nm  &\leq C(\mathcal{N}(0)+\v^{\f{1}{8}})+C\|\na(\phi,\tilde\phi)(t)\|^2+C\int_0^t\|\na(w-u,T-\t)\|^2d\tau
 +C\int_0^t\int\int\f{|\pa^{\alpha}G|^2}{\Mloc_*}dvdxd\tau,
\end{align}
where $C>0$ is a generic constant independent of time $\mathcal{T}>0$ and $\pa^{\alpha}$ denotes $\pa_t$ or $\na_x$.
\end{proposition}

\begin{proposition} \label{Prop3.2}
Let $(\r,u,\t,n,w,T)$ be the strong solution of the system \eqref{twophase} for the Euler--VFP system, with the initial data $(\r_0, u_0,\t_0)$ given in \eqref{initial} and $(n_0,w_0,T_0)$ given in \eqref{n0}. Let $G$ be the strong solution of the equation \eqref{non-fluid} with the initial data \eqref{G-initial}. Then, under the a priori assumption \eqref{assumption}, the following estimate holds for $t\in[0,\mathcal{T}]:$
\begin{align}\label{high-1}
\nm &  \|\na(\r,u,\t,n,w,T)(t)\|^2
+\int_0^t(\|\na\pa^{\alpha}(\r,n,w,T)\|^2+\|\na( u_t,\t_t)\|^2+\|\na(w-u,T-\t)\|^2)d\tau\\
\nm &
+\int\int\f{|\pa^{\alpha}G(t)|^2}{\Mloc_*}dvdx
+\int_0^t\int\int\f{|\na_v\pa^{\alpha}G|^2+(1+|v|^2)|\pa^{\alpha}G|^2}{\Mloc_*}dvdxd\tau
\\
&\leq C(\mathcal{N}(0)+\v^{\f12})+C\|\na^2(\r,n)(t)\|^2+C\int_0^t\|\na^2(w-u,T-\t)\|^2d\tau\\
\nm &
+C(\d+\v^{\f12})\int_0^t(\|\pa^{\alpha}(\phi,\tilde\phi,\tilde\psi,\tilde\z)\|^2+\|(\psi_t,\z_t)\|^2+\|(\tilde\psi-\psi,\tilde\z-\z)\|^2)d\tau\\
\nm &+C\int_0^t\int\int\f{|\na_x\pa^{\alpha}G|^2}{\mb{M_*}}dvdxd\tau
+C(\d+\v^{\f12})\int_0^t\int\int\f{|(\na_v\na^2_xG,\na_v\wt{G},\wt{G})|^2}{\mb{M_*}}dvdxd\tau,
\end{align}
where $C>0$ is a generic constant independent of time $\mathcal{T}>0$, and $\pa^{\alpha}$ denotes $\pa_t$ or $\na_x$.
\end{proposition}

\begin{proposition} \label{Prop3.3}
Let $(\r,u,\t,f)$ be the strong solution of the 3D full compressible Euler--VFP system \eqref{EVFP} with the initial data \eqref{initial}. Let $(n,w,T)$ be the strong solution of system \eqref{mix-e} with the initial data \eqref{n0}, and let $G$ be the strong solution to the equation \eqref{non-fluid} with the initial data \eqref{G-initial}. Then, under the a priori assumption \eqref{assumption}, the following estimate holds for $t\in[0,\mathcal{T}]:$
\begin{align}\label{high-2}
\nm & \|\na^2(\r,u,\t)(t)\|_{H^1}^2+\|\na^2(n,w,T)(t)\|^2+\|\na^2(\r_t,u_t,\t_t)(t)\|^2\\
\nm &+\int\int \f{|(\na_x^2\pa^{\alpha}f,\na_x\pa^{\alpha}G)(t)|^2}{\Mloc_*}dv dx
+\int_0^t(\|\na^2\pa^{\alpha}(\r,n,w,T)\|^2+\|\na^2( u_t,\t_t)\|^2)d\tau\\
&+\int_0^t(\|\na^2\pa^{\alpha}(w-u,T-\t)\|^2+\|\na^2(w-u,T-\t)\|^2)d\tau\\
\nm &
+\int_0^t\int\int\f{|\na_v(\na^2_x\pa^{\alpha}G,\na_x\pa^{\alpha}G)|^2
+(1+|v|^2)|(\na^2_x\pa^{\alpha}G,\na_x\pa^{\alpha}G)|^2}{\Mloc_*}dvdxd\tau
\\
\nm &\leq C(\mathcal{N}(0)+\v)
+C(\d+\v^{\f12})\int_0^t(\|\na\pa^{\alpha}(\r,n,w,T)\|^2+\|\na(u_t,\t_t)\|^2+\|\na(w-u,T-\t)\|^2)d\tau\\
\nm &+C(\d+\v^{\f12})\int_0^t\|\pa^{\alpha}(\phi,\tilde\phi,\tilde\psi,\tilde\z)\|^2d\tau
+C(\d+\v)\int_0^t\int\int\f{|(\na_v\pa^{\alpha}G,\pa^{\alpha}G,\na_v\wt{G},\wt{G})|^2}{\mb{M_*}}dvdxd\tau,
\end{align}
where $C>0$ is a generic constant independent of time $\mathcal{T}>0$ and $\pa^{\alpha}$ denotes $\pa_t$ or $\na_x$.
\end{proposition}

Theorem \ref{priori} follows from Propositions \ref{Prop3.1}, \ref{Prop3.2} and \ref{Prop3.3} in the following way.

\noindent{{\emph{Proof of Theorem \ref{priori}:}}} First, combining \eqref{high-1} and \eqref{high-2}, and choosing $\d$ and $\v$ sufficiently small yields,
\begin{align}\label{high-3}\nm &\|\na(\r,u,\t)(t)\|_{H^2}^2+\|\na(n,w,T)(t)\|_{H^1}^2+\|\na^2(\r_t,u_t,\t_t)(t)\|^2\\ \nm
&+\int\int \f{|(\na_x^2\pa^{\alpha}f,\na_x\pa^{\alpha}G,\pa^{\alpha}G)(t)|^2}{\Mloc_*}dv dx+\int_0^t(\|\na\pa^{\alpha}(\r,n,w,T)\|_{H^1}^2+\|\na (u_t,\t_t)\|_{H^1}^2)d\tau\\
&+\int_0^t(\|\na(w-u,T-\t)\|_{H^1}^2+\|\na^2\pa^{\alpha}(w-u,T-\t)\|^2)d\tau\\ \nm
&+\int_0^t\int\int\f{|\na_v(\na^2_x\pa^{\alpha}G,\na_x\pa^{\alpha}G,\pa^{\alpha}G)|^2
+(1+|v|^2)|(\na^2_x\pa^{\alpha}G,\na_x\pa^{\alpha}G,\pa^{\alpha}G)|^2}{\Mloc_*}dvdxd\tau\\ \nm
&\leq C(\mathcal{N}(0)+\v^{\f12})+C(\d+\v^{\f12})\int_0^t\int\int\f{|(\na_v\wt{G},\wt{G})|^2}{\mb{M_*}}dvdxd\tau\\ \nm
&+C(\d+\v^{\f12})\int_0^t\big(\|\pa^{\alpha}(\phi,\tilde\phi,\tilde\psi,\tilde\z)\|^2
+\|(\psi_t,\z_t)\|^2+\|(\tilde\psi-\psi,\tilde\z-\z)\|^2\big)d\tau.
\end{align}
Then, combining \eqref{lower} and \eqref{high-3}, and using the identities
$$
\|\pa^{\alpha}\psi\|_{H^2}^2=\|\na\psi\|_{H^2}^2+\|\psi_t\|_{H^2}^2
\leq \|\na(\psi-\tilde\psi)\|_{H^2}^2+\|\na\tilde\psi\|_{H^2}^2+\|\psi_t\|_{H^2}^2
$$
and
$$
\|\pa^{\alpha}\z\|_{H^2}^2=\|\na\z\|_{H^2}^2+\|\z_t\|_{H^2}^2
\leq \|\na(\z-\tilde\z)\|_{H^2}^2+\|\na\tilde\z\|_{H^2}^2+\|\z_t\|_{H^2}^2,
$$
%and by \eqref{initial-relation}, we can replace $\mathcal{N}(0)$ by $C(\mathcal{E}^{e}(0)+\v)$,
we obtain the desired inequality \eqref{full-es} by choosing $\d$ and $\v$ suitably small. The proof of Theorem \ref{priori} is complete.

 \hfill $\Box$

\subsection{Basic energy estimates}
\label{SSpf-prop-1}

In this subsection, we prove  Proposition \ref{Prop3.1} on the basic energy estimates.

\noindent{{\emph{Proof of Proposition \ref{Prop3.1} }}}:
We first establish the $L^2$-entropy estimates for the macroscopic perturbation $(\phi,\psi,\z,\tilde\phi,\tilde\psi,\tilde\z)$. Let $(\eta,q)$ and $(\tilde\eta,\tilde q)$ denote the relative entropy-entropy flux pairs as
$$
\left\{
\ba
\eta&=\bar\theta\left\{-\rho
S+\bar\rho\bar{S}+\nabla_\mathbf{U}(\rho
S)\Big|_{\mathbf{U}=\bar{\mathbf{U}}}\cdot(\mathbf{U}-\bar{\mathbf{U}})\right\},\\
q_j&=\bar\theta\left\{-\rho u_j
S+\bar\rho\bar{u}_j\bar{S}+\nabla_\mathbf{U}(\rho
S)\Big|_{\mathbf{U}=\bar{\mathbf{U}}}\cdot(\mathbf{V}_j-\bar{\mathbf{V}}_j)\right\},\quad j=1,2,3,
\ea
\right.
$$
and
$$
\left\{
\ba
\tilde\eta&=\bar\theta\left\{-\f32 n
\wt S+\f32\bar n\overline{\wt S}+\f32\nabla_\mathbf{X}(n
\wt S)\Big|_{\mathbf{X}=\bar{\mathbf{X}}}\cdot(\mathbf{X}-\bar{\mathbf{X}})\right\},\\
\tilde q_j&=\bar\theta\left\{-n w_j
\wt S+\bar n\bar{u}_j\overline{\wt S}+\nabla_\mathbf{X}(n
\wt S)\Big|_{\mathbf{X}=\bar{\mathbf{X}}}\cdot(\mathbf{Y}_j-\bar{\mathbf{Y}}_j)\right\},\quad j=1,2,3,
\ea
\right.
$$
where
$$
\begin{array}{l}
\displaystyle \mathbf{U}=(\rho,\r u_1,\r u_2,\r u_3,\r E)^t,\quad
\displaystyle \mathbf{V}=(\r u,\r u u_1+p\mathbb{I}_1,\r uu_2+p\mathbb{I}_2,\r u u_3+p\mathbb{I}_3,\r u E+pu)^t,
\end{array}
$$
and
$$
\begin{array}{l}
\displaystyle \mathbf{X}=(n,n w_1, n w_2, n w_3,n \wt E)^t,\quad
\displaystyle \mathbf{Y}=(n w, n w w_1+\tilde p\mathbb{I}_1,n ww_2+\tilde p\mathbb{I}_2,n w w_3+\tilde p\mathbb{I}_3, nw \wt E+ \tilde p w)^t
\end{array}
$$
with $\mathbb{I}_1=(1,0,0)^t, \mathbb{I}_2=(0,1,0)^t, \mathbb{I}_3=(0,0,1)^t$. Then, a direct calculation yields
$$
\displaystyle \pa_{\r}(\rho S)=S+\frac{|u|^2}{2\theta}-\f52,\qquad
\displaystyle \pa_{m_i}(\rho S)=-\frac{u_i}{\theta},~i=1,2,3,\qquad
\displaystyle \pa_{\mathcal{E}}(\rho S)=\frac{1}{\theta},
$$
$$
\displaystyle \pa_{n}(n \wt S)=\wt S+\frac{|w|^2}{3T}-\f53,\qquad
\displaystyle \pa_{\tilde m_i}(n \wt S)=-\frac{2w_i}{3T},~i=1,2,3,\qquad
\displaystyle \pa_{\wt{\mathcal{E}}}(n \wt S)=\frac{2}{3T},
$$
where $m_i=\r u_i$, $\tilde m_i=n w_i~(i=1,2,3)$,  $\mathcal{E}=\r E$, and $\wt{\mathcal{E}} =n \wt E$,  then it holds
$$
\left\{\begin{aligned}
 \eta&=\f32\rho\theta-\bar\theta\rho
S+\rho\Big[\Big(\bar{S}-\f52\Big)\bar\theta+\frac{|u-\bar u|^2}{2}\Big]
+\bar\rho\bar\theta\\
&=\rho\bar\t\Psi\left(\f{\bar\r}{\rho}\right)+\f32\rho\bar\t\Psi\left(\f{\t}{\bar\t}\right)
+\f12\rho|u-\bar u|^2,\\
q&=u\eta +(u-\bar u)(\rho\theta-\bar\rho\bar\theta),
\end{aligned}
\right.
$$
and
$$
\left\{\begin{aligned}
 \tilde \eta&=\f32nT-\f32\bar\theta n\wt S+n\Big[\Big(\f32\overline{\wt S}-\f52\Big)\bar\theta+\frac{|w-\bar u|^2}{2}\Big]
+\bar n\bar\theta\\
&=n\bar\t\Psi\left(\f{\bar n}{n}\right)+\f32n\bar\t\Psi\left(\f{T}{\bar\t}\right)
+\f12n|w-\bar u|^2,\\
\tilde q&=w\tilde\eta +(w-\bar u)(nT-\bar n\bar\theta),
\end{aligned}
\right.
$$
where $\Psi(\cdot)$ is the convex function
$$
\Psi(s)=s-\ln s-1.
$$

Next, a direct computation yields
\be\label{eta-1}
\begin{array}{l}
\di\quad \eta_t+\dv q-\big[\nabla_{(\bar\rho,\bar{u},\bar{S})}\eta\cdot(\bar\rho,\bar{u},\bar{S})_t
+\nabla_{(\bar\rho,\bar{u},\bar{S})}q\cdot(\bar\rho,\bar{u},\bar{S})_{x_1}\big]\\
\di=3\r n(T-\t)+\r nw\cdot(w-u)-\f{\bar\t}{\t}(3\r n(T-\t)+\r n|w-u|^2)-\r n\bar u_1(w_1-u_1),
\end{array}
\ee
and
\begin{align}\label{eta-2}
\nm &\quad \tilde\eta_t+\dv\tilde q-\big[\nabla_{(\bar n,\bar{u},\overline{\wt S})}\tilde\eta\cdot(\bar n,\bar{u},\overline{\wt S})_t
+\nabla_{(\bar n,\bar{u},\overline{\wt S})}\tilde q\cdot(\bar n,\bar{u},\overline{\wt S})_{x_1}\big]\\
&=-3\r n(T-\t)-\r nw\cdot(w-u)+\f{\bar\t}{T}3\r n(T-\t)+\r n\bar u_1(w_1-u_1)\\
\nm &+\f{5}{6}\Big[\dv(\f{n}{\r}\tilde\z\na\tilde\z)-\f{n\bar\t}{\r T}|\na\tilde\z|^2
+\bar\t_{x_1}\f{\tilde\z}{\r}(\pa_1\tilde\phi+\f{2n}{T}\pa_1\tilde\z-\f{n}{\r}\pa_1\phi)
+\f{\tilde\z}{\r}(\f{n}{T}\bar\t_{x_1}^2+\bar n_{x_1}\bar\t_{x_1}-\f{n}{\r}\bar\r_{x_1}\bar\t_{x_1}+n\bar\t_{x_1x_1})\Big]\\
\nm &+\Big[\dv\big[\f{nT}{2\r}\tilde\psi\cdot\mathbb{S}(\tilde\psi)\big]
-\f{n \bar\t}{2\r}\big(\f{|\na\tilde\psi+(\na\tilde\psi)^t|^2}{2}-\f23(\dv\tilde\psi)^2\big)
+\bar u_{1x_1}\f{2\tilde\psi_1}{3\r}(T\pa_1\tilde\phi+n\pa_1\tilde\z-\f{nT}{\r}\pa_1\phi)\\
\nm & +\bar u_{1x_1}\f{n\tilde\z}{\r}(\f53\pa_1\tilde\psi_1-\f13\dv\tilde\psi)
+\f{2\tilde\psi_1}{3\r}(T\bar n_{x_1}\bar u_{1x_1}+n\bar\t_{x_1}\bar u_{1x_1}-\f{nT}{\r}\bar\r_{x_1}\bar u_{1x_1}
+nT\bar u_{1x_1x_1})+\f{n\tilde\z}{\r}\f23\bar u_{1x_1}^2\Big]\\
\nm &+\Big[-\dv\big[\f{\tilde\z}{T}\int\f12|v|^2v\Pi dv\big]+\big(\int\f12|v|^2v\Pi dv\big)\cdot\na\big(\f{\tilde\z}{T}\big)
-\sum_{i,j=1}^3\pa_j\big(\f{\bar\t}{T}\tilde\psi_i\int v_iv_j\Pi dv\big)\\
\nm &+\sum_{i,j=1}^3\pa_j\big(\f{\bar\t}{T}\tilde\psi_i\big)\int v_iv_j\Pi dv+\dv\big[\f{\tilde\z\bar u_1}{T}\int v_1v\Pi dv\big]
-\big(\int v_1 v\Pi dv\Big)\cdot\na \big(\f{\tilde\z\bar u_1}{T}\big)
\Big].
\end{align}
By Lemma \ref{appu} (1), there exists a positive constant $C>0$ such that (cf. \cite{LiWangWang2018})
$$
\ba
&\quad -\big[\nabla_{(\bar\rho,\bar{u},\bar{S})}\eta\cdot(\bar\rho,\bar{u},\bar{S})_t
+\nabla_{(\bar\rho,\bar{u},\bar{S})}q\cdot(\bar\rho,\bar{u},\bar{S})_{x_1}\big]\\
& =\bar u_{1x_1}\Big[\rho\psi_1^2+\f23\rho\bar\t \Psi\Big(\f{\bar\rho}\rho\Big)
+\rho\bar\t\Psi\Big(\f{\t}{\bar\t}\Big)\Big]
+\bar\t_{x_1}\rho\psi_1\Big(\ln\f{\bar\r}{\r}+\f{3}{2}\ln\f{\t}{\bar\t}\Big)
 \geq C^{-1} \bar u_{1x_1}(\phi^2+\psi_1^2+\z^2),
\ea
$$
and
$$
\ba
&\quad
-\big[\nabla_{(\bar n,\bar{u},\overline{\wt S})}\tilde\eta\cdot(\bar n,\bar{u},\overline{\wt S})_t
+\nabla_{(\bar n,\bar{u},\overline{\wt S})}\tilde q\cdot(\bar n,\bar{u},\overline{\wt S})_{x_1}\big]\\
& =\bar u_{1x_1}\Big[n\tilde\psi_1^2+\f23 n\bar\t \Psi\Big(\f{\bar n}n\Big)
+n\bar\t\Psi\Big(\f{T}{\bar\t}\Big)\Big]
+\bar\t_{x_1}n\tilde\psi_1\Big(\ln\f{\bar n}{n}+\f{3}{2}\ln\f{T}{\bar\t}\Big)
 \geq C^{-1} \bar u_{1x_1}(\tilde\phi^2+\tilde\psi_1^2+\tilde\z^2).
\ea
$$
Adding \eqref{eta-1} and \eqref{eta-2}, and integrating the resulting equation by parts over $[0,t]\times\Omega$, we obtain
\begin{align}\label{E1}\nm
&\quad \int (\eta+\tilde\eta)dx\big|_{\tau=0}^{\tau=t}+C^{-1}\int_0^t\int \big(\bar u_{1x_1}|(\phi,\psi,\z,\tilde\phi,\tilde\psi,\tilde\z)|^2
+\f{3\bar\t}{\t T}\r n|T-\t|^2+\f{\bar\t}{\t}\r n|w-u|^2\big)dxd\tau\\ \nm
&+\int_0^t\int\Big[\f{5}{6}\f{n\bar\t}{\r T}|\na\tilde\z|^2+\f{n\bar\t}{2\r}
\Big(\f{|\na\tilde\psi+(\na\tilde\psi)^t|^2}{2}-\f23(\dv\tilde\psi)^2\Big)\Big]dxd\tau\\
&\leq \int_0^t\int \f{5}{6}\Big[\bar\t_{x_1}\f{\tilde\z}{\r}(\pa_1\tilde\phi+\f{2n}{T}\pa_1\tilde\z-\f{n}{\r}\pa_1\phi)
+\f{\tilde\z}{\r}(\f{n}{T}\bar\t_{x_1}^2+\bar n_{x_1}\bar\t_{x_1}-\f{n}{\r}\bar\r_{x_1}\bar\t_{x_1}+n\bar\t_{x_1x_1})\Big]dxd\tau\\ \nm
&+\int_0^t\int \Big[\bar u_{1x_1}\f{2\tilde\psi_1}{3\r}(T\pa_1\tilde\phi+n\pa_1\tilde\z-\f{nT}{\r}\pa_1\phi)
+\bar u_{1x_1}\f{n\tilde\z}{\r}(\f53\pa_1\tilde\psi_1-\f13\dv\tilde\psi)\\ \nm
&\qquad\qquad +\f{2\tilde\psi_1}{3\r}(T\bar n_{x_1}\bar u_{1x_1}+n\bar\t_{x_1}\bar u_{1x_1}-\f{nT}{\r}\bar\r_{x_1}\bar u_{1x_1}
+nT\bar u_{1x_1x_1})+\f{n\tilde\z}{\r}\f23\bar u_{1x_1}^2\Big]dxd\tau\\ \nm
&+\int_0^t\int \Big[\big(\int\f12|v|^2v\Pi dv\big)\cdot\na\big(\f{\tilde\z}{T}\big)
+\sum_{i,j=1}^3\big(\int v_iv_j\Pi dv\big)\pa_j\big(\f{\bar\t}{T}\tilde\psi_i\big)-\big(\int v_1 v\Pi dv\big)\cdot\na \big(\f{\tilde\z\bar u_1}{T}\big) \Big]dxd\tau\\ \nm
& :=I_1+I_2+I_3.
\end{align}
We estimate the right hand side of \eqref{E1} term by term. First, by Lemma \ref{appu}, the a priori assumption \eqref{assumption}, and Sobolev inequality, we have
\begin{align}\label{I1}
\nm &\quad\int_0^t\int\f{5\tilde\z}{6\r}(\f{n}{T}\bar\t_{x_1}^2+\bar n_{x_1}\bar\t_{x_1}-\f{n}{\r}\bar\r_{x_1}\bar\t_{x_1}+n\bar\t_{x_1x_1})dxd\tau\\
&\leq C\int_0^t\int_{\mathbb{T}^2}(\|\bar\t_{x_1x_1}\|_{L^1_{x_1}}+\|\bar\t_{x_1}\|_{L^2_{x_1}}^2)\|\tilde\z\|_{L^{\infty}_{x_1}}dx_2dx_3d\tau\\
\nm &\leq C\v^{\f18}\int_0^t\int_{\mathbb{T}^2}(1+\tau)^{-\f78}\|\tilde\z\|^{\f12}_{L_{x_1}^2}
\|\pa_1\tilde\z\|^{\f12}_{L_{x_1}^2}dx_2dx_3d\tau
 \leq C\v^{\f18}\int_0^t(1+\tau)^{-\f78}\|\tilde\z\|^{\f12}
\|\pa_1\tilde\z\|^{\f12}d\tau\\
\nm & \leq C\v^{\f18}\int_0^t(1+\tau)^{-\f78}\|\pa_1\tilde\z\|^{\f12}d\tau
\leq C\v^{\f18}\int_0^t\|\pa_1\tilde\z\|^2d\tau+C\v^{\f18}.
\end{align}
By Cauchy inequality and the above estimate, we have
\begin{align}\label{I2}
\nm &\int_0^t\int \f{5}{6}\bar\t_{x_1}\f{\tilde\z}{\r}(\pa_1\tilde\phi+\f{2n}{T}\pa_1\tilde\z-\f{n}{\r}\pa_1\phi) dxd\tau\\
&\leq \v^{\f{1}{8}}\int_0^t\|\pa_1(\phi,\tilde\phi,\tilde\z)\|^2d\tau+C\v^{-\f{1}{8}}\int_0^t\int\bar\t_{x_1}^{2}\tilde\z^2dxd\tau\\
\nm &\leq \v^{\f{1}{8}}\int_0^t\|\pa_1(\phi,\tilde\phi,\tilde\z)\|^2d\tau
+C\v^{-\f{1}{8}}\big(\v^{\f14}\int_0^t\|\pa_1\tilde\z\|^2d\tau+\v^{\f14}\int_0^t(1+\tau)^{-\f32}\|\tilde\z\|^2d\tau\big)\\
\nm &\leq \v^{\f{1}{8}}\int_0^t\|\pa_1(\phi,\tilde\phi,\tilde\z)\|^2d\tau+C\v^{\f{1}{8}}.
\end{align}
Similarly to \eqref{I1} and \eqref{I2}, we obtain
\be\label{I1I2}
I_1+I_2\leq C\v^{\f{1}{8}}\int_0^t\|\na(\phi,\tilde\phi,\tilde\psi,\tilde\z)\|^2d\tau+C\v^{\f{1}{8}}.
\ee
To handle $I_3$, we need the following estimate for the microscopic term. Indeed, from \eqref{Pi}, we have
$$
\begin{array}{ll}
\di \int v\otimes v\Pi dv=\int v\otimes v \Lop_{[w,T]}^{-1}[\r^{-1}G_t+\r^{-1}\mb {P}_1(v\cdot\na_xG)
-(w-u)\cdot\na_vG+(T-\t)\Delta_vG] dv.
\end{array}
$$
By Cauchy inequality and Corollary \ref{co-2}, we have
\be\label{L-1}
\ba
\int v\otimes v \Lop_{[w,T]}^{-1}\big(\r^{-1}G_t\big) dv &\leq C\Big(\int\f{1+|v|^2}{\Mloc_*}|\Lop_{[w,T]}^{-1}G_t|^2dv\Big)^{\f12}
\Big(\int(1+|v|^2)^{-1}|v|^4\Mloc_*dv\Big)^{\f12}\\
 &\leq C\Big(\int\f{(1+|v|^2)^{-1}}{\Mloc_*}|G_t|^2dv\Big)^{\f12}.
\ea
\ee
Likewise,
\be\label{L-2}
\ba
& \int v\otimes v \Lop_{[w,T]}^{-1}[\r^{-1}\mb{P}_1(v\cdot\na_xG)] dv
\leq C\Big(\int\f{1+|v|^2}{\Mloc_*}\big|\Lop_{[w,T]}^{-1}[\mb{P}_1(v\cdot\na_xG)]\big|^2dv\Big)^{\f12}\\
&\qquad\qquad\leq C\Big(\int\f{(1+|v|^2)^{-1}}{\Mloc_*}|\mb{P}_1(v\cdot\na_xG)|^2dv\Big)^{\f12}
 \leq C\Big(\int\f{|\na_xG|^2}{\Mloc_*}dv\Big)^{\f12}.
\ea
\ee
Since $w-u=(w-\bar u)-(u-\bar u)=\tilde\psi-\psi$, we have
\be\label{L-3}
\ba
 \int v\otimes v \Lop_{[w,T]}^{-1}[(w-u)\cdot\na_vG] dv
&\leq C|\psi-\tilde\psi|\Big(\int\f{(1+|v|^2)^{-1}}{\Mloc_*}|\na_vG|^2dv\Big)^{\f12}\\
& \leq C\d\Big(\int\f{|\na_v\wt{G}|^2}{\Mloc_*}dv\Big)^{\f12}+C|\psi-\tilde\psi|\bar u_{1x_1}.
\ea
\ee
It follows from \eqref{l2} that
\be\label{3.21}
\ba
\Lop^{-1}_{[w,T]}[(T-\t)\Delta_vG]&=\Lop^{-1}_{[w,T]}\big[\f{T-\t}{T}(\Lop_{[w,T]}G-3G-(v-w)\cdot\na_vG)\big]\\
&=\f{T-\t}{T}\big[G-3\Lop^{-1}_{[w,T]}G-\Lop^{-1}_{[w,T]}[(v-w)\cdot\na_vG]\big],
\ea
\ee
and since $T-\t=(T-\bar\t)-(\t-\bar\t)=\tilde\z-\z$, we have
\begin{align}\label{L-4}
\nm \int v\otimes v \Lop^{-1}_{[w,T]}[(T-\t)\Delta_vG]dv
&=\f{\tilde\z-\z}{T}\int v\otimes v\big[G-3\Lop^{-1}_{[w,T]}G-\Lop^{-1}_{[w,T]}[(v-w)\cdot\na_vG]\big]dv\\
&\leq C\d\Big(\int\f{|(\na_v\wt{G},\wt{G})|^2}{\Mloc_*}dv\Big)^{\f12}+C|\tilde\z-\z|\bar u_{1x_1}.
\end{align}
Combining \eqref{L-1}, \eqref{L-2}, \eqref{L-3} and \eqref{L-4} yields
\begin{align}\label{non-1}
\nm \int v\otimes v\Pi dv&\leq C\Big(\int\f{|(G_t,\na_xG)|^2}{\Mloc_*}dv\Big)^{\f12}
+C\d\Big(\int\f{|(\na_v\wt{G},\wt{G})|^2}{\Mloc_*}dv\Big)^{\f12}\\
&+C(|\tilde\psi-\psi|+|\tilde\z-\z|)\bar u_{1x_1},
\end{align}
which together with Cauchy inequality yields
\begin{align}\label{non-2}
\nm & I_3 \leq (\f18+C\d)\int_0^t\|\na(\tilde\psi,\tilde\z)\|^2d\tau+C\int_0^t\|\bar u_{1x_1}(\tilde\psi,\tilde\z)\|^2d\tau\\
&\qquad+C\int_0^t\int\Big[\big|\int \f{1}{2}|v|^2 v\Pi dv\big|^2+\big|\sum_{i,j=1}^3\int v_i v_j\Pi dv\big|^2\Big]dxd\tau\\
\nm &\quad \leq\big(\f18+C(\d+\v^{\f{1}{4}})\big)\int_0^t\|\na(\tilde\psi,\tilde\z)\|^2d\tau
+C\v^{\f14}+C\v\int_0^t(\|\tilde\psi-\psi\|^2+\|\tilde\z-\z\|^2)d\tau\\
\nm &\qquad +C\int_0^t\int\int\f{|\pa^{\alpha}G|^2}{\Mloc_*}dvdxd\tau
 +C\d\int_0^t\int\int\f{|(\na_v\wt{G},\wt{G})|^2}{\Mloc_*}dvdxd\tau.
\end{align}
By \eqref{rho-bd} and \eqref{t-bd}, there exists a constant $C_0>0$ independent of time such that
$$
 C_0^{-1}|(\phi,\psi,\z)|^2\leq \eta\leq  C_0|(\phi,\psi,\z)|^2,\quad
 C_0^{-1}|(\tilde\phi,\tilde\psi,\tilde\z)|^2\leq \tilde\eta\leq C_0|(\tilde\phi,\tilde\psi,\tilde\z)|^2.
$$
Note that integrating by parts over $\Omega=\mathbb{R}\times\mathbb{T}^2$ yields
$$
\ba
\int\f{n\bar\t}{2\r}
\Big(\f{|\na\tilde\psi+(\na\tilde\psi)^t|^2}{2}-\f23(\dv\tilde\psi)^2\Big)dx
&\geq C\int\Big( |\na\tilde\psi|^2+\na\tilde\psi\cdot(\na\tilde\psi)^t
-\f23(\dv\tilde\psi)^2\Big) dx\\
&=C\int|\na\tilde\psi|^2dx+C\int\f{(\dv\tilde\psi)^2}{3}dx\geq C\|\na\tilde\psi\|^2
\ea
$$
for some constant $C>0$.

Substituting \eqref{I1I2} and \eqref{non-2} into \eqref{E1}, and choosing $\d$ and $\v$ suitably small, we obtain
\begin{align}\label{low-1}\nm
& \|(\phi,\psi,\z,\tilde\phi,\tilde\psi,\tilde\z)(t)\|^2
+\int_0^t(\|\sqrt{\bar u_{1x_1}}(\phi,\psi_1,\z,\tilde\phi,\tilde\psi_1,\tilde\z)\|^2+\|\na(\tilde\psi,\tilde\z)\|^2)d\tau\\
&+\int_0^t(\|\tilde\psi-\psi\|^2+\|\tilde\z-\z\|^2)d\tau
\leq C\|(\phi_0,\psi_0,\z_0,\tilde\phi_0,\tilde\psi_0,\tilde\z_0)\|^2
  +C\v^{\f{1}{8}}\\ \nm
&+C\v^{\f{1}{8}}\int_0^t\|\na(\phi,\tilde\phi)\|^2d\tau
+C\int_0^t\int\int\f{|\pa^{\alpha}G|^2}{\Mloc_*}dvdxd\tau
+C\d\int_0^t\int\int\f{|(\na_v\wt{G},\wt{G})|^2}{\Mloc_*}dvdxd\tau.
\end{align}

It remains to control $\na\phi$ and $\na\tilde\phi$. To this end, we rewrite \eqref{perturb-0} as
\begin{align}\label{perturb-1}
\begin{cases}
& \phi_t+u\cdot\na\phi+\r\dv\psi+\psi\cdot\na\bar\r+\phi\dv\bar u=0,\\
&\psi_t+u\cdot\na\psi+\f{\t}{\r}\na\phi+\na\z+\psi\cdot\na\bar u+\f{1}{\r}\big(\z-\f{\bar\t}{\bar\r}\phi\big)\na\bar\r
=n(w-u),\\
&\f32\z_t+\f32u\cdot\na\z+\t\dv\psi+\f32\psi\cdot\na\bar\t+\z\dv\bar u=3n(T-\t)+n|w-u|^2,\\
&\tilde\phi_t+w\cdot\na\tilde\phi+ n\dv\tilde\psi+\tilde\psi\cdot\na\bar n+\tilde\phi\dv\bar u=0,\\
& \tilde\psi_t+w\cdot\na\tilde\psi+\f{T}{n}\na\tilde\phi+\na\tilde\z+\tilde\psi\cdot\na\bar u
+\f{1}{n}\big(\tilde\z-\f{\bar\t}{\bar n}\tilde\phi\big)\na\bar n\\
&\qquad\qquad\qquad=-\r(w-u)
-\f{1}{n}\int v\otimes v\cdot\na_xGdv,\\
&\f32\tilde\z_t+\f32w\cdot\na\tilde\z+T\dv\tilde\psi+\f32\tilde\psi\cdot\na\bar\t+\tilde\z\dv\bar u=-3\r(T-\t)\\
&\qquad\qquad\qquad+\f{1}{n}\Big(-\int \f12|v|^2v\cdot\na_xG dv+w\cdot\int v\otimes v\cdot\na_xG dv\Big).
\end{cases}
\end{align}
Inspired by \cite{WangWang2019}, we multiply \eqref{perturb-1}$_2$ by $\na\phi$
and \eqref{perturb-1}$_5$ by $\na\tilde\phi$, add the resulting equations, and integrate the sum by parts over $[0,t]\times\Omega$. This gives
\begin{align}\label{phi-1}
\nm &\quad \int(\psi\cdot\na\phi+\tilde\psi\cdot\na\tilde\phi)dx\Big|_{\tau=0}^{\tau=t}
+\int_0^t\int\big(\f{\t}{\r}|\na\phi|^2+\f{T}{n}|\na\tilde\phi|^2\big)dxd\tau\\
\nm &=\int_0^t\int\big[\dv\psi(u\cdot\na\phi+\r\dv\psi+\psi\cdot\na\bar\r+\phi\dv\bar u)
-(\psi\cdot\na\bar u+\f{1}{\r}\big(\z-\f{\bar\t}{\bar\r}\phi\big)\na\bar\r)\cdot\na\phi \big]dxd\tau\\
\nm &+\int_0^t\int\big[\dv\tilde\psi(w\cdot\na\tilde\phi+n\dv\tilde\psi+\tilde\psi\cdot\na\bar n+\tilde\phi\dv\bar u)
-(\tilde\psi\cdot\na\bar u++\f{1}{n}\big(\tilde\z-\f{\bar\t}{\bar n}\tilde\phi\big)\na\bar n)\cdot\na\tilde\phi \big]dxd\tau\\
&+\int\int\big[(w-u)\cdot\big(n \na\phi-\r \na\tilde\phi\big)-u\cdot\na\psi\cdot\na\phi-w\cdot\na\tilde\psi\cdot\na\tilde\phi\big]dxd\tau\\
\nm &-\int_0^t\int(\na\z\cdot\na\phi+\na\tilde\z\cdot\na\tilde\phi)dxd\tau
-\int_0^t\int\f{\na\tilde\phi}{n}\cdot(\int v\otimes v\cdot\na_xGdv) dxd\tau,
\end{align}
where by \eqref{perturb-1}$_1$ and \eqref{perturb-1}$_4$, we have used the following facts
\begin{equation*}
\ba
\int\psi_t\cdot\na\phi \,dx&=\f{d}{dt}\int \psi\cdot\na\phi \,dx-\int\psi\cdot\na\phi_t\,dx
=\f{d}{dt}\int \psi\cdot\na\phi \,dx+\int \dv\psi\,\phi_t \,dx\\
&=\f{d}{dt}\int \psi\cdot\na\phi \,dx-\int\dv\psi(u\cdot\na\phi+\r\dv\psi+\psi\cdot\na\bar\r+\phi\dv\bar u)dx
\ea
\end{equation*}
and
\begin{equation*}
\ba
\int\tilde\psi_t\cdot\na\tilde\phi \,dx&=\f{d}{dt}\int \tilde\psi\cdot\na\tilde\phi \,dx-\int\tilde\psi\cdot\na\tilde\phi_t\,dx
=\f{d}{dt}\int \tilde\psi\cdot\na\tilde\phi \,dx+\int \dv\tilde\psi\,\tilde\phi_t \,dx\\
&=\f{d}{dt}\int \tilde\psi\cdot\na\tilde\phi \,dx-\int\dv\tilde\psi(w\cdot\na\tilde\phi+n\dv\tilde\psi
+\tilde\psi\cdot\na\bar n+\tilde\phi\dv\bar u)dx.
\ea
\end{equation*}
By the Cauchy inequality and the following facts
\be\label{na-psi}
\left\{
\ba
&\int_0^t\|\na\psi\|^2d\tau\leq 2\int_0^t(\|\na(\tilde\psi-\psi)\|^2+\|\na\tilde\psi\|^2)d\tau
=2\int_0^t(\|\na(w-u)\|^2+\|\na\tilde\psi\|^2)d\tau,\\
& \int_0^t\|\na\z\|^2d\tau\leq 2\int_0^t(\|\na(\tilde\z-\z)\|^2+\|\na\tilde\z\|^2)d\tau
=2\int_0^t(\|\na(T-\t)\|^2+\|\na\tilde\z\|^2)dxd\tau,
\ea
\right.
\ee
from \eqref{phi-1} and \eqref{na-psi}, we have
\begin{equation}\label{phi-x}
\ba
\int_0^t\|\na(\phi,\tilde\phi)\|^2d\tau&\leq C\|(\psi_0,\na\phi_0,\tilde\psi_0,\na\tilde\phi_0)\|^2
+C\|(\psi,\tilde\psi)(t)\|^2+C\|\na(\phi,\tilde\phi)(t)\|^2\\
&+C\int_0^t(\|\na(\tilde\psi,\tilde\z)\|^2+\|w-u\|_{H^1}^2+\|\na(T-\t)\|^2)d\tau\\
&+C\v\int_0^t\|\sqrt{\bar u_{1x_1}}(\phi,\psi_1,\z,\tilde\phi,\tilde\psi_1,\tilde\z)\|^2d\tau
+C\int_0^t\int\int\f{|\na_xG|^2}{\Mloc_*}dvdxd\tau.
\ea
\end{equation}

We next establish estimates for the time derivatives of the macroscopic perturbation $(\phi,\psi,\z,\tilde\phi,\tilde\psi,\tilde\z)$. Multiplying \eqref{perturb-1}$_1$, \eqref{perturb-1}$_2$, \eqref{perturb-1}$_3$ by $\phi_t$, $\psi_t$ and $\z_t$, respectively,
then integrating the resulting equations by parts over $[0,t]\times\Omega$, we obtain
$$
\int_0^t\|(\phi_t,\psi_t,\z_t)\|^2 d\tau \leq C\int_0^t(\|\na(\phi,\psi,\z)\|^2+\|w-u\|^2+\|T-\t\|^2)d\tau
+C\int_0^t \|\bar u_{1x_1}(\phi,\psi_1,\z)\|^2d\tau,
$$
which together with \eqref{na-psi} yields
\begin{equation}\label{phi-t}
\ba
\int_0^t\|(\phi_t,\psi_t,\z_t)\|^2 d\tau &\leq C\int_0^t(\|\na(\phi,\tilde\psi,\tilde\z)\|^2+\|w-u\|_{H^1}^2+\|T-\t\|_{H^1}^2)d\tau\\
&+C\v\int_0^t\|\sqrt{\bar u_{1x_1}}(\phi,\psi_1,\z)\|^2d\tau.
\ea
\end{equation}
Similarly, multiplying \eqref{perturb-1}$_4$, \eqref{perturb-1}$_5$, \eqref{perturb-1}$_6$ by $\tilde\phi_t$, $\tilde\psi_t$, $\tilde\z_t$, respectively, and integrating the resulting equations by parts over $[0,t]\times\Omega$, we obtain
\begin{equation}\label{phi-t'}
\ba
\int_0^t\|(\tilde\phi_t,\tilde\psi_t,\tilde\z_t)\|^2 d\tau&
\leq C\int_0^t(\|\na(\tilde\phi,\tilde\psi,\tilde\z)\|^2+\|w-u\|^2+\|T-\t\|^2)d\tau\\
&+C\v\int_0^t\|\sqrt{\bar u_{1x_1}}(\tilde\phi,\tilde\psi_1,\tilde\z)\|^2d\tau
+C\int_0^t\int\int\f{|\na_xG|^2}{\Mloc_*}dvdxd\tau.
\ea
\end{equation}
Combining \eqref{low-1}, \eqref{phi-x}, \eqref{phi-t} and \eqref{phi-t'}  yields
\be\label{low-2}
\ba
&\|(\phi,\psi,\z,\tilde\phi,\tilde\psi,\tilde\z)(t)\|^2+\int_0^t(\|\sqrt{\bar u_{1x_1}}(\phi,\psi_1,\z,\tilde\phi,\tilde\psi_1,\tilde\z)\|^2
+\|\pa^{\alpha}(\phi,\tilde\phi,\tilde\psi,\tilde\z)\|^2+\|(\psi_t,\z_t)\|^2)d\tau\\
&+\int_0^t(\|w-u\|^2+\|T-\t\|^2)d\tau
 \leq C(\mathcal{N}(0)+\v^{\f{1}{8}})+C\|\na(\phi,\tilde\phi)(t)\|^2+C\int_0^t\|\na(w-u)\|^2d\tau
\\
&+C\int_0^t\|\na(T-\t)\|^2d\tau+C\d\int_0^t\int\int\f{|(\na_v\wt{G},\wt{G})|^2}{\Mloc_*}dvdxd\tau
+C\int_0^t\int\int\f{|\pa^{\alpha}G|^2}{\Mloc_*}dvdxd\tau.
\ea
\ee

It remains to estimate the microscopic terms on the right-hand side of \eqref{low-2}. Multiplying \eqref{tilde-G} by $\widetilde{G}/\mb{M_*}$, and integrating it over $[0,t]\times\Omega\times\mathbb{R}^3$, we have
\begin{align}\label{G0}
&\quad \nm \int\int\f{|\wt{G}|^2}{2\Mloc_*}dvdx\Big|_{\tau=0}^{\tau=t}
-\int_0^t\int\int\r\f{\wt{G}}{\Mloc_*}\Lop_{[w,T]}\wt{G}\,dvdxd\tau\\
&=\int_0^t\int\int\Big[-\f{1}{T}\mb{P}_1\big[v\cdot\big(\na\tilde\psi\cdot v+\f{|v-w|^2}{2T}\na\tilde\z\big)\Mloc\big]
-\mb{P}_1(v\cdot\na_xG)\\ \nm
&\qquad\qquad\qquad-\bar{\mathcal{G}}_t+\r(w-u)\cdot\na_vG-\r(T-\t)\Delta_vG\Big]\f{\wt{G}}{\Mloc_*}dvdxd\tau.
\end{align}
By Cauchy inequality and Lemma \ref{appu}, we have
\be\label{I4I5}
\ba
&\quad\int_0^t\int\int\big[-\f{1}{T}\mb{P}_1\big[v\cdot\big(\na\tilde\psi\cdot v+\f{|v-w|^2}{2T}\na\tilde\z\big)\Mloc\big]
-\mb{P}_1(v\cdot\na_xG)-\bar{\mathcal{G}}_t\big]\f{\wt{G}}{\Mloc_*}dvdxd\tau\\
& \leq \f{\bar\sigma}{8}\int_0^t\int\int\f{1+|v|^2}{\Mloc_*}|\wt{G}|^2dvdxd\tau
+C\int_0^t\|\na(\tilde\psi,\tilde\z)\|^2d\tau+C\int_0^t\int\int\f{|\na_xG|^2}{\Mloc_*}dvdxd\tau\\
&\quad+C\v\int_0^t\|(\tilde\phi_t,\tilde\psi_t,\tilde\z_t)\|^2d\tau+C\v^{\f12}
\ea
\ee
and
\be\label{I6}
\ba
&\quad \int_0^t\int\int \r (w-u)\cdot\na_vG\f{\wt{G}}{\Mloc_*}dvdxd\tau\\
 &\leq (\f{\bar\sigma}{8}+C\d)\int_0^t\int\int\f{|\na_v\wt{G}|^2+(1+|v|^2)|\wt{G}|^2}{\Mloc_*}dvdxd\tau
+C\v\int_0^t\|w-u\|^2d\tau.
\ea
\ee
Integrating by parts with respect to $v\in\mathbb{R}^3$, we obtain
\begin{align}\label{I6'}
\nm& -\int_0^t\int\int \r (T-\t)\Delta_vG\f{\wt{G}}{\Mloc_*}dvdxd\tau
=\int_0^t\int \r(T-\t)\int\na_vG\cdot\na_v\big(\f{\wt{G}}{\Mloc_*}\big)dv \,dxd\tau\\
&=\int_0^t\int \r(T-\t)\int(\na_v\wt{G}+\na_v\bar{\mathcal{G}})\cdot
\big(\f{\na_v\wt{G}}{\Mloc_*}+\f{v-w_*}{T_*}\f{\wt{G}}{\Mloc_*}\big)dv \,dxd\tau\\
\nm &\leq (\f{\bar\sigma}{8}+C\d)\int_0^t\int\int\f{|\na_v\wt{G}|^2+(1+|v|^2)|\wt{G}|^2}{\Mloc_*}dvdxd\tau
+C\v\int_0^t\|T-\t\|^2d\tau.
\end{align}
Substituting \eqref{I4I5}, \eqref{I6}, \eqref{I6'} into \eqref{G0}, and using Lemma \ref{co-1} with $\delta$ chosen sufficiently small, we obtain
\be\label{EG}
\begin{array}{l}
\di \int\int \f{|\widetilde{G}(t)|^2}{\Mloc_*}dxdv
 +\int_0^t\int\int\frac{|\na_v\wt{G}|^2+(1+|v|^2)|\wt{G}|^2}{\Mloc_*}dvdxd\tau\\
\di \leq C\int\int \f{|\widetilde{G}(0)|^2}{\Mloc_*}dxdv+C\v^{\f12}
+C\int_0^t\|\na(\tilde\psi,\tilde\z)\|^2d\tau\\
\di+C\v\int_0^t(\|(\tilde\phi_t,\tilde\psi_t,\tilde\z_t)\|^2+\|w-u\|^2+\|T-\t\|^2)d\tau
+C\int_0^t\int\int\f{|\na_xG|^2}{\Mloc_*}dvdxd\tau.
\end{array}
\ee
Multiplying \eqref{low-2} by a sufficiently large constant $C$, adding the resulting inequality to \eqref{EG}, and choosing $\d$ and $\v$ suitably small, we obtain the desired estimate \eqref{lower}. This completes the proof of Proposition \ref{Prop3.1}.

 \hfill $\Box$

%%%%%%%%%%%%%%%%%%%%%%%%%%%%%%%%%%%%%%%%%%%%%%%%%%%%%%%%%%%%%%%%%%%%%%%%%%%%%%%%%%%
\subsection{Higher-order estimate (I)}
\label{SSpf-prop-2}

In this subsection, we establish the first-order derivative estimates. To do this, we use the systems \eqref{twophase}$_{1,2,3}$ for $(\r,u,\t)$ and \eqref{mix-e} for $(n,w,T)$ with diffusion terms, in order to obtain estimates for $\int_0^t\|\na^2(u,\t)\|d\tau$ and $\int_0^t\|\na^2(w,T)\|d\tau$. We work with these systems rather than the perturbation system for
$(\phi,\psi,\z,\tilde\phi,\tilde\psi,\tilde\z)$ since the rarefaction wave itself has time decay rates
$$
\|(\bar\rho,\bar u_1,\bar\t)_{x_1x_1}\|^2\sim(1+t)^{-2},\quad
\|(\bar n,\bar u_1,\bar\t)_{x_1x_1}\|^2\sim(1+t)^{-2},
$$
which are integrable in time.

The proof of Proposition \ref{Prop3.2} consists of the following three lemmas.

\begin{lemma}\label{le3.1}
Under the same assumptions as in Proposition \ref{Prop3.2}, the following estimate holds for $t\in[0,\mathcal{T}]:$
\begin{align}\label{1-es}
\nm &\quad\|\na(\r,u,\t,n,w,T)(t)\|^2+\int_0^t(\|\na^2(w,T)\|^2+\|\na(w-u,T-\t)\|^2)d\tau\\
&\leq C(\mathcal{N}(0)+\v^{\f12})+C\int_0^t\int\int\f{|\na_x\pa^{\alpha}G|^2}{\Mloc_*}dvdxd\tau\\
\nm &+C(\d+\v)\int_0^t(\|\na(\phi,\tilde\phi,\tilde\psi,\tilde\z)\|^2+\|(\z_t,\tilde\phi_t,\tilde\z_t)\|^2
+\|\na^2 (\r,n)\|^2+\|(w-u,T-\t)\|^2)d\tau\\
\nm &+C(\d+\v)\int_0^t\int\int\Big(\f{|(\na_x\pa^{\alpha}G,\na_v\na_xG,\pa^{\alpha}G)|^2}{\Mloc_*}
+\f{|(\na_v\wt{G},\wt{G})|^2}{\Mloc_*}\Big)dvdxd\tau,
\end{align}
where $C>0$ is a generic constant independent of time $\mathcal{T}>0$, and $\pa^{\alpha}$ denotes $\pa_t$ or $\na_x$.
\end{lemma}

\noindent{{{\emph{Proof}}}}: By \eqref{rho-bd}, we rewrite the system \eqref{twophase} as
\be\label{1-ns}
\begin{cases}
\r_t+u\cdot\na\r+\r\dv u=0,\\
\di u_t+u\cdot\na u+\frac{\t}{\r}\na\r+\na\t=n(w-u),\\
\di \f32\t_t+\f32u\cdot\na\t+\t\dv u=3n(T-\t)+n|w-u|^2,
\end{cases}
\ee
and
\begin{align}\label{1-nw}
\begin{cases}
&n_t+w\cdot\na n+n\dv w=0,\\
&w_t+w\cdot\na w+\frac{T}{n}\na n+\na T=-\r (w-u)+\f{T}{2\r}(\Delta w+\f13\na\dv w)\\
&\qquad\qquad\qquad\qquad\qquad\qquad
+\na(\f{n T}{2\r})\cdot\f{\mathbb{S}(w)}{n}-\f{1}{n}\int v\otimes v\cdot\na_x\Pi dv,\\
&\f32T_t+\f32w\cdot\na T+T\dv w=-3\r(T-\t)+\f{5T}{6\r}\Delta T+\na(\f{5nT}{6\r})\cdot\f{\na T}{n}\\
&+\f{T}{\r}\big[\f{(\na w+(\na w)^t)^2}{4}-\f13(\dv w)^2\big]
+\f{1}{n}\big(-\int\f12|v|^2v\cdot\na_x\Pi dv+w\cdot\int v\otimes v\cdot\na_x \Pi dv\big).
\end{cases}
\end{align}
Applying $\pa_{i}:=\pa_{x_i}~(i=1,2,3)$ to system \eqref{1-ns} and \eqref{1-nw}, respectively, we obtain
\begin{align}\label{1-ns-x}
\begin{cases}
&\pa_i\r_t+u\cdot\na\pa_i\r+\r\pa_i\dv u+\pa_i u\cdot\na\r+\pa_i\r\dv u=0,\\
&\pa_i u_t+u\cdot\na\pa_i u+\frac{\t}{\r}\na\pa_i\r+\na\pa_i\t+\pa_i u\cdot\na u+\pa_i\Big(\f{\t}{\r}\Big)\na\r\\
& = n \pa_i(w-u)+\pa_i n (w-u),\\
& \f32\pa_i\t_t+\f32u\cdot\na\pa_i\t+\t\pa_i\dv u+\f32\pa_i u\cdot\na\t+\pa_i\t\dv u\\
&=3n\pa_i(T-\t)+3\pa_in(T-\t)+\pa_i n |w-u|^2+n\pa_i|w-u|^2,
\end{cases}
\end{align}
and
\begin{align}\label{1-nw-x}
\begin{cases}
&\pa_i n_t+w\cdot\na\pa_i n+n\pa_i\dv w+\pa_i w\cdot\na n+\pa_i n\dv w=0,\\
&\pa_i w_t+w\cdot\na\pa_i w+\frac{T}{n}\na\pa_i n+\na\pa_i T+\pa_i w\cdot\na w+\pa_i\Big(\f{T}{n}\Big)\na n\\
&=-\r\pa_i (w-u) -\pa_i\r (w-u)+\f{T}{2\r}(\Delta\pa_i w+\f13\na\pa_i\dv w)+\pa_i\big(\na(\f{n T}{2\r})\cdot\f{\mathbb{S}(w)}{n}\big)\\
&\quad+\pa_i (\f{T}{2\r}) (\Delta w+\f13\na\dv w)-\pa_i\Big(\f{1}{n}\int v\otimes v\cdot\na_x\Pi dv\Big),\\
&\f32\pa_i T_t+\f32w\cdot\na\pa_i T+T\pa_i\dv w+\f32\pa_i w\cdot\na T+\pa_i T\dv w\\
&=-3\r\pa_i(T-\t)-3\pa_i\r (T-\t)
+\f{5T}{6\r}\Delta\pa_i T\\
&\quad +\pa_i (\f{5T}{6\r})\Delta T+\pa_i\big(\na(\f{5nT}{6\r})\cdot\f{\na T}{n}\big)
+\pa_i\Big(\f{T}{\r}\big[\f{(\na w+(\na w)^t)^2}{4}-\f13(\dv w)^2\big]\Big)\\
&\quad +\pa_i\Big(\f{1}{n}\big(-\int\f12|v|^2v\cdot\na_x\Pi dv+w\cdot\int v\otimes v\cdot\na_x \Pi dv\big)\Big).
\end{cases}
\end{align}
Multiplying $\eqref{1-ns-x}_1$, $\eqref{1-ns-x}_2$, and $\eqref{1-ns-x}_3$ by $\f{\t}{\r n}\pa_i\r$, $\f{\r}{n}\pa_i u$, and
$\f{\r}{\t n}\pa_i \t$, respectively, then adding the resulting equations and integrating by parts over $[0,t]\times\Omega$, we obtain
\begin{align}\label{1-ns-ji} \nm
\nm&\int\big[\f{\t}{\r n}\f{(\pa_i\r)^2}{2}+\f{\r}{n}\f{|\pa_i u|^2}{2}+\f{3\r}{2\t n}\f{(\pa_i\t)^2}{2}\big]dx\Big|_{\tau=0}^{\tau=t}
-\int_0^t\int\r\big[\pa_i(w-u)\cdot\pa_i u +\f{3}{\t}\pa_i(T-\t)\pa_i\t\big]dxd\tau\\
\nm&=\int_0^t\int \Big[(\t_t+u\cdot\na\t)\big[\f{1}{\r n}\f{(\pa_i\r)^2}{2}-\f{3\r}{2\t^2n}\f{(\pa_i\t)^2}{2}\big]
-(n_t+u\cdot\na n)\big[\f{\t}{\r n^2}\f{(\pa_i\r)^2}{2}+\f{\r}{n^2}\f{|\pa_i u|^2}{2}\\
\nm&+\f{3\r}{2\t n^2}\f{(\pa_i\t)^2}{2}\big]+\pa_i\r\pa_i u\cdot\na(\f{\t}{n})+\pa_i\t\pa_i u\cdot\na(\f{\r}{n})
-\f{\r}{n}\pa_i u\cdot\na u\cdot\pa_i u-\f{\r}{n}\pa_i(\f{\t}{\r})\pa_i u\cdot\na\r\\
&+\f{\r}{n}\pa_i n\pa_i u\cdot(w-u)
+\big[3\pa_i n (T-\t)+\pa_i n |w-u|^2+ n \pa_i|w-u|^2\big]\f{\r}{\t n}\pa_i\t\\
\nm &-(\f32\pa_i u\cdot\na\t+\pa_i\t\dv u)\f{\r}{n\t}\pa_i\t-\pa_iu\cdot\na\r\f{\t}{\r n}\pa_i\r\Big]dxd\tau=:I_4,
\end{align}
where we have used the following facts,
\begin{align} \nm
&\quad\int(\pa_i\r_t+u\cdot\na\pa_i\r+\r\pa_i\dv u+\pa_i\r\dv u)\f{\t}{\r n}\pa_i\r\,dx
\\ \nm
&=\f{d}{dt}\int\f{\t}{\r n}\f{(\pa_i\r)^2}{2}dx-\int\big[\t_t+u\cdot\na\t-\f{\t}{n}(n_t+u\cdot\na n)\big]\f{(\pa_i\r)^2}{2\r n}dx\\ \nm
&\quad-\int\f{\t}{n}\pa_i u\cdot\na\pa_i\r dx-\int\pa_i\r \pa_i u\cdot\na(\f{\t}{n})dx,
\end{align}
\begin{align}\nm
&\quad\int(\pa_i u_t+u\cdot\na\pa_i u+\f{\t}{\r}\na\pa_i\r+\na\pa_i\t)\cdot\f{\r}{n}\pa_i u\,dx\\ \nm
&=\f{d}{dt}\int\f{\r}{n}\f{|\pa_i u|^2}{2}dx+\int(n_t+u\cdot\na n)\f{\r|\pa_i u|^2}{2n^2}dx
+\int\f{\t}{n}\pa_i u\cdot\na\pa_i\r \, dx\\ \nm
& \quad -\int\f{\r}{n}\pa_i\t\pa_i\dv u\,dx-\int\pa_i\t\pa_i u\cdot\na(\f{\r}{n})dx,
\end{align}
and
\begin{align} \nm
&\quad\int(\f32\pa_i\t_t+\f32u\cdot\na\pa_i\t+\t\pa_i\dv u)\f{\r\pa_i\t}{\t n}\,dx
\\ \nm
&=\f{d}{dt}\int\f{3\r}{2\t n}\f{(\pa_i\t)^2}{2}dx+\f32\int\big[\f{\r}{\t^2 n}(\t_t+u\cdot\na\t)+\f{\r}{\t n^2}(n_t+u\cdot\na n)\big]\f{(\pa_i\t)^2}{2}dx\\ \nm
&\quad +\int\f{\r}{n}\pa_i \t\pa_i\dv u \,dx.
\end{align}
Similarly, multiplying $\eqref{1-nw-x}_1$, $\eqref{1-nw-x}_2$, and $\eqref{1-nw-x}_3$ by $\f{T}{n^2}\pa_i n$, $\pa_i w$, and $\f{\pa_i T}{T}$, respectively, then adding the resulting equations and integrating by parts over $[0,t]\times\Omega$, we obtain
\begin{align}\label{1-nw-ji} \nm
& \int\big[\f{T}{n^2}\f{(\pa_i n)^2}{2}+\f{|\pa_i w|^2}{2}+\f{3}{2T}\f{(\pa_i T)^2}{2}\big]dx\Big|_{\tau=0}^{\tau=t}
+\int_0^t\int\r\big[\pa_i(w-u)\cdot\pa_i w+\f{3}{T}\pa_i(T-\t)\pa_i T\big] dxd\tau\\ \nm
&\quad +\int_0^t\int\big(\f{T}{2\r}\big(|\na\pa_i w|^2+\f{(\pa_i\dv w)^2}{3}\big)+\f{5}{6\r}|\na\pa_i T|^2\big)dxd\tau\\ \nm
& =\int_0^t\int \Big\{\dv w\big(\f{T}{n^2}\f{(\pa_i n)^2}{2}+\f{|\pa_i w|^2}{2}+\f{3}{2T}\f{(\pa_i T)^2}{2}\big)
+(T_t+w\cdot\na T)\big[\f{(\pa_i n)^2}{2n^2}-\f{3(\pa_i T)^2}{4T^2}\big]\\ \nm
&\quad  +\pa_i w\cdot\Big[\na(\f{T}{n})\pa_i n-\f{T}{n^2}\pa_i n\na n-\pa_i w\cdot\na w
-\pa_i(\f{T}{n})\na n-\na (\f{T}{2\r})\cdot\na\pa_i w-\f{1}{3}\pa_i\dv w\na (\f{T}{2\r})
\\
&\quad +\pa_i (\f{T}{2\r})(\Delta w+\f{\na\dv w}{3})
+\pa_i\big(\na(\f{ nT}{2\r})\cdot\f{\mathbb{S}(w)}{n}\big)-\pa_i\r(w-u)
\Big]-\na\pa_i T\cdot\na(\f{5}{6\r})\pa_i T\\ \nm
&
\quad +\Big[-\f32\pa_iw\cdot\na T-\pa_iT\dv w-3\pa_i\r(T-\t)
+\pa_i(\f{5T}{6\r})\Delta T+\pa_i\big(\na(\f{5nT}{6\r})\cdot\f{\na T}{n}\big)\\ \nm
&\quad +\pa_i \big(\f{T}{\r}\big(\f{(\na w+(\na w)^t)^2}{4}-\f{(\dv w)^2}{3}\big)\big)\Big]\f{\pa_i T}{T}\Big\}dxd\tau
+\int_0^t\int\Big[\f{\pa_i\pa_i w}{n}\cdot\int v\otimes v\cdot\na_x \Pi dv\\ \nm
&
\quad +\pa_i(\f{\pa_i T}{T})\f{1}{n}\big(\int\f12|v|^2v\cdot\na_x\Pi dv-w\cdot\int v\otimes v\cdot\na_x\Pi dv\big)
  \Big]dxd\tau=:I_5+I_6,
\end{align}
where we have used the following facts,
\begin{align} \nm
&\quad\int(\pa_i n_t+w\cdot\na\pa_i n+\pa_i n\dv w+n\pa_i\dv w)\f{T}{n^2}\pa_i n\,dx
\\ \nm
&=\f{d}{dt}\int\f{T}{n^2}\f{(\pa_i n)^2}{2}dx-\int\f{(\pa_i n)^2}{2n^2}\big(T_t+w\cdot\na T+T\dv w\big)dx \\ \nm
& \quad -\int\pa_i w\cdot\na(\f{T}{n})\pa_i n \,dx-\int \f{T}{n}\pa_i w\cdot\na\pa_i n\, dx,
\end{align}
\begin{align} \nm
&\quad \int(\pa_i w_t+w\cdot\na\pa_i w+\f{T}{n}\na\pa_i n+\na\pa_i T)\cdot\pa_i w dx\\ \nm
&=\f{d}{dt}\int \f{|\pa_i w|^2}{2} dx-\int\dv w\f{|\pa_i w|^2}{2}dx
+\int\f{T}{n}\pa_i w\cdot\na\pa_i n \, dx-\int\pa_i T\pa_i \dv w \,dx,
\end{align}
and
\begin{align} \nm
&\quad \int(\f32\pa_i T_t+\f32w\cdot\na\pa_i T+T\pa_i \dv w)\f{\pa_i T}{T}\,dx\\ \nm
&=\f{d}{dt}\int \f{3}{2T}\f{(\pa_i T)^2}{2}dx-\f32\int\f{(\pa_i T)^2}{2}\Big(\f{\dv w}{T}-\f{T_t+w\cdot\na T}{T^2}\Big)dx
+\int \pa_i\dv w\pa_i T\, dx.
\end{align}
Adding \eqref{1-ns-ji} and \eqref{1-nw-ji}, and summing the resulting equation over $i=1,2,3$, we obtain
\begin{align}\label{1-long} \nm
& \int\big(\f{\t}{\r n}\f{|\na\r|^2}{2}+\f{\r}{n}\f{|\na u|^2}{2}+\f{3\r}{2\t n}\f{|\na\t|^2}{2}
+\f{T}{n^2}\f{|\na n|^2}{2}+\f{|\na w|^2}{2}+\f{3}{2T}\f{(\na T)^2}{2}\big)dx\Big|_{\tau=0}^{\tau=t}\\
&\quad +\int_0^t\int\big(\f{T}{2\r}\big(|\na^2 w|^2+\f{|\na\dv w|^2}{3}\big)+\f{5}{6\r}|\na^2 T|^2
+\r|\na(w-u)|^2+\f{3\r}{T}|\na(T-\t)|^2\big)dxd\tau\\ \nm
&=\sum_{i=1}^3\int_0^t\int\f{3\r}{\t T}(T-\t)\pa_i\t\pa_i(T-\t)dxd\tau+I_4+I_5+I_6.
\end{align}
 The first term on the right-hand side of \eqref{1-long} is absorbed as follows. By Cauchy's inequality,
$$
\f{3\r}{\t T}|T-\t||\pa_i\t||\pa_i(T-\t)|
\leq \f{\r}{2T}|\pa_i(T-\t)|^2+C|T-\t|^2|\pa_i\t|^2,
$$
and, since $\|\na\t\|_{L^\infty}\leq C(\d+\v)$ by \eqref{assum-1},
$$
\sum_{i=1}^3\int_0^t\int\f{3\r}{\t T}(T-\t)\pa_i\t\pa_i(T-\t)dxd\tau
\leq \f12\int_0^t\int\f{3\r}{T}|\na(T-\t)|^2dxd\tau
+C(\d+\eps)^2\int_0^t\|T-\t\|^2d\tau,
$$
where the first term on the right is absorbed by the corresponding term on the
left-hand side of \eqref{1-long}, and the second is of the form already present
on the right-hand side of \eqref{1-es}.
The estimates for $I_4$, $I_5$ and $I_6$ are similar to those in Lemma 4.1 in \cite{LiWangWang2022}. Specifically, we have
$$
\ba
I_4&\leq C\int_0^t\int(|(\t_t,n_t)|+|\na(u,\t,n)|)|\na(\r,u,\t)|^2dxd\tau\\
&\quad +C\int_0^t\int(|w-u|+|T-\t|)|\na(\r,n)||\na(u,\t)|dxd\tau\\
&\quad +C\int_0^t\int|w-u||\na(w-u)||\na\t|dxd\tau=:I_{4,1}+I_{4,2}+I_{4,3}.
\ea
$$
By Lemma \ref{appu} and the assumption \eqref{assum-1}, we obtain
\begin{align} \nm
I_{4,1}&\leq C\int_0^t\int(|(\z_t,\tilde\phi_t)|+|\na(\psi,\z,\tilde\phi)|)|\na(\phi,\psi,\z)|^2dxd\tau
+C\int_0^t\int\bar u_{1x_1}|\na(\phi,\psi,\z)|^2dxd\tau\\ \nm
&\quad+C\int_0^t\int(|(\z_t,\tilde\phi_t)|+|\na(\psi,\z,\tilde\phi)|)\bar u_{1x_1}^2 dxd\tau
+C\int_0^t\int|\bar u_{1x_1}|^3dxd\tau\\ \nm
&\leq C\int_0^t\int\|\na(\phi,\psi,\z)\|_{L^{\infty}}(\|(\z_t,\tilde\phi_t)\|+\|\na(\psi,\z,\tilde\phi)\|) \|\na(\phi,\psi,\z)\|dxd\tau\\ \nm
&\quad+C\int_0^t\|\bar u_{1x_1}\|_{L^{\infty}_{x_1}}(\|(\z_t,\tilde\phi_t)\|
+\|\na(\psi,\z,\tilde\phi)\|+\|\na(\phi,\psi,\z)\|)d\tau+C\int_0^t\|\bar u_{1x_1}\|^3_{L^3_{x_1}}d\tau\\ \nm
&\leq C(\d+\v)\int_0^t(\|(\z_t,\tilde\phi_t)\|^2+\|\na(\phi,\psi,\z,\tilde\phi)\|^2)d\tau
+C\v^{\f23}\int_0^t(1+\tau)^{-\f43}d\tau\\ \nm
&\leq C(\d+\v)\int_0^t(\|(\z_t,\tilde\phi_t)\|^2+\|\na(\phi,\tilde\phi,\tilde\psi,\tilde\z)\|^2
+\|\na(w-u,T-\t)\|^2)d\tau+C\v^{\f23}.
\end{align}
Similarly, we have
\begin{align}\nm
I_{4,2}&\leq C\int_0^t\int|(w-u,T-\t)|(|\na(\phi,\tilde\phi)||\na(\psi,\z)|
+\bar u_{1x_1}|\na(\phi,\psi,\z,\tilde\phi)|+\bar u_{1x_1}^2) dxd\tau\\ \nm
&\leq C\int_0^t\|w-u,T-\t\|_{L^{\infty}}\|\na(\phi,\tilde\phi)\|\|\na(\psi,\z)\|d\tau\\ \nm
&\quad+C\int_0^t\|\bar u_{1x_1}\|_{L^{\infty}_{x_1}}\|(w-u,T-\t)\|\|\na(\phi,\psi,\z,\tilde\phi)\|d\tau\\ \nm
&\quad+C\int_0^t\|\bar u_{1x_1}\|_{L^{\infty}_{x_1}}\|(w-u,T-\t)\|^2d\tau
+C\int_0^t\|\bar u_{1x_1}\|^3_{L^3_{x_1}}d\tau\\ \nm
&\leq C(\d+\v)\int_0^t(\|\na(\phi,\psi,\z,\tilde\phi)\|^2+\|(w-u,T-\t)\|^2)d\tau
+C\v^{\f23}\int_0^t(1+\tau)^{-\f43}d\tau\\ \nm
&\leq C(\d+\v)\int_0^t(\|\na(\phi,\tilde\phi,\tilde\psi,\tilde\z)\|^2
+\|(w-u,T-\t)\|_{H^1}^2)d\tau+C\v^{\f23},
\end{align}
and
$$
I_{4,3}\leq C\int_0^t\|\na\t\|_{L^{\infty}}\|w-u\|\|\na(w-u)\|d\tau
\leq C(\d+\v)\int_0^t\|w-u\|^2_{H^1}d\tau.
$$
Thus, we obtain
\be\label{I4}
I_4\leq C(\d+\v)\int_0^t(\|(\z_t,\tilde\phi_t)\|^2+\|\na(\phi,\tilde\phi,\tilde\psi,\tilde\z)\|^2
+\|(w-u,T-\t)\|_{H^1}^2)d\tau+C\v^{\f23}.
\ee
For $I_5$, we have
$$
\ba
I_5&\leq C\int_0^t\int(|T_t|+|\na(w,T)|)|\na(n,w,T)|^2dxd\tau\\
&\quad +\int_0^t\int|\na(\r,n,w,T)|^2(|\na^2(\r,n,w,T)|+|w-u|)dxd\tau
=:I_{5,1}+I_{5,2}.
\ea
$$
Similarly to $I_{4,1}$, we have
$$
I_{5,1}\leq C(\d+\v)\int_0^t(\|\tilde\z_t\|^2+\|\na(\tilde\phi,\tilde\psi,\tilde\z)\|^2)d\tau
+C\v^{\f23},
$$
and similarly to $I_{4,2}$, we have
$$
\ba
I_{5,2}&\leq \int_0^t\int(|\na(\phi,\tilde\phi,\tilde\psi,\tilde\z)|^2+\bar u_{1x_1}^2)(|\na^2(\r,n,w,T)|+|w-u|)dxd\tau\\
&\leq C\int_0^t\|\na(\phi,\tilde\phi,\tilde\psi,\tilde\z)\|_{L^{\infty}_{x_1}}
\|\na(\phi,\tilde\phi,\tilde\psi,\tilde\z)\| (\|\na^2(\r,n,w,T)\|+\|w-u\|)d\tau\\
&\quad+C\int_0^t\|\bar u_{1x_1}\|_{L^{\infty}_{x_1}}(\|\na^2(\r,n,w,T)\|^2+\|w-u\|^2)d\tau
+C\int_0^t\|\bar u_{1x_1}\|^3_{L^3_{x_1}}d\tau\\
&\leq C(\d+\v)\int_0^t(\|\na(\phi,\tilde\phi,\tilde\psi,\tilde\z)\|^2+\|\na^2(\r,n,w,T)\|^2)d\tau
+C\v^{\f23}.
\ea
$$
Therefore, we obtain
\be\label{I5}
I_5\leq C(\d+\v)\int_0^t(\|\tilde\z_t\|^2+\|\na(\phi,\tilde\phi,\tilde\psi,\tilde\z)\|^2
+\|\na^2(\r,n,w,T)\|^2+\|w-u\|^2)d\tau+C\v^{\f23}.
\ee
It remains to handle the last term $I_6$ on the right-hand side of \eqref{1-long}. By Lemma \ref{ex}, we have
\be\label{1-Gt}
\ba
& \int v\otimes v\cdot\na_x \Lop^{-1}_{[w,T]}(\r^{-1}G_t) dv\\
=&\int v\otimes v\cdot\big[\Lop^{-1}_{[w,T]}\na_x(\r^{-1}G_t)
+\na w\cdot\Lop_{[w,T]}^{-1}[\na_v\Lop^{-1}_{[w,T]}(\r^{-1}G_t)]\big]dv\\
&\quad  -\f{\na T}{T}\int v\otimes v\cdot
\Lop^{-1}_{[w,T]}[\r^{-1}G_t-3\Lop^{-1}_{[w,T]}(\r^{-1}G_t)-(v-w)\cdot\na_v\Lop^{-1}_{[w,T]}(\r^{-1}G_t)] dv.
\ea
\ee
By the assumption \eqref{assumption} and Corollary \ref{co-2}, we have
\be\label{3.60}
\ba
&\quad \int v\otimes v\cdot\big[\Lop^{-1}_{[w,T]}\na_x(\r^{-1}G_t)
+\na w\cdot\Lop_{[w,T]}^{-1}[\na_v\Lop^{-1}_{[w,T]}(\r^{-1}G_t)]\big]dv\\
&= \int v\otimes v\cdot\big[\r^{-1}\Lop^{-1}_{[w,T]}(\na_xG_t)
+\na(\r^{-1})\Lop^{-1}_{[w,T]}G_t
+\r^{-1}\na w\cdot\Lop_{[w,T]}^{-1}(\na_v\Lop^{-1}_{[w,T]}G_t)\big]dv\\
& \leq C\Big(\int\f{1+|v|^2}{\Mloc_*}|\Lop^{-1}_{[w,T]}(\na_xG_t)|^2dv\Big)^{\f12}
+C|\na\r|\Big(\int\f{1+|v|^2}{\Mloc_*}|\Lop^{-1}_{[w,T]}G_t|^2dv\Big)^{\f12}\\
&\quad+ C|\na w|\Big(\int\f{1+|v|^2}{\Mloc_*}|\Lop^{-1}_{[w,T]}[\na_v(\Lop_{[w,T]}^{-1}G_t)]|^2dv\Big)^{\f12}\\
&\leq C\Big(\int\f{(1+|v|^2)^{-1}}{\Mloc_*}|\na_xG_t|^2dv\Big)^{\f12}
+C(\d+\v)\Big(\int\f{(1+|v|^2)^{-1}}{\Mloc_*}|G_t|^2dv\Big)^{\f12},
\ea
\ee
and
\be\label{3.61}
\ba
&\quad\f{\na T}{T}\int v\otimes v\cdot
\Lop^{-1}_{[w,T]}[\r^{-1}G_t-3\Lop^{-1}_{[w,T]}(\r^{-1}G_t)-(v-w)\cdot\na_v\Lop^{-1}_{[w,T]}(\r^{-1}G_t)] dv\\
& \leq C|\na T|\Big(\int\f{1+|v|^2}{\Mloc_*}|\Lop^{-1}_{[w,T]}
[G_t-3\Lop^{-1}_{[w,T]}G_t-(v-w)\cdot\na_v(\Lop^{-1}_{[w,T]}G_t)]|^2dv\Big)^{\f12}\\
&\leq C(\d+\v)\Big(\int\f{(1+|v|^2)^{-1}}{\Mloc_*}|G_t|^2dv\Big)^{\f12}.
\ea
\ee
Substituting \eqref{3.60} and \eqref{3.61} into \eqref{1-Gt}, we obtain
\begin{align}\label{1-Gt'}
 \int v\otimes v\cdot\na_x \Lop^{-1}_{[w,T]}(\r^{-1}G_t) dv
&\leq C\Big(\int\f{(1+|v|^2)^{-1}}{\Mloc_*}|\na_xG_t|^2dv\Big)^{\f12}\\
\nm &+C(\d+\v)\Big(\int\f{(1+|v|^2)^{-1}}{\Mloc_*}|G_t|^2dv\Big)^{\f12}.
\end{align}
Likewise,
\be\label{1-Gx}
\int v\otimes v\cdot\na_x \Lop^{-1}_w[\r^{-1}\mb{P}_1(v\cdot\na_xG)] dv
\leq C\Big(\int\f{|\na^2_xG|^2}{\Mloc_*}dv\Big)^{\f12}
+C(\d+\v)\Big(\int\f{|\na_xG|^2}{\Mloc_*}dv\Big)^{\f12}.
\ee
By Lemma \ref{ex} and the assumption \eqref{assum-1}, we have
\begin{align}\label{1-Gv}
\nm &\quad \int v\otimes v\cdot\na_x \Lop_{[w,T]}^{-1}[(w-u)\cdot\na_vG]dv\\
\nm &=\int v\otimes v\cdot \Lop_{[w,T]}^{-1}[\na(w-u)\cdot\na_vG+(w-u)\cdot\na_v\na_xG]dv\\
&\quad+\int v\otimes v\cdot\na w\cdot\Lop_{[w,T]}^{-1}[\na_v\Lop_{[w,T]}^{-1}((w-u)\cdot\na_vG)]dv\\
\nm &\quad-\int v\otimes v\cdot\f{\na T}{T}\Lop^{-1}_{[w,T]}\big[
(w-u)\cdot\na_vG-3\Lop^{-1}_{[w,T]}[(w-u)\cdot\na_vG]\\
\nm &\qquad-(v-w)\cdot\na_v\Lop^{-1}_{[w,T]}[(w-u)\cdot\na_vG]\big]dv\\
\nm &\leq C(|\na(w-u)|+|\na (w,T)||w-u|)\Big(\int\f{|\na_vG|^2}{\Mloc_*}dv\Big)^{\f12}
+C|w-u|\Big(\int\f{|\na_v\na_xG|^2}{\Mloc_*}dv\Big)^{\f12}\\
\nm &\leq C\d\Big(\int\f{|(\na_v\na_xG,\na_v\wt{G})|^2}{\Mloc_*}dv\Big)^{\f12}
+C(|\na(w-u)|+|\na (w,T)||w-u|)\bar u_{1x_1}.
\end{align}
From \eqref{3.21}, we deduce that
\be\label{3.67}
\ba
 \na_x\Lop^{-1}_{[w,T]}[(T-\t)\Delta_vG]
&=\na_x(\f{T-\t}{T})\big(G-3\Lop^{-1}_{[w,T]}G-\Lop^{-1}_{[w,T]}[(v-w)\cdot\na_vG]\big)\\
&\quad +\f{T-\t}{T}\big(\na_xG-3\na_x\Lop^{-1}_{[w,T]}G-\na_x\Lop^{-1}_{[w,T]}[(v-w)\cdot\na_vG]\big),
\ea
\ee
where it follows from Lemma \ref{ex} that
\be\label{3.68}
\ba
\na_x\Lop^{-1}_{[w,T]}G&=\Lop^{-1}_{[w,T]}(\na_xG)+\na w\cdot\Lop^{-1}_{[w,T]}[\na_v(\Lop^{-1}_{[w,T]}G)]\\
&\quad -\f{\na T}{T}\Lop^{-1}_{[w,T]}[G-3\Lop^{-1}_{[w,T]}G-(v-w)\cdot\na_v(\Lop^{-1}_{[w,T]}G)]
\ea
\ee
and
\begin{align}\label{3.69}
&\quad\na_x\Lop^{-1}_{[w,T]}[(v-w)\cdot\na_vG]\\
\nm &=\Lop^{-1}_{[w,T]}[(v-w)\cdot\na_v\na_xG-\na w\cdot\na_vG]
+\na w\cdot\Lop^{-1}_{[w,T]}\big[\na_v\Lop^{-1}_{[w,T]}[(v-w)\cdot\na_vG]\big]\\
\nm &\quad -\f{\na T}{T}\Lop^{-1}_{[w,T]}\big[(v-w)\cdot\na_vG-3\Lop^{-1}_{[w,T]}[(v-w)\cdot\na_vG]
-(v-w)\cdot\na_v\Lop^{-1}_{[w,T]}[(v-w)\cdot\na_vG]\big].
\end{align}
Substituting \eqref{3.68} and \eqref{3.69} into \eqref{3.67} and using Cauchy inequality together with Corollary \ref{co-2}, we obtain
\begin{align}\label{3.70}
&\quad \int v\otimes v\cdot\na_x\Lop^{-1}_{[w,T]}[(T-\t)\Delta_vG]dv\\
\nm&\leq C|T-\t|\Big(\int\f{|(\na_v\na_xG,\na_xG)|^2}{\Mloc_*} dv\Big)^{\f12} +C\big(|\na(T-\t)|+|\na(w,T)||T-\t|\big)\Big(\int\f{|(\na_vG,G)|^2}{\Mloc_*} dv\Big)^{\f12}\\
\nm & \leq C\d\Big(\int\f{|(\na_v\na_xG,\na_xG,\na_v\wt{G},\wt{G})|^2}{\Mloc_*} dv\Big)^{\f12}+C\big(|\na(T-\t)|+|\na(w,T)||T-\t|)\bar u_{1x_1}.
\end{align}
Combining \eqref{1-Gt'}, \eqref{1-Gx}, \eqref{1-Gv} and \eqref{3.70} gives
\begin{align}\label{Pi-x}
\nm & \int v\otimes v\cdot\na_x\Pi dv
 \leq C\Big(\int\f{|\na_x\pa^{\alpha}G|^2}{\Mloc_*}dv\Big)^{\f12}
+C(\d+\v)\Big(\int\f{|(\na_v\na_xG,\pa^{\alpha}G,\na_v\wt{G},\wt{G})|^2}{\Mloc_*}dv\Big)^{\f12}\\
&\quad+C(|\na(w-u)|+|\na(w,T)||w-u)|)\bar u_{1x_1}+C(|\na(T-\t)|+|\na(w,T)||T-\t|)\bar u_{x_1},
\end{align}
which together with Cauchy inequality and the assumption \eqref{assumption} yields
\begin{align}\label{es-Pi-x}
\nm I_6&\leq \f14\int_0^t\|\na^2(w,T)\|^2d\tau+C\int_0^t\|\na T\|^4_{L^4}d\tau\\
&\quad+C\int_0^t\int\big(\big|\int \f12|v|^2v\cdot\na_x\Pi dv\big|^2+\big|\int v\otimes v\cdot\na_x\Pi dv\big|^2\big)dxd\tau\\
\nm &\leq \big(\f14+C(\d+\v)\big)\int_0^t\|\na^2(w,T)\|^2d\tau+C\d\int_0^t\|\na\tilde\z\|^2d\tau
+C\v\int_0^t\|(w-u,T-\t)\|_{H^1}^2d\tau+C\v^{\f12}\\
\nm &\quad+C\int_0^t\int\int\f{|\na_x\pa^{\alpha}G|^2}{\Mloc_*}dvdxd\tau
+C(\d+\v)\int_0^t\int\int\f{|(\na_v\na_xG,\pa^{\alpha}G,\na_v\wt{G},\wt{G})|^2}{\Mloc_*}dvdxd\tau.
\end{align}
Substituting \eqref{I4}, \eqref{I5} and \eqref{es-Pi-x} into \eqref{1-long} yields the desired inequality \eqref{1-es}. The proof of Lemma \ref{le3.1} is complete.

\hfill $\Box$

\begin{lemma}\label{le3.2}
Under the same assumptions as in Proposition \ref{Prop3.2}, the following estimate holds for $t\in[0,\mathcal{T}]:$
\be\label{rho-2x}
\begin{array}{l}
\di \int_0^t\|\na^2(\r,n)\|^2d\tau \leq C(\mathcal{N}(0)+\v^{\f12})
+C\|\na(u,w)(t)\|^2+C\|\na^2(\r,n)(t)\|^2
\\
\di\quad+C\int_0^t(\|\na^2(w,T)\|^2+\|\na(w-u)\|_{H^1}^2+\|\na^2(T-\t)\|^2)d\tau\\
\di\quad+C(\d+\v)\int_0^t(\|\na(\phi,\tilde\phi,\tilde\psi,\tilde\z)\|^2+\|w-u\|^2+\|\na(T-\t)\|^2)d\tau\\
\di\quad+C\int_0^t\int\int\f{|\na_x^2G|^2}{\Mloc_*}dvdxd\tau
+C(\d+\v)\int_0^t\int\int\f{|\na_xG|^2}{\mb{M_*}}dvdxd\tau,
\end{array}
\ee
where $C>0$ is a generic constant independent of time $\mathcal{T}>0$.
\end{lemma}

\noindent {\it Proof.}  
We rewrite \eqref{1-nw-x}$_2$ with the non-fluid effect as
\be\label{nw-x-G}
\ba
&\quad\pa_i w_t+w\cdot\na\pa_i w+\frac{T}{n}\na\pa_i n+\na\pa_i T+\pa_i w\cdot\na w+\pa_i(\f{T}{n})\na n\\
&=-\r\pa_i (w-u)-\pa_i\r(w-u)
+\f{\pa_i n}{n^2}\int v\otimes v\cdot\na_xGdv-\f{1}{n}\int v\otimes v\cdot\na_x\pa_iGdv.
\ea
\ee
Similarly to \eqref{phi-1}, multiplying \eqref{1-ns-x}$_2$ by $\na\pa_i\r$ and \eqref{nw-x-G} by $\na \pa_i n$,
then adding the resulting equations and integrating it over $[0,t]\times\Omega$, we have
\begin{align}\label{phi-2}
\nm &\quad \sum_{i=1}^3\int(\pa_i u\cdot\na\pa_i\r+\pa_iw\cdot\na\pa_i n) dx\Big|_{\tau=0}^{\tau=t}
+\int_0^t\int\big(\f{\t}{\r}|\na^2\r|^2+\f{T}{n}|\na^2n|^2\big)dxd\tau\\
\nm &=\sum_{i=1}^3\int_0^t\int\Big[\pa_i\div u(u\cdot\na\pa_i\r+\r\pa_i\div u+\pa_i u\cdot\na\r+\pa_i\r\,\div u)
+\pa_i n(w-u)\cdot\na\pa_i\r\\
&\quad-(u\cdot\na\pa_i u+\na\pa_i\t+\pa_i u\cdot\na u+\pa_i(\f{\t}{\r})\na\r)\cdot\na\pa_i\r
+n\pa_i(w-u)\cdot\na\pa_i\r\Big] dxd\tau\\
\nm &
+\sum_{i=1}^3\int_0^t\int\Big[\pa_i\div w(w\cdot\na\pa_i n+n\pa_i\div w+\pa_i w\cdot\na n+\pa_i n\,\div w)
-\r\pa_i(w-u)\cdot\na\pa_i n\\
\nm &\quad-\pa_i\r(w-u)\cdot\na\pa_i n
-(w\cdot\na\pa_i w+\na\pa_i T+\pa_i w\cdot\na w+\pa_i(\f{T}{n})\na n)\cdot\na\pa_i n\\
\nm &\quad+\f{\pa_i n}{n^2}\na\pa_i n\cdot\int v\otimes v\cdot\na_xGdv
-\f{\na\pa_i n}{n}\cdot\int v\otimes v\cdot\na_x\pa_iGdv
\Big]dxd\tau=:I_7+I_8.
\end{align}
The terms $I_7$ and $I_8$ can be estimated similarly to Lemma 4.2 in \cite{LiWangWang2022}. By Sobolev inequality and the assumption \eqref{assumption}, we have
\begin{align}\label{I7}
I_7&\leq  \nm C\int_0^t\int\big[(|\na^2(u,\t)|+|\na(w-u)|)|\na^2\r|+(|\na^2(\r,u)|+|\na\r||\na u|)|\na^2 u|\big]dxd\tau\\ \nm
&\quad +C\int_0^t\int(|\na(\r,n)||w-u||\na^2\r|+|\na(\r,u,\t)|^2|\na^2\r|)dxd\tau\\
&\leq \f14\int_0^t\Big\|\sqrt{\f{\t}{\r}}\na^2\r\Big\|^2d\tau
+C\int_0^t(\|\na^2(u,\t)\|^2+\|\na(w-u)\|^2)d\tau\\ \nm
&\quad+C(\d+\v)\int_0^t(\|\na^2\r\|^2+\|w-u\|^2+\|\na(\r,u,\t)\|_{L^4}^4)d\tau\\ \nm
&\leq \big(\f14+C(\d+\v)\big)\int_0^t\Big\|\sqrt{\f{\t}{\r}}\na^2\r\Big\|^2d\tau
+C\int_0^t(\|\na^2(w,T)\|^2+\|\na(w-u)\|_{H^1}^2+\|\na^2(T-\t)\|^2)d\tau
\\ \nm
&\quad+C(\d+\v)\int_0^t(\|\na(\phi,\tilde\psi,\tilde\z)\|^2+\|w-u\|^2+\|\na(T-\t)\|^2)d\tau
+C\v^{\f12}.
\end{align}
Likewise,
\begin{align}\label{I8}\nm
\nm I_8&\leq C\int_0^t\int\big[(|\na^2(w,T)|+|\na(w-u)|+|\na(n,w,T)|^2)|\na^2n|
+(|\na^2(n,w)|+|\na n||\na w|)|\na^2 w|\big]dxd\tau\\
\nm &\quad +C\int_0^t\int |\na n||\na^2 n|\Big(\int\f{|\na_xG|^2}{\Mloc_*}dv\Big)^{\f12}dxd\tau
+C\int_0^t\int|\na^2 n|\Big(\int\f{|\na_x^2G|^2}{\Mloc_*}dv\Big)^{\f12}dxd\tau\\
&\leq
\big(\f14+C(\d+\v)\big)\int_0^t\Big\|\sqrt{\f{T}{n}}\na^2 n\Big\|^2d\tau+C\int_0^t(\|\na^2(w,T)\|^2+\|\na(w-u)\|^2)d\tau+C\v^{\f12}\\
\nm &+C(\d+\v)\int_0^t\|\na(\tilde\phi,\tilde\psi,\tilde\z)\|^2d\tau+C(\d+\v)\int_0^t\int\int\f{|\na_xG|^2}{\Mloc_*}dvdxd\tau
+C\int_0^t\int\int\f{|\na_x^2G|^2}{\Mloc_*}dvdxd\tau.
\end{align}
Substituting \eqref{I7} and \eqref{I8} into \eqref{phi-2}, and summing over $i=1,2,3$ yields the desired inequality \eqref{rho-2x}. The proof of Lemma \ref{le3.2} is complete.

\hfill $\Box$

\begin{lemma}\label{le3.3}
Under the same assumptions as in Proposition \ref{Prop3.2}, the following estimate holds for $t\in[0,\mathcal{T}]:$
\be\label{eG-1x}
\ba
&\quad \int\int\f{|\pa^{\alpha}G(t)|^2}{\Mloc_*} dv dx
+\int_0^t\int\int\f{|\na_v\pa^{\alpha}G|^2+(1+|v|^2)|\pa^{\alpha}G|^2}{\Mloc_*} dv dxd\tau\\
&\leq C\int\int\f{|\pa^{\alpha}G(0)|^2}{\Mloc_*} dv dx+C\v^{\f12}
+C\int_0^t\|\na\pa^{\alpha} (w,T)\|^2d\tau\\
&\quad +C(\d+\v^{\f12})\int_0^t(\|\pa^{\alpha}(\phi,\tilde\phi,\tilde\psi,\tilde\z)\|^2+\|\na(w-u,T-\t)\|^2+\|(\psi_t,\z_t)\|^2)d\tau\\
&\quad +C\int_0^t\int\int\f{|\na_x\pa^{\alpha}G|^2}{\Mloc_*} dvdxd\tau
+C(\d+\v^{\f12})\int_0^t\int\int\f{|(\na_v\na^2_xG,\na_v\wt{G},\wt{G})|^2}{\Mloc_*} dvdxd\tau,
\ea
\ee
where $C>0$ is a generic constant independent of time $\mathcal{T}>0$ and $\pa^{\alpha}$ denotes $\pa_t$ or $\na_x$.
\end{lemma}

\noindent {\it Proof.}  Applying $\pa^{\alpha}$ to  \eqref{non-fluid} yields
\begin{equation}\label{G-1x}
\ba
\pa^{\alpha}G_t -\r\Lop_{[w,T]}(\pa^{\alpha}G)&=
-\pa^{\alpha}\mb{P}_1(v\cdot\na_x\Mloc)
-\pa^{\alpha}\mb{P}_1(v\cdot\na_xG)-\r\pa^{\alpha}u\cdot\na_vG\\
&\quad +\r(w-u)\cdot\na_v\pa^{\alpha}G
+\r\pa^{\alpha}\t\Delta_vG-\r(T-\t)\Delta_v\pa^{\alpha}G\\
&\quad +\pa^{\alpha}\r\Lop_{[w,T]}G+\pa^{\alpha}\r(w-u)\cdot\na_vG
-\pa^{\alpha}\r(T-\t)\Delta_vG.
\ea
\end{equation}
Multiplying \eqref{G-1x} by $\pa^{\alpha}{G}/\Mloc_*$,
and integrating the resulting equation by parts over $[0,t]\times\Omega\times\mathbb{R}^3$, we have
\begin{align}\label{es-G-1x}
\nm & \int\int\f{|\pa^{\alpha}G|^2}{2\Mloc_*}dvdx\Big|_{\tau=0}^{\tau=t}
-\int_0^t\int\int\r\f{\pa^{\alpha}G}{\Mloc_*}\Lop_{[w,T]}(\pa^{\alpha}G)dvdxd\tau
\\
&=\int_0^t\int\int\Big\{ -\pa^{\alpha}\mb{P}_1(v\cdot\na_x\Mloc)
-\pa^{\alpha}\mb{P}_1(v\cdot\na_xG)-\r\pa^{\alpha}u\cdot\na_vG\\
\nm &\quad +\r(w-u)\cdot\na_v\pa^{\alpha}G
+\r\pa^{\alpha}\t\Delta_vG-\r(T-\t)\Delta_v\pa^{\alpha}G+\pa^{\alpha}\r\Lop_{[w,T]}G\\
\nm &\quad +\pa^{\alpha}\r(w-u)\cdot\na_vG-\pa^{\alpha}\r(T-\t)\Delta_vG
\Big\}\f{\pa^{\alpha}G}{\Mloc_*}dvdxd\tau:=\sum_{i=9}^{17}I_i.
\end{align}
The terms $I_i~(i=9,\cdots,17)$ on the right-hand side of \eqref{es-G-1x} are estimated as follows.
From \eqref{P1-vm}, $\mb{P}_1(v\cdot\na_x\Mloc)$ does not contain $\na n$. Thus, by Cauchy inequality and the assumption \eqref{assum-1}, we have
$$
\ba
|I_{9}|& \leq C\int_0^t\int \big(|\na\pa^{\alpha} (w,T)|+|\na (w,T)||\pa^{\alpha}(n,w,T)|\big)
\Big(\int\f{|\pa^{\alpha}G|^2}{\Mloc_*}dv\Big)^{\f12}dxd\tau\\
& \leq \big(\f{\bar\s}{8}+C(\d+\v)\big)\int_0^t\int\int\f{|\pa^{\alpha}G|^2}{\Mloc_*} dvdxd\tau
+C\int_0^t\|\na\pa^{\alpha}(w,T)\|^2d\tau\\
&\quad +C(\d+\v)\int_0^t\|\pa^{\alpha}(\tilde\phi,\tilde\psi,\tilde\z)\|^2d\tau+C\v.
\ea
$$
Similarly to \eqref{I4I5}, we have
$$
\ba
|I_{10}|& \leq \f{\bar\s}{8}\int_0^t\int\int\f{1+|v|^2}{\Mloc_*}|\pa^{\alpha}G|^2 dvdxd\tau
+C\int_0^t\int\int\f{|\na_x\pa^{\alpha}G|^2}{\Mloc_*} dvdxd\tau\\
&\quad +C(\d+\v)\int_0^t\int\int\f{|(\pa^{\alpha}G,\na_xG)|^2}{\Mloc_*} dvdxd\tau.
\ea
$$
By Cauchy inequality, the assumption \eqref{assum-1}, and the fact that $\|\pa^{\alpha}u\|\leq C(\d+\v^{\f12})$, we deduce that
\begin{align} \nm
|I_{11}|&\leq  C\int_0^t\int|\pa^{\alpha}u|\Big(\int\f{|\na_vG|^2}{\Mloc_*}dv\Big)^{\f12}
\Big(\int\f{|\pa^{\alpha}G|^2}{\Mloc_*}dv\Big)^{\f12}dxd\tau\\ \nm
&\leq C\int_0^t\Big\|\int\f{|\na_vG|^2}{\Mloc_*}dv\Big\|^{\f12}_{L^{\infty}_x}\|\pa^{\alpha}u\|
\Big(\int\int\f{|\pa^{\alpha}G|^2}{\Mloc_*}dvdx\Big)^{\f12}d\tau\\ \nm
&\leq \big(\f{\bar\sigma}{8}+C(\d+\v^{\f12})\big)\int_0^t\int\int\f{|\pa^{\alpha}G|^2}{\Mloc_*}dvdxd\tau\\ \nm
&\quad+C(\d+\v^{\f12})\int_0^t\int\int\f{|(\na_v\wt{G},\na_v\na_xG,\na_v\na^2_xG)|^2}{\Mloc_*}dvdxd\tau\\ \nm
&\quad+C\v^{\f12}\int_0^t(\|\na\tilde\psi\|^2+\|\na(w-u)\|^2+\|\psi_t\|^2)d\tau+C\v^{\f12},
\end{align}
where in the last inequality we have used three-dimensional Sobolev inequality
\be\label{im}
\begin{array}{l}
\di \Big\|\int \f{|\na_vG|^2}{\Mloc_*} dv\Big\|_{L^{\infty}_x}=\int\big\|\f{\na_vG}{\sqrt{\Mloc_*}}\Big\|^2_{L_x^{\infty}}dv
\leq C\int\big\|\f{\na_vG}{\sqrt{\Mloc_*}}\Big\|^2_{H_x^2}dv\\
\di \leq C\int\int \f{|\na_vG|^2+|\na_v\na_xG|^2+|\na_v\na_x^2G|^2}{\Mloc_*} dv dx.
\end{array}
\ee
By Cauchy inequality and the assumption \eqref{assum-1}, we have
$$
\ba
|I_{12}|&\leq C\int_0^t\int\|\tilde\psi-\psi\|_{L^{\infty}}\Big(\int\f{|\pa^{\alpha}G|^2}{\mathbf{M}_*}dv\Big)^{\f12}
\Big(\int\f{|\na_v\pa^{\alpha}G|^2}{\mathbf{M}_*}dv\Big)^{\f12}dxd\tau\\
&\leq C(\d+\v)\int_0^t\int\int\f{|(\pa^{\alpha}G,\na_v\pa^{\alpha}G)|^2}{\Mloc_*}dvdxd\tau.
\ea
$$
Integrating by parts in $v$ over $\mathbb{R}^3$,  using  $\|\pa^{\alpha}\t\|\leq\|\pa^{\alpha}\z\|+\|\pa^{\alpha}\bar\t\|\leq C(\d+\v^{\f12})$,  and applying \eqref{im}, we obtain
\begin{align}\nm
I_{13}&=-\int_0^t\int  \r \pa^{\alpha}\t\int\na_vG\cdot\na_v\Big(\f{\pa^{\alpha}G}{\Mloc_*}\Big)dv\,dxd\tau\\ \nm
&=-\int_0^t\int \r \pa^{\alpha}\t\int\na_vG\cdot\Big(\f{\na_v\pa^{\alpha}G}{\Mloc_*}
+\f{v-w_*}{T_*}\f{\pa^{\alpha}G}{\Mloc_*}\Big)dv\,dxd\tau\\ \nm
&\leq C\int_0^t\|\pa^{\alpha}\t\|\Big\|\int\f{|\na_vG|^2}{\Mloc_*}dv\Big\|^{\f12}_{L^{\infty}_x}
\Big(\int\int\f{|\na_v\pa^{\alpha}G|^2+(1+|v|^2)|\pa^{\alpha}G|^2}{\Mloc_*}dvdx\Big)^{\f12}d\tau\\ \nm
&\leq \big(\f{\bar\sigma}{8}+C(\d+\v^{\f12})\big)\int_0^t\int\int\f{|\na_v\pa^{\alpha}G|^2+(1+|v|^2)|\pa^{\alpha}G|^2}{\Mloc_*}dvdxd\tau\\ \nm
&\quad+C(\d+\v^{\f12})\int_0^t\int\int
\f{|(\na_v\wt{G},\na_v\na_xG,\na_v\na^2_xG)|^2}{\Mloc_*}dvdxd\tau\\ \nm
&\quad+C\v^{\f12}\int_0^t(\|\na\tilde\z\|^2+\|\na(T-\t)\|^2+\|\z_t\|^2)d\tau+C\v^{\f12}.
\end{align}
Integration by parts in $v$ over $\mathbb{R}^3$ and Cauchy inequality lead to
$$
\ba
|I_{14}|&\leq  C\d\int_0^t\int\int\f{|\na_v\pa^{\alpha}G|^2+(1+|v|^2)|\pa^{\alpha}G|^2}{\Mloc_*}dvdxd\tau.
\ea
$$
Similarly to $I_{13}$, we have
\begin{align}\nm
I_{15}&=-\int_0^t\int  \nm \pa^{\alpha}\r\int\big((v-w)G+T\na_vG\big)\cdot\na_v\Big(\f{\pa^{\alpha}G}{\Mloc_*}\Big)dv\,dxd\tau\\ \nm
&\leq \big(\f{\bar\sigma}{8}+C(\d+\v^{\f12})\big)\int_0^t\int\int\f{|\na_v\pa^{\alpha}G|^2+(1+|v|^2)|\pa^{\alpha}G|^2}{\Mloc_*}dvdxd\tau\\ \nm
&\quad+C(\d+\v^{\f12})\int_0^t\int\int
\f{|(\na_v\wt{G},\na_v\na_xG,\na_v\na^2_xG)|^2}{\Mloc_*}dvdxd\tau\\ \nm
&\quad+C(\d+\v^{\f12})\int_0^t\int\int
\f{|(\wt{G},\na_xG,\na^2_xG)|^2}{\Mloc_*}dvdxd\tau+C\v^{\f12}\int_0^t\|\pa^{\alpha}\phi\|^2d\tau+C\v^{\f12},
\end{align}
and
\begin{align}\nm
I_{17}&=\int_0^t\int  \nm \pa^{\alpha}\r(T-\t)\int\na_vG\cdot\na_v\Big(\f{\pa^{\alpha}G}{\Mloc_*}\Big)dv\,dxd\tau\\ \nm
&\leq \big(\f{\bar\sigma}{8}+C(\d+\v^{\f12})\big)\int_0^t\int\int\f{|\na_v\pa^{\alpha}G|^2+(1+|v|^2)|\pa^{\alpha}G|^2}{\Mloc_*}dvdxd\tau\\ \nm
&\quad+C(\d+\v^{\f12})\int_0^t\int\int
\f{|(\na_v\wt{G},\na_v\na_xG,\na_v\na^2_xG)|^2}{\Mloc_*}dvdxd\tau
+C\v^{\f12}\int_0^t\|\pa^{\alpha}\phi\|^2d\tau+C\v^{\f12},
\end{align}
Similarly to $I_{11}$, we have
\begin{align} \nm
|I_{16}|&\leq  \|w-u\|_{L^{\infty}}\int_0^t\int|\pa^{\alpha}\r|\Big(\int\f{|\na_vG|^2}{\Mloc_*}dv\Big)^{\f12}
\Big(\int\f{|\pa^{\alpha}G|^2}{\Mloc_*}dv\Big)^{\f12}dxd\tau\\ \nm
&\leq \big(\f{\bar\sigma}{8}+C(\d+\v^{\f12})\big)\int_0^t\int\int\f{|\pa^{\alpha}G|^2}{\Mloc_*}dvdxd\tau
+C\v^{\f12}\int_0^t\|\pa^{\alpha}\phi\|^2d\tau+C\v^{\f12}\\ \nm
&\quad+C(\d+\v^{\f12})\int_0^t\int\int\f{|(\na_v\wt{G},\na_v\na_xG,\na_v\na^2_xG)|^2}{\Mloc_*}dvdxd\tau.
\end{align}
Substituting the above estimates into \eqref{es-G-1x}, choosing $\d$ and $\v$ suitably small yields the desired inequality \eqref{eG-1x}. This  completes the proof of Lemma \ref{le3.3}.

\hfill $\Box$

\noindent{{\emph{Proof of Proposition \ref{Prop3.2}}}}:
We first multiply \eqref{1-es} by a sufficiently large constant $C$, add the resulting inequality to \eqref{rho-2x}, and obtain
\begin{align}\label{1x}\nm
&\|\na(\r,u,\t,n,w,T)(t)\|^2+\int_0^t(\|\na^2(\r,n,w,T)\|^2+\|\na(w-u,T-\t)\|^2)d\tau
\\
& \leq C(\mathcal{N}(0)+\v^{\f12})+C\|\na^2(\r,n)(t)\|^2+C\int_0^t\|\na^2(w-u,T-\t)\|^2d\tau\\ \nm
&\quad +C(\d+\v)\int_0^t(\|\na(\phi,\tilde\phi,\tilde\psi,\tilde\z)\|^2+\|(\z_t,\tilde\phi_t,\tilde\z_t)\|^2
+\|(w-u,T-\t)\|^2)d\tau\\ \nm
&\quad +C(\d+\v)\int_0^t\int\int\f{|(\na_v\na_xG,\pa^{\alpha}G,\na_v\wt{G},\wt{G})|^2}{\Mloc_*}dvdxd\tau
+C\int_0^t\int\int\f{|\na_x\pa^{\alpha}G|^2}{\Mloc_*}dvdxd\tau.
\end{align}

Next, we rewrite \eqref{1-nw-x}$_3$ with non-fluid diffusive effect as
\begin{align}\label{3.79}
&\quad\f32\pa_i T_t+\f32w\cdot\na\pa_i T+T\pa_i\dv w+\f32\pa_i w\cdot\na T+\pa_i T\dv w\\ \nm
&=-3\r\pa_i(T-\t)-3\pa_i\r(T-\t)
+\pa_i\Big(\f{1}{n}\big(-\int\f12|v|^2v\cdot\na_xG dv+w\cdot\int v\otimes v\cdot\na_x G dv\big)\Big).
\end{align}
Multiplying \eqref{1-ns-x}$_1$, \eqref{1-ns-x}$_2$, \eqref{1-ns-x}$_3$, \eqref{1-nw-x}$_1$, \eqref{nw-x-G}, and \eqref{3.79} by $\pa_i\r_t$, $\pa_i u_t$, $\pa_i \t_t$, $\pa_i n_t$,  $\pa_i w_t$
and $\pa_i T_t~(i=1,2,3)$, respectively, summing the resulting equations over $i=1,2,3$, and integrating by parts over $[0,t]\times\Omega$, we obtain
\begin{align}\label{phi-xt}
\int_0^t\|\na(\r_t, u_t,\t_t, n_t, w_t, T_t)\|^2 d\tau& \leq C\int_0^t(\|\na^2(\r,n,w,T)\|^2+\|\na(w-u,T-\t)\|_{H^1}^2)d\tau
\\ \nm
&+C(\d+\v)\int_0^t(\|\na(\phi,\tilde\phi,\tilde\psi,\tilde\z)\|^2+\|(w-u,T-\t)\|^2)d\tau+C\v^{\f12} \\ \nm
&+C\int_0^t\int\int\frac{|\na^2_xG|^2}{\mathbf{M}_*}d v dxd\tau
+C(\d+\v)\int_0^t\int\int\f{|\na_xG|^2}{\Mloc_*}dvdxd\tau.
\end{align}
Multiplying \eqref{1x} by a sufficiently large constant $C$, and adding the resulting inequality to \eqref{phi-xt}, we obtain
\begin{align}\label{fluid-1x}
& \nm \|\na(\r,u,\t,n,w,T)(t)\|^2+\int_0^t(\|\na\pa^{\alpha}(\r,n,w,T)\|^2+\|\na( u_t,\t_t)\|^2
+\|\na(w-u,T-\t)\|^2)d\tau\\
&\leq  C(\mathcal{N}(0)+\v^{\f12})+C\|\na^2(\r,n)(t)\|^2+C\int_0^t\|\na^2(w-u,T-\t)\|^2d\tau
\\ \nm
&\quad +C(\d+\v)\int_0^t(\|\pa^{\alpha}(\phi,\tilde\phi,\tilde\psi,\tilde\z)\|^2+\|\z_t\|^2+\|(w-u,T-\t)\|^2)d\tau\\ \nm
&\quad +C(\d+\v)\int_0^t\int\int\f{|(\na_v\na_xG,\pa^{\alpha}G,\na_v\wt{G},\wt{G})|^2}{\mb{M_*}}dvdxd\tau
+C\int_0^t\int\int\f{|\na_x\pa^{\alpha}G|^2}{\mb{M_*}}dvdxd\tau.
\end{align}
Multiplying \eqref{fluid-1x} by a sufficiently large constant $C$,
adding the resulting inequality to \eqref{eG-1x}, we obtain the desired inequality \eqref{high-1}. This completes the proof of Proposition \ref{Prop3.2}.

\hfill $\Box$

%%%%%%%%%%%%%%%%%%%%%%%%%%%%%%%%%%%%%%%%%%%%%%%%%%%%%%%%%%%%%%%%%%%%%%%%%%%%%%%

\subsection{Higher-order estimate (II)}
\label{SSpf-prop-3}

The proof of Proposition \ref{Prop3.3} is divided into the following several lemmas.

\begin{lemma}\label{le3.4}
Under the same assumptions as in Proposition \ref{Prop3.3}, the following estimate holds for $t\in[0,\mathcal{T}]:$
\be\label{e-2x}
\begin{array}{l}
\di \|\na^2(\r,u,\t,n,w,T)(t)\|^2+\int_0^t(\|\na^3(w,T)\|^2+\|\na^2(w-u,T-\t)\|^2)d\tau
\\
\di\leq C(\mathcal{N}(0)+\v)
+C\int_0^t\int\int\f{|\na_x^2\pa^{\alpha}G|^2}{\Mloc_*}dvdxd\tau+C(\d+\v)\int_0^t\|\na^3(\r,n)\|^2d\tau\\
\di+C(\d+\v)\int_0^t(\|\na^2(\r,n,w,T)\|^2+\|\na(\phi,\tilde\phi,\tilde\psi,\tilde\z)\|^2+\|\na(w-u,T-\t)\|^2)d\tau\\
\di+C(\d+\v)\int_0^t\int\int\f{|(\na_v\na_x^2G,\na_x\pa^{\alpha}G,\na_v\na_xG,
\pa^{\alpha}G,\na_v\wt{G},\wt{G})|^2}{\mb{M_*}}
dvdxd\tau,
\end{array}
\ee
where $C>0$ is a generic constant independent of time $\mathcal{T}>0$ and $\pa^{\alpha}$ denotes $\pa_t$ or $\na_x$.
\end{lemma}

\noindent {\it Proof.}  Applying $\pa_j~(j=1,2,3)$ to \eqref{1-ns-x} and \eqref{1-nw-x}, we have
\be\label{2-e-xx}
\begin{cases}
\di \pa_j\pa_i\r_t+u\cdot\na\pa_j\pa_i\r+\r\pa_j\pa_i\dv u
+(\pa_j u\cdot\na\pa_i\r+\pa_i u\cdot\na\pa_j\r+\dv u\pa_j\pa_i\r)\\
\quad+(\pa_j\r\pa_i\dv u+\pa_j\pa_i u\cdot\na\r+\pa_i\r\pa_j\dv u)=0,\\
\di \pa_j\pa_i u_t+u\cdot\na\pa_j\pa_i u+\frac{\t}{\r}\na\pa_j\pa_i\r+\na\pa_j\pa_i\t
+\pa_j u\cdot\na\pa_i u+\pa_j\pa_i u\cdot\na u\\
\di\quad+\pa_i u\cdot\na\pa_j u+\pa_j\big(\f{\t}{\r}\big)\na\pa_i\r+\pa_j\pa_i\big(\f{\t}{\r}\big)\na\r
+\pa_i\big(\f{\t}{\r}\big)\na\pa_j\r\\
\di\quad=n\pa_j\pa_i(w-u)+\pa_j\big(\pa_i n (w-u)\big)+\pa_j n\pa_i(w-u),\\
\di \f32\pa_j\pa_i\t_t+\f32u\cdot\na\pa_j\pa_i\t+\t\pa_j\pa_i\dv u+\f32\pa_j u\cdot\na\pa_i\t+\f32\pa_j\pa_i u\cdot\na\t
+\f32\pa_i u\cdot\na\pa_j\t\\
\di\quad+\pa_j\t\pa_i\dv u+\pa_j\pa_i\t\dv u+\pa_i\t\pa_j\dv u
=3n\pa_j\pa_i(T-\t)\\
\di\quad +\pa_j\big(3\pa_i n (T-\t)\big)+3\pa_j n\pa_i(T-\t)+\pa_j\pa_i\big(n |w-u|^2\big),
\end{cases}
\ee
and
\begin{align}\label{2-nw-xx}
\begin{cases}
& \pa_j\pa_i n_t+w\cdot\na\pa_j\pa_i n+ n\pa_j\pa_i\dv w
+(\pa_j w\cdot\na\pa_i n+\pa_i w\cdot\na\pa_j n+\dv w\pa_j\pa_i n)\\
&\quad+(\pa_j n\pa_i\dv w+\pa_j\pa_i w\cdot\na n+\pa_i n\pa_j\dv w)=0,\\
& \pa_j\pa_i w_t+w\cdot\na\pa_j\pa_i w+\frac{T}{n}\na\pa_j\pa_i n+\na\pa_j\pa_i T
+\pa_j w\cdot\na\pa_i w+\pa_j\pa_i w\cdot\na w\\
&\quad+\pa_i w\cdot\na\pa_j w+\pa_j\big(\f{T}{n}\big)\na\pa_i n+\pa_j\pa_i\big(\f{T}{n}\big)\na n
+\pa_i\big(\f{T}{n}\big)\na\pa_j n =-\r\pa_j\pa_i(w-u)\\
&\quad-\pa_j\big(\pa_i\r(w-u)\big)-\pa_j\r\pa_i(w-u)
+\f{T}{2\r}(\Delta\pa_j\pa_i w+\f13\na\pa_j\pa_i\dv w)+\pa_j\big(\f{T}{2\r}\big)(\Delta\pa_i w+\f13\na\pa_i\dv w)\\
&\quad+\pa_j\Big(\pa_i (\f{T}{2\r})(\Delta w+\f13\na\dv w)+\pa_i\big(\na(\f{nT}{2\r})\cdot\f{\mathbb{S}(w)}{n}\big)\Big)
-\pa_j\pa_i\Big(\f{1}{n}\int v\otimes v\cdot\na_x\Pi dv\Big),\\
&\f32\pa_j\pa_i T_t+\f32w\cdot\na\pa_j\pa_i T+T\pa_j\pa_i\dv w+\f32\pa_jw\cdot\na\pa_i T+\f32\pa_j\pa_i w\cdot\na T+\f32\pa_iw\cdot\na\pa_j T
\\
&\quad+\pa_j T\pa_i\dv w+\pa_j\pa_i T\dv w+\pa_i T\pa_j\dv w=-3\r\pa_j\pa_i(T-\t)
-3\pa_j\r\pa_i(T-\t)-3\pa_j\big(\pa_i\r(T-\t)\big)\\
&\quad
+\pa_j\Big(\pa_i(\f{5T}{6\r})\Delta T+\pa_i\big(\na(\f{5nT}{6\r})\cdot\f{\na T}{n}\big)
+\pa_i\big(\f{T}{\r}\big(\f{(\na w+(\na w)^t)^2}{4}-\f13(\dv w)^2\big)\big)\Big)\\
&\quad+\f{5T}{6\r}\Delta\pa_j\pa_i T+\pa_j (\f{5T}{6\r})\Delta\pa_i T
+\pa_j\pa_i\Big(\f{1}{n}\big(-\int\f12|v|^2v\cdot\na_x\Pi dv+w\cdot\int v\otimes v\cdot\na_x \Pi dv\big)\Big).
\end{cases}
\end{align}
Multiplying $\eqref{2-e-xx}_1$, $\eqref{2-e-xx}_2$, $\eqref{2-e-xx}_3$, $\eqref{2-nw-xx}_1$, $\eqref{2-nw-xx}_2$, and $\eqref{2-nw-xx}_3$ by $\f{\t}{\r n}\pa_j\pa_i\r$, $\f{\r}{n}\pa_j\pa_i u$, $\f{\r}{\t n}\pa_j\pa_i\t$,  $\f{T}{n^2}\pa_j\pa_i n$,
$\pa_j\pa_i w$, and $\f{\pa_j\pa_i T}{T}$, respectively, then summing the resulting equations over $i,j=1,2,,3$, and integrating by parts over $[0,t]\times\Omega$, we obtain
$$
\ba
& \int\Big(\f{\t}{\r n}\f{|\na^2\r|^2}{2}+\f{\r}{n}\f{|\na^2 u|^2}{2}+\f{3\r}{2\t n}\f{|\na^2\t|^2}{2}
+\f{T}{n^2}\f{|\na^2 n|^2}{2}+\f{|\na^2 w|^2}{2}+\f{3}{2T}\f{|\na^2 T|^2}{2}\Big)dx\Big|_{\tau=0}^{\tau=t}\\
&+\int_0^t\int\big(\f{T}{2\r}\big(|\na^3 w|^2+\f{|\na^2\dv w|^2}{3}\big)+\f{5}{6\r}|\na^3 T|^2
+\r|\na^2(w-u)|^2+\f{3\r}{T}|\na^2(T-\t)|^2\big)dxd\tau\\
&=:\sum_{i,j=1}^3\int_0^t\int\f{3\r}{\t T}(T-\t)\pa_j\pa_i\t\pa_j\pa_i(T-\t)dxd\tau+\sum_{j=1}^3J_i,
\ea
$$
where
\begin{align}\nm
J_1&\leq C\int_0^t\int(|\na(\r,u,\t)||\na^2(\r,u,\t)|^2+|\na\r|^2|\na\t||\na^2u|)dxd\tau\\ \nm
&\quad +C\int_0^t\int(|\na n||\na(w-u,T-\t)|+|(w-u,T-\t)||\na^2 n|)|\na^2(u,\t)|dxd\tau\\ \nm
&\quad +C \int_0^t\int(|w-u||\na^2(w-u)|+|\na(w-u)|^2)|\na^2\t|dxd\tau,\\ \nm
J_2&\leq C\int_0^t\int(|\na(\r,n,w,T)||\na^2(n,w,T)|^2+|\na n|^2|\na T||\na^2 w|)dxd\tau\\ \nm
&\quad +C \int_0^t\int|\na(\r,n,w,T)||\na^3(\r,n,w,T)||\na^2(w,T)|dxd\tau\\ \nm
&\quad +C \int_0^t\int(|\na^2(\r,n,w,T)|^2+|\na(\r,n,w,T)|^4)|\na^2(w,T)|dxd\tau\\ \nm
&\quad +C \int_0^t\int(|\na^2\r||(w-u,T-\t)|+|\na\r||\na(w-u,T-\t)|)|\na^2(w,T)|dxd\tau,\\ \nm
J_3&=\sum_{i,j=1}^3\int_0^t\int\pa_j^2\pa_iw\pa_i\Big(\f{1}{n}\int v\otimes v\cdot\na_x\Pi dv\Big)dxd\tau\\ \nm
&\quad -\sum_{i,j=1}^3\int_0^t\int\pa_j\big(\f{\pa_j\pa_i T}{T}\big)
\pa_i\Big(\f{1}{n}\big(-\int\f12|v|^2v\cdot\na_x\Pi dv+w\cdot\int v\otimes v\cdot\na_x \Pi dv\big)\Big)dxd\tau.
\end{align}
By the assumption \eqref{assumption} and Cauchy inequality, we have
$$
\ba
J_1&\leq C\int_0^t\|\na(\r,u,\t)\|_{L^{\infty}}\|\na^2(\r,u,\t)\|^2d\tau
+C\int_0^t\|\na\t\|_{L^{\infty}}\|\na\r\|_{L^4}^2\|\na^2 u\|d\tau\\
&\quad+C\int_0^t(\|\na n\|_{L^{\infty}}\|\na(w-u,T-\t)\|
+\|(w-u,T-\t)\|_{L^{\infty}}\|\na^2 n\|)\|\na^2(u,\t)\|d\tau\\
&\quad+C\int_0^t(\|w-u\|_{L^{\infty}}\|\na^2(w-u)\|+\|\na(w-u)\|_{L^{\infty}}\|\na(w-u)\|)\|\na^2 \t\|d\tau\\
&\leq C(\d+\v)\int_0^t(\|\na^2(\r,u,\t,n)\|^2+\|\na(w-u,T-\t)\|_{H^1}^2)\\
&\leq C(\d+\v)\int_0^t(\|\na^2(\r,n,w,T)\|^2+\|\na(w-u,T-\t)\|_{H^1}^2).
\ea
$$
By Sobolev inequality, the assumption \eqref{assumption}, and \eqref{assum-3}, we have
\begin{align}\label{L3-assump}
\nm \|\na^2(\r,n,w,T)\|_{L^3}&\leq C\|\na^2(\r,n,w,T)\|^{\f12}\|\na^2(\r,n,w,T)\|_{L^6}^{\f12}\\
&\leq C\|\na^2(\r,n,w,T)\|_{H^1}\leq C(\d+\v).
\end{align}
Using Cauchy inequality, the assumption \eqref{assumption}, and \eqref{L3-assump}, we have
\begin{align} \nm
J_2&\leq C\int_0^t\|\na(\r,n,w,T)\|_{L^{\infty}}\|\na^2(n,w,T)\|^2d\tau
+C\int_0^t\|\na T\|_{L^{\infty}}\|\na n\|_{L^4}^2\|\na^2 w\|d\tau\\ \nm
&\quad+C\int_0^t\|\na(\r,n,w,T)\|_{L^{\infty}}\|\na^3(\r,n,w,T)\|\|\na^2(w,T)\|d\tau\\ \nm
&\quad+C\int_0^t(\|\na^2(\r,n,w,T)\|^2_{L^3}\|\na^2(w,T)\|_{L^3}
+\|(\bar\r_{x_1},\bar n_{x_1},\bar u_{1x_1},\bar\t_{x_1})\|_{L^8}^4\|\na^2(w,T)\|)d\tau\\ \nm
&\quad+C\int_0^t
\|\na(\phi,\tilde\phi,\tilde\psi,\tilde\z)\|_{L^{\infty}}^3\|\na(\phi,\tilde\phi,\tilde\psi,\tilde\z)\|\|\na^2(w,T)\|d\tau\\ \nm
&\quad+C\int_0^t (\|(w-u,T-\t)\|_{L^{\infty}}\|\na^2\r\|+\|\na\r\|_{L^{\infty}}\|\na(w-u,T-\t)\|)\|\na^2(w,T)\|d\tau\\ \nm
&\leq C(\d+\v)\int_0^t(\|\na^2(\r,n,w,T)\|_{H^1}^2+\|\na(w-u,T-\t)\|^2+\|\na(\phi,\tilde\phi,\tilde\psi,\tilde\z)\|^2)+C\v.
\end{align}
Cauchy inequality and \eqref{assum-1} yield
$$
\ba
J_3&\leq \f{1}{8}\int_0^t\|\na^3(w,T)\|^2d\tau+C(\d+\v)\int_0^t\|\na^2T\|^2d\tau\\
&\quad+C\int_0^t\int\Big(\Big|\int\f12|v|^2v\cdot\na_x^2\Pi dv\Big|^2
+\Big| \int v\otimes v\cdot\na_x^2 \Pi dv\Big|^2\Big)dxd\tau\\
&\quad+C(\d+\v)\int_0^t\int\Big(\Big|\int\f12|v|^2v\cdot\na_x\Pi dv\Big|^2
+\Big| \int v\otimes v\cdot\na_x \Pi dv\Big|^2\Big)dxd\tau.
\ea
$$
By a similar but more laborious computation as in \eqref{Pi-x}, we have
\begin{align}\label{complex} \nm
&\quad \Big|\int\f12|v|^2v\cdot\na_x^2\Pi dv\Big|
+\Big| \int v\otimes v\cdot\na_x^2 \Pi dv\Big|\\ \nm
&\leq C\Big(\int\f{|\na^2_x\pa^{\alpha}G|^2}{\Mloc_*}dv\Big)^{\f12}
+C|w-u|\Big(\int\f{|\na_v\na_x^2G|^2}{\Mloc_*} dv\Big)^{\f12}
+C|\na(w,T)|\Big(\int\f{|\na_x\pa^{\alpha}G|^2}{\Mloc_*} dv\Big)^{\f12}\\ \nm
&\quad+C|\na^2(w,T)|\Big(\int\f{|\pa^{\alpha}G|^2}{\Mloc_*} dv\Big)^{\f12}
+C(|\na(w-u,T-\t)|+|\na(w,T)||(w-u,T-\t)|+|\na(w,T)|^2)\\
&\quad
\cdot\Big(\int\f{|(\na_v\na_xG,\pa^{\alpha}G)|^2}{\Mloc_*} dv\Big)^{\f12}+C\big(|\na^2(w-u,T-\t)|+|\na(w,T)||\na(w-u,T-\t)|\\ \nm
&\quad +(|\na^2(w,T)|+|\na(w,T)|^2)|(w-u,T-\t)|\big)
\cdot\Big(\int\f{|(\na_vG,G)|^2}{\Mloc_*} dv\Big)^{\f12}\\ \nm
&\leq C\Big(\int\f{|\na^2_x\pa^{\alpha}G|^2}{\Mloc_*}dv\Big)^{\f12}
+C(\d+\v)\Big(\int\f{|(\na_v\na^2_xG,\na_x\pa^{\alpha}G,\na_v\na_xG,\pa^{\alpha}G,\na_v\wt{G},\wt{G})|^2}{\Mloc_*}dv\Big)^{\f12}\\ \nm
&\quad+C(\d+\v)|\na^2(w,T)|+C(|\na(w,T)||\na(w-u,T-\t)|+|\na(w,T)|^2|(w-u,T-\t)|)\bar u_{1x_1}\\ \nm
&\quad+C(|\na^2(w-u,T-\t)|+|\na^2(w,T)||(w-u,T-\t)|)\Big(\int\f{|(\na_vG,G)|^2}{\Mloc_*} dv\Big)^{\f12},
\end{align}
where in the last inequality we have used \eqref{assum-1} and \eqref{assum-2}. Similarly to \eqref{assum-2}, by Sobolev inequality, we have
$$
\begin{aligned}
\Big(\int\f{|(\na_vG,G)|^2}{\Mloc_*} dv\Big)^{\f12}
&\leq C\Big(\int\int\f{|(\na_v\wt{G},\na_v\na_xG,\na_v\na_x^2G)|^2}{\Mloc_*} dvdx\Big)^{\f12}\\
&\quad+C\Big(\int\int\f{|(\wt{G},\na_xG,\na_x^2G)|^2}{\Mloc_*} dvdx\Big)^{\f12}+C\|\bar u_{1x_1}\|_{L^2_{x_1}}\\
&\leq C\Big(\int\int\f{|(\na_v\wt{G},\na_v\na_xG,\na_v\na_x^2G)|^2}{\Mloc_*} dvdx\Big)^{\f12}+C(\d+\v^{\f12}),
\end{aligned}
$$
which together with \eqref{complex} yields
$$
\ba
&\int_0^t\int\Big(\Big|\int\f12|v|^2v\cdot\na_x^2\Pi dv\Big|^2
+\Big| \int v\otimes v\cdot\na_x^2 \Pi dv\Big|^2\Big)dxd\tau
\leq C\int_0^t\int\int\f{|\na^2_x\pa^{\alpha}G|^2}{\Mloc_*}dvdxd\tau\\
&\qquad\quad
+C(\d+\v)\int_0^t\int\int\Big(\f{|(\na_v\na^2_xG,\na_v\na_xG,\na_v\wt{G})|^2}{\Mloc_*}
+\f{|(\na_x\pa^{\alpha}G,\pa^{\alpha}G,\wt{G})|^2}{\Mloc_*}\Big)dvdxd\tau\\
&\qquad\quad+C(\d+\v)\int_0^t(\|\na^2(w,T)\|^2+\|\na(w-u,T-\t)\|_{H^1}^2+\|\na(\tilde\psi,\tilde\z)\|^2)d\tau+C\v,
\ea
$$
where we have used the fact that $$\|\na^2(w-u,T-\t)\|=\|\na^2(\tilde\psi-\psi,\tilde\z-\z)\|\leq C\d,
\quad\|\na^2(w,T)\|\leq C(\d+\v).
$$
Combining the estimates for $J_1$, $J_2$ and $J_3$ yields \eqref{e-2x}. The proof of Lemma \ref{le3.4} is complete.

\hfill $\Box$

\begin{lemma}\label{le3.5}
Under the same assumptions as in Proposition \ref{Prop3.3}, the following estimate holds for $t\in[0,\mathcal{T}]:$
\begin{align}\label{3-rho} \nm
 \int_0^t\|\na^3(\r,n)\|^2d\tau&\leq C(\mathcal{N}(0)+\v)+C\|\na^2(u,w)(t)\|^2
+C\|\na^3(\r,n)(t)\|^2\\
&\quad +C\int_0^t(\|\na^3(w,T)\|^2+\|\na^2(w-u,T-\t)\|_{H^1}^2)d\tau\\ \nm
&\quad +C(\d+\v)\int_0^t(\|\na^2(\r,n,w,T)\|^2+\|\na(w-u)\|^2+\|\na(\phi,\tilde\phi)\|^2)d\tau\\ \nm
&\quad +C\int_0^t\int\int\f{|\na_x^3G|^2}{\Mloc_*}dvdxd\tau
+C(\d+\v)\int_0^t\int\int\f{|\na_x^2G|^2}{\Mloc_*}dvdxd\tau,
\end{align}
where $C>0$ is a generic constant independent of time $\mathcal{T}>0$.
\end{lemma}

\noindent {\it Proof.} Applying $\pa_j~(j=1,2,3)$ to \eqref{nw-x-G} yields
\be\label{nw-2x-G}
\ba
& \pa_j\pa_i w_t+w\cdot\na\pa_j\pa_i w+\frac{T}{n}\na\pa_j\pa_i n+\na\pa_j\pa_i T
+\pa_j w\cdot\na\pa_i w+\pa_j\pa_i w\cdot\na w+\pa_i w\cdot\na\pa_j w\\
&+\pa_j(\f{T}{n})\na\pa_i n+\pa_j\pa_i(\f{T}{n})\na n+\pa_i(\f{T}{n})\na\pa_j n\\
&=-\r\pa_j\pa_i(w-u)-\pa_j\big(\pa_i\r(w-u)\big)-\pa_j\r\pa_i(w-u)
-\pa_j\pa_i\Big(\f{1}{n}\int v\otimes v\cdot\na_xGdv\Big).
\ea
\ee
Multiplying \eqref{2-e-xx}$_2$ by $\na\pa_j\pa_i\r$ and \eqref{nw-2x-G} by $\na\pa_j\pa_i n$, integrating the resulting equations by parts over $[0,t]\times\Omega$, and adding them together, we obtain the desired inequality \eqref{3-rho} similarly to Lemma \ref{le3.2}; see also Lemma 4.5 in \cite{LiWangWang2022}.  The proof of Lemma \ref{le3.5} is complete.

\hfill $\Box$

\begin{lemma}\label{le3.6}
Under the same assumptions as in Proposition \ref{Prop3.3}, the following estimate holds for $t\in[0,\mathcal{T}]:$
\begin{align}\label{eG-2x}
\nm
&\quad \int\int\f{|\na_x\pa^{\alpha}G(t)|^2}{\Mloc_*} dv dx
+\int_0^t\int\int\f{|\na_v\na_x\pa^{\alpha}G|^2+(1+|v|^2)|\na_x\pa^{\alpha}G|^2}{\Mloc_*} dv dxd\tau\\
\nm &\leq C\int\int\f{|\na_x\pa^{\alpha}G(0)|^2}{\Mloc_*} dv dx+C\v+C\int_0^t\|\na^2\pa^{\alpha} (w,T)\|^2d\tau
+C\int_0^t\int\int\f{|\na^2_x\pa^{\alpha}G|^2}{\Mloc_*} dvdxd\tau\\
\nm &\quad +C(\d+\v^{\f12})\int_0^t(\|\na\pa^{\alpha}(n,w,T)\|^2+\|\na(u_t,\t_t)\|^2+\|\na^2(w-u,T-\t)\|^2
+\|\pa^{\alpha}(\tilde\phi,\tilde\psi,\tilde\z)\|^2)d\tau\\
&\quad +C(\d+\v)\int_0^t\int\int\f{|(\na_v\pa^{\alpha}G,\na_v\wt{G})|^2}{\Mloc_*} dvdxd\tau,
\end{align}
where $C>0$ is a generic constant independent of time $\mathcal{T}>0$, and $\pa^{\alpha}$ denotes $\pa_t$ or $\na_x$.
\end{lemma}

\noindent {\it Proof.} First, applying $\pa_i~(i=1,2,3)$ to \eqref{G-1x}, then multiplying the resulting equation by $\pa_i\pa^{\alpha}{G}/\Mloc_*$, and integrating it over $[0,t]\times\Omega\times\mathbb{R}^3$, we obtain the desired inequality \eqref{eG-2x} similarly to Lemma \ref{le3.3}; see also Lemma 4.6 in \cite{LiWangWang2022}.  The proof of Lemma \ref{le3.6} is complete.

\hfill $\Box$

\

Finally, to obtain the third-order derivatives (in $x$ and/or $t$) for $G$, as well as the third-order spatial derivative estimates for the densities $\na^3(\r,n)$, we need to work with the original full compressible Euler--VFP system \eqref{EVFP}.

\begin{lemma}\label{le3.7}
Under the same assumptions as in Proposition \ref{Prop3.3}, the following estimate holds for $t\in[0,\mathcal{T}]:$
\begin{align}\label{es-all} \nm
& \|\na^2\pa^{\alpha}(\r,u,\t)(t)\|^2+\int\int\f{|\na_x^2\pa^{\alpha}f(t)|^2}{\Mloc_*}dv dx
+\int_0^t\|\na^2\pa^{\alpha}(w-u,T-\t)\|^2d\tau\\ \nm
&
+\int_0^t\int\int\f{|\na_v\na_x^2\pa^{\alpha}G|^2+(1+|v|^2)|\na_x^2\pa^{\alpha}G|^2}{\Mloc_*}dvdxd\tau\\
&\leq C(\mathcal{N}(0)+\v)
+C(\d+\v^{\f12}+\eta_0)\int_0^t\|\na^2\pa^{\alpha}(\r,n,w,T)\|^2d\tau\\ \nm
&\quad +C(\d+\v^{\f12})\int_0^t(\|\na\pa^{\alpha}(\r,n,w,T)\|^2+\|\na(u_t,\t_t)\|^2+\|(\psi_t,\z_t)\|^2)d\tau\\ \nm
&\quad +C(\d+\v^{\f12})\int_0^t(\|\na(w-u,T-\t)\|_{H^1}^2+\|\pa^{\alpha}(\tilde\phi,\tilde\psi,\tilde\z)\|^2)d\tau\\ \nm
&\quad +C(\d+\v)\int_0^t\int\int\f{|(\na_v\na_x\pa^{\alpha}G,\na_v\pa^{\alpha}G,\na_v\wt{G})|^2}{\Mloc_*}dvdxd\tau,
\end{align}
where $C>0$ is a generic constant independent of time $\mathcal{T}>0$ and $\pa^{\alpha}$ denotes $\pa_t$ or $\na_x$.
\end{lemma}

\noindent {\it Proof.}  Applying $\pa^{\alpha}$ \eqref{2-e-xx} yields
\begin{align}\label{2-e-xxt}
\begin{cases}
 & \pa_j\pa_i\pa^{\alpha}\r_t+u\cdot\na\pa_j\pa_i\pa^{\alpha}\r+\r\pa_j\pa_i\pa^{\alpha}\dv u
+\pa^{\alpha}u\cdot\na\pa_j\pa_i\r+\pa^{\alpha}\r\pa_j\pa_i\dv u\\
 &\quad+\pa^{\alpha}\big(
\pa_j u\cdot\na\pa_i\r+\pa_i u\cdot\na\pa_j\r+\dv u\pa_j\pa_i\r+\pa_j\r\pa_i\dv u+\pa_j\pa_i u\cdot\na\r+\pa_i\r\pa_j\dv u
\big)=0,\\[2mm]
 & \pa_j\pa_i\pa^{\alpha} u_t+u\cdot\na\pa_j\pa_i\pa^{\alpha} u+\f{\t}{\r}\na\pa_j\pa_i\pa^{\alpha}\r
+\na\pa_j\pa_i\pa^{\alpha}\t=-\pa^{\alpha}u\cdot\na\pa_j\pa_i u-\pa^{\alpha}\big(\f{\t}{\r}\big)\na\pa_j\pa_i\r\\
 &\quad-\pa^{\alpha}\Big(\pa_j u\cdot\na\pa_i u+\pa_j\pa_i u\cdot\na u+\pa_i u\cdot\na\pa_j u
+\pa_j\big(\f{\t}{\r}\big)\na\pa_i\r+\pa_j\pa_i\big(\f{\t}{\r}\big)\na\r+\pa_i\big(\f{\t}{\r}\big)\na\pa_j\r\Big)\\
 &\quad+\pa^{\alpha}\Big(\pa_j\big(\pa_i n(w-u)\big)+\pa_j n \pa_i(w-u)\Big)
+n\pa_j\pa_i\pa^{\alpha}(w-u)+\pa^{\alpha} n\pa_j\pa_i(w-u),\\[3mm]
& \f32\pa_j\pa_i\pa^{\alpha}\t_t+\f32u\cdot\na\pa_j\pa_i\pa^{\alpha}\t+\t\pa_j\pa_i\pa^{\alpha}\dv u
=-\f32\pa^{\alpha}u\cdot\na\pa_j\pa_i\t-\pa^{\alpha}\t\pa_j\pa_i\dv u\\
 &\quad-\pa^{\alpha}\Big(
\f32\pa_j u\cdot\na\pa_i\t+\f32\pa_j\pa_i u\cdot\na\t+\f32\pa_i u\cdot\na\pa_j\t
+\pa_j\t\pa_i\dv u+\pa_j\pa_i\t\dv u+\pa_i\t\pa_j\dv u\Big)
\\
 &\quad+\pa^{\alpha}\Big(\pa_j\big(3\pa_i n(T-\t)\big)
 +3\pa_j n \pa_i(T-\t)+\pa_j\pa_i\big(n|w-u|^2\big)\Big)
+3n\pa_j\pa_i\pa^{\alpha}(T-\t)\\
 &\quad+3\pa^{\alpha}n\pa_j\pa_i(T-\t).
\end{cases}
\end{align}
Multiplying \eqref{2-e-xxt}$_1$, \eqref{2-e-xxt}$_2$, and \eqref{2-e-xxt}$_3$ by $\f{\t}{\r T}\pa_j\pa_i\pa^{\alpha}\r$,
$\f{\r}{T}\pa_j\pa_i\pa^{\alpha} u$, and $\f{\r}{\t T}\pa_j\pa_i\pa^{\alpha}\t$, respectively,
then adding the resulting equations and integrating it by parts over $[0,t]\times\Omega\times\mathbb{R}^3$, we obtain
\begin{align}\label{u-xxt}\nm
&\quad\int\Big(\f{\t}{\r T}\f{(\pa_j\pa_i\pa^{\alpha}\r)^2}{2}+\f{\r}{T}\f{|\pa_j\pa_i\pa^{\alpha}u|^2}{2}
+\f{3\r}{2\t T}\f{(\pa_j\pa_i\pa^{\alpha}\t)^2}{2}\Big)dx\Big|_{\tau=0}^{\tau=t}\\\nm
&-\int_0^t\int \f{\r n}{T}\pa_j\pa_i\pa^{\alpha}(w-u)\cdot\pa_j\pa_i\pa^{\alpha}u \, dxd\tau
-\int_0^t\int\f{3\r n}{\t T}\pa_j\pa_i\pa^{\alpha}(T-\t)\pa_j\pa_i\pa^{\alpha}\t \, dxd\tau\\\nm
&
=\int_0^t\int\Big(
\big((\f{\t}{\r T})_t+\dv(\f{ u \t}{\r T})\big)\f{(\pa_j\pa_i\pa^{\alpha}\r)^2}{2}
+\big((\f{\r}{ T})_t+\dv(\f{\r u}{ T})\big)\f{(\pa_j\pa_i\pa^{\alpha} u)^2}{2}\\
&\quad +\f32\big(\big(\f{\r}{\t T}\big)_t+\dv\big(\f{\r u}{\t T}\big)\big)\f{(\pa_j\pa_i\pa^{\alpha}\t)^2}{2}
+\pa_j\pa_i\pa^{\alpha}\pa^{\alpha}u\cdot\big[\na\big(\f{\t}{T}\big)\r\pa_j\pa_i\pa^{\alpha}\r+\na\big(\f{\r}{T}\big)\pa_j\pa_i\pa^{\alpha}\t\big]\\ \nm
&\quad +\big(-\pa^{\alpha}u\cdot\na\pa_j\pa_i u-\pa^{\alpha}\big(\f{\t}{\r}\big)\na\pa_j\pa_i\r+\pa^{\alpha} n \pa_j\pa_i(w-u)
+\pa^{\alpha}(\cdots)\big)\f{\r}{T}\pa_j\pa_i\pa^{\alpha}u\\ \nm
&\quad +\big(-\f32\pa^{\alpha}u\cdot\na\pa_j\pa_i\t-\pa^{\alpha}\t\pa_j\pa_i\dv u+3\pa^{\alpha} n\pa_j\pa_i(T-\t)+\pa^{\alpha}(\cdots)\big)
\f{\r}{\t T}\pa_j\pa_i\pa^{\alpha}\t\\ \nm
&\quad -\big(\pa^{\alpha}u\cdot\na\pa_j\pa_i\r+\pa^{\alpha}\r\pa_j\pa_i\dv u+\pa^{\alpha}(\cdots)\big)\f{\t}{\r T}\pa_j\pa_i\pa^{\alpha}\r
\Big)dxd\tau.
\end{align}
The terms on the right-hand side of \eqref{u-xxt} can be controlled by using \eqref{assum-1}-\eqref{assum-3}, similarly to Lemma \ref{le3.4}. We illustrate the argument with one typical triple nonlinear term, for instance, the first term in $\pa^{\alpha}(\cdots)\f{\t}{\r T}\pa_j\pa_i\pa^{\alpha}\r$. By Sobolev inequality and assumption \eqref{assumption}, we have
\begin{align}\label{typical} \nm
\nm&\int_0^t\int -\pa^{\alpha}\pa_j u\cdot\na\pa_i\r\f{\t}{\r T}\pa_j\pa_i\pa^{\alpha}\r\, dxd\tau
\leq C\int_0^t\|\na^2\pa^{\alpha}\r\|\|\na\pa^{\alpha} u\|_{L^4}\|\na^2\r\|_{L^4}d\tau\\
& \leq C\int_0^t\|\na^2\pa^{\alpha}\r\|\|\na\pa^{\alpha} u\|^{\f14}\|\na\pa^{\alpha} u\|^{\f34}_{H^1}\|\na^2\r\|^{\f14}\|\na^2\r\|^{\f34}_{H^1}d\tau\\
\nm & \leq C\int_0^t\|\na^2\pa^{\alpha}\r\|\|\na\pa^{\alpha} u\|\|\na^2\r\|d\tau
+C\int_0^t\|\na^2\pa^{\alpha}\r\|\|\na\pa^{\alpha} u\|^{\f14}\|\na^2\pa^{\alpha} u\|^{\f34}
\|\na^2\r\|^{\f14}\|\na^3\r\|^{\f34}d\tau\\
\nm &\quad+ C\int_0^t\|\na^2\pa^{\alpha}\r\|\|\na\pa^{\alpha} u\|\|\na^2\r\|^{\f14}\|\na^3\r\|^{\f34}d\tau
+C\int_0^t\|\na^2\pa^{\alpha}\r\|\|\na\pa^{\alpha} u\|^{\f14}\|\na^2\pa^{\alpha} u\|^{\f34}
\|\na^2\r\|d\tau\\
\nm &
 \leq C(\d+\v)\int_0^t(\|\na^2\pa^{\alpha}(\r,u)\|^2+\|\na\pa^{\alpha}u\|^2)d\tau,
%& \leq C(\d+\v)\int_0^t(\|\na^2\pa^{\alpha}(\r,w)\|^2+\|\na^2\pa^{\alpha}(w-u)\|^2+\|\na u_t\|^2+\|\na^2(w,w-u)\|^2)d\tau.
\end{align}
where in the last inequality we have used  $\|\na^2\r\|_{H^1}\leq C(\d+\v)$.

Summing \eqref{u-xxt} over $i,j=1,2,3$ and combining the result with  \eqref{typical}, we have
\begin{align}\label{2xt}
\nm &\quad \int\Big(\f{\t}{\r T}\f{|\na^2\pa^{\alpha}\r|^2}{2}+\f{\r}{T}\f{|\na^2\pa^{\alpha}u|^2}{2}
+\f{3\r}{2\t T}\f{|\na^2\pa^{\alpha}\t|^2}{2}\Big)dx\Big|_{\tau=t}
\\
&-\int_0^t\int \Big(\f{\r n}{T}\na^2\pa^{\alpha}(w-u)\cdot\na^2\pa^{\alpha}u
+\f{3\r n}{\t T}\na^2\pa^{\alpha}(T-\t)\cdot\na^2\pa^{\alpha}\t\Big) dxd\tau\\
\nm &\leq C\|\na^2\pa^{\alpha}(\r,u,\t)(0)\|^2
+C(\d+\v)\int_0^t(\|\na^2\pa^{\alpha}(\r,u,\t,n)\|^2+\|\na\pa^{\alpha}(\r,u,\t,n)\|^2)d\tau
\\
\nm &\quad+C(\d+\v)\int_0^t\|\pa^{\alpha}(w-u,T-\t)\|_{H^1}^2d\tau.
\end{align}
The initial data satisfy
$$
\ba
\|\na^2(\r_t,u_t,\t_t)(0)\|^2&\leq C(\|\na^3(\r_0,u_0,\t_0)\|^2+\|\na^2(w_0-u_0,T_0-\t_0)\|^2)\\
&+C(\d+\v)(\|\na^2(\r_0,u_0,\t_0)\|^2+\|\na(w_0-u_0,T_0-\t_0)\|^2).
\ea
$$
Moreover,
\begin{align}\nm
\|\na\pa^{\alpha}(u,\t)\|&=\|\na^2(u,\t)\|+\|\na(u_t,\t_t)\|\\ \nm
&\leq \|\na^2(w,T)\|+\|\na^2(w-u,T-\t)\|+\|\na(u_t,\t_t)\|,\\ \nm
\|\pa^{\alpha}(w-u,T-\t)\|_{H^1}&=\|\na\pa^{\alpha}(w,T)\|+\|\na\pa^{\alpha}(u,\t)\|
+\|\pa^{\alpha}(w-u,T-\t)\|\\ \nm
&\leq \|\na\pa^{\alpha}(w,T)\|+\|\na(w-u,T-\t)\|_{H^1}+\|\na(u_t,\t_t)\|
+\|\pa_t(\psi,\z,\tilde\psi,\tilde\z)\|,
\end{align}
which together with \eqref{2xt} leads to
\begin{align}\label{2xt-re}
\nm&\quad \int\Big(\f{\t}{\r T}\f{|\na^2\pa^{\alpha}\r|^2}{2}+\f{\r}{T}\f{|\na^2\pa^{\alpha}u|^2}{2}
+\f{3\r}{2\t T}\f{|\na^2\pa^{\alpha}\t|^2}{2}\Big)dx\Big|_{\tau=t}
\\
\nm &-\int_0^t\int \Big(\f{\r n}{T}\na^2\pa^{\alpha}(w-u)\cdot\na^2\pa^{\alpha}u
+\f{3\r n}{\t T}\na^2\pa^{\alpha}(T-\t)\cdot\na^2\pa^{\alpha}\t\Big) dxd\tau\\
&\leq C(\mathcal{N}(0)+\v)
+C(\d+\v)\int_0^t(\|\na^2\pa^{\alpha}(\r,n,w,T)\|^2+\|\na\pa^{\alpha}(\r,n,w,T)\|^2)d\tau
\\
\nm &\quad+C(\d+\v)\int_0^t(\|\na^2\pa^{\alpha}(w-u,T-\t)\|^2+\|\na(w-u,T-\t)\|_{H^1}^2)d\tau\\
\nm &\quad+C(\d+\v)\int_0^t(\|\na(u_t,\t_t)\|^2+\|\pa_t(\psi,\z,\tilde\psi,\tilde\z)\|^2)d\tau.
\end{align}

Next, we rewrite \eqref{EVFP}$_4$ as
$$
\begin{array}{l}
\di\quad f_t+v\cdot\na_x f=\r\Lop_{[w,T]} f+\r(w-u)\cdot\na_v f-\r(T-\t)\Delta_v f\\
=\r\Lop_{[w,T]}G+\r(w-u)\cdot\na_vG+\r(w-u)\cdot\na_v\Mloc-\r(T-\t)\Delta_vG
-\r(T-\t)\Delta_v\Mloc.
\end{array}
$$
Applying $\pa_j\pa_i\pa^{\alpha}$ to the above equation, we have
\begin{align}\label{f-3x} \nm
& \pa_j\pa_i\pa^{\alpha}f_t+v\cdot\na_x\pa_j\pa_i\pa^{\alpha}f
=\r\Lop_{[w,T]}(\pa_j\pa_i\pa^{\alpha}G)\\ \nm
&-\r\pa^{\alpha}u\cdot\na_v\pa_j\pa_iG-\r(\pa_j u\cdot\na_v\pa_i\pa^{\alpha}G+\pa_i u\cdot\na_v\pa_j\pa^{\alpha}G)\\ \nm
&-\r\pa_j\pa_i u\cdot\na_v\pa^{\alpha}G
-\r(\pa_j\pa^{\alpha}u\cdot\na_v\pa_iG+\pa_i\pa^{\alpha}u\cdot\na_v\pa_jG)
-\r\pa_j\pa_i\pa^{\alpha}u\cdot\na_vG\\ \nm
&+\r\pa^{\alpha}(w-u)\cdot\na_v\pa_j\pa_i\Mloc
+\r(\pa_j(w-u)\cdot\na_v\pa_i\pa^{\alpha}\Mloc+\pa_i(w-u)\cdot\na_v\pa_j\pa^{\alpha}\Mloc)\\\nm
&+\r\pa_j\pa_i(w-u)\cdot\na_v\pa^{\alpha}\Mloc
+\r(\pa_j\pa^{\alpha}(w-u)\cdot\na_v\pa_i\Mloc+\pa_i\pa^{\alpha}(w-u)\cdot\na_v\pa_j\Mloc)\\
&+\r\pa_j\pa_i\pa^{\alpha}(w-u)\cdot\na_v\Mloc+\r(w-u)\cdot\na_v\pa_j\pa_i\pa^{\alpha}f\\ \nm
&+\r\pa^{\alpha}\t\Delta_v\pa_j\pa_iG+\r(\pa_j \t\Delta_v\pa_i\pa^{\alpha}G+\pa_i \t\Delta_v\pa_j\pa^{\alpha}G)\\ \nm
&+\r\pa_j\pa_i\t\Delta_v\pa^{\alpha}G
+\r(\pa_j\pa^{\alpha}\t\Delta_v\pa_iG+\pa_i\pa^{\alpha}\t\Delta_v\pa_jG)
+\r\pa_j\pa_i\pa^{\alpha}\t\Delta_vG\\ \nm
&-\r\pa^{\alpha}(T-\t)\na_v\pa_j\pa_i\Mloc
-\r(\pa_j(T-\t)\Delta_v\pa_i\pa^{\alpha}\Mloc+\pa_i(T-\t)\Delta_v\pa_j\pa^{\alpha}\Mloc)\\ \nm
&-\r\pa_j\pa_i(T-\t)\Delta_v\pa^{\alpha}\Mloc
-\r(\pa_j\pa^{\alpha}(T-\t)\Delta_v\pa_i\Mloc+\pa_i\pa^{\alpha}(T-\t)\Delta_v\pa_j\Mloc)\\ \nm
&-\r\pa_j\pa_i\pa^{\alpha}(T-\t)\Delta_v\Mloc-\r(T-\t)\Delta_v\pa_j\pa_i\pa^{\alpha}f
+(\cdots)\\ \nm
&=:\r\Lop_{[w,T]}(\pa_j\pa_i\pa^{\alpha}G)
+\sum_{i=1}^{11}\tilde K_i+\sum_{i=1}^{11}\tilde L_i+(\cdots),
\end{align}
where $(\cdots)$ denotes the terms containing derivatives of the density $\r$.

Multiplying \eqref{f-3x} by $\pa_j\pa_i\pa^{\alpha}f/\Mloc_*$
and integrating it by parts over $[0,t]\times\Omega\times\mathbb{R}^3$, we obtain
\begin{align}\label{3f}
&\quad\int\int\f{|\pa_j\pa_i\pa^{\alpha}f|^2}{2\Mloc_*}dv dx\Big|_{\tau=0}^{\tau=t}
-\int_0^t\int\int\r \f{\pa_j\pa_i\pa^{\alpha}G}{\Mloc_*}\Lop_{[w,T]}(\pa_j\pa_i\pa^{\alpha}G)dvdxd\tau\\ \nm
&=:\int_0^t\int\int\r\Lop_{[w,T]}(\pa_j\pa_i\pa^{\alpha}G)\f{\pa_j\pa_i\pa^{\alpha}\Mloc}{\Mloc_*}dvdxd\tau
+\sum_{i=1}^{11} K_i+\sum_{i=1}^{11} L_i+\int_0^t\int\int  (\cdots)\f{\pa_j\pa_i\pa^{\alpha}f}{\Mloc_*}\, dv dx d\tau,
\end{align}
where
$$
K_i:=\int_0^t\int\int\tilde  K_i\f{\pa_j\pa_i\pa^{\alpha}f}{\Mloc_*}\, dv dx d\tau, \qquad
L_i:=\int_0^t\int\int\tilde  L_i\f{\pa_j\pa_i\pa^{\alpha}f}{\Mloc_*}\, dv dx d\tau.
$$
Since $\Mloc_t$ and $\na_x\Mloc$ belong to $\mathfrak{N}$, the term $\mb{P}_1(\pa_j\pa_i\pa^{\alpha}\Mloc)$ does not contain $\pa_j\pa_i\pa^{\alpha}(n,w,T)$. By the assumption \eqref{assum-1}, we have
\begin{align} \nm
& \quad \int_0^t\int\int\r\Lop_{[w,T]}(\pa_j\pa_i\pa^{\alpha}G)\f{\pa_j\pa_i\pa^{\alpha}\Mloc}{\Mloc_*}dvdxd\tau\\ \nm
&= \int_0^t\int\int \r\Lop_{[w,T]}(\pa_j\pa_i\pa^{\alpha}G)\f{\mb{P}_1(\pa_j\pa_i\pa^{\alpha}\Mloc)}{\Mloc}dvdxd\tau
\\ \nm
&\quad
+\int_0^t\int\int\r\Lop_{[w,T]}(\pa_j\pa_i\pa^{\alpha}G)
\Big(\f{1}{\Mloc_*}-\f{1}{\Mloc}\Big)\pa_j\pa_i\pa^{\alpha}\Mloc dvdxd\tau\\ \nm
&=-\int_0^t\int\int\r\big((v-w)\pa_j\pa_i\pa^{\alpha}G+T\na_v\pa_j\pa_i\pa^{\alpha}G\big)
\cdot\na_v\Big(\f{\mb{P}_1(\pa_j\pa_i\pa^{\alpha}\Mloc)}{\Mloc}\Big)dvdxd\tau\\ \nm
&\quad-\int_0^t\int\int\r\big((v-w)\pa_j\pa_i\pa^{\alpha}G+T\na_v\pa_j\pa_i\pa^{\alpha}G\big)
\cdot\na_v\Big[\Big(\f{1}{\Mloc_*}-\f{1}{\Mloc}\Big)\pa_j\pa_i\pa^{\alpha}\Mloc\Big]dvdxd\tau\\ \nm
& \leq C(\d+\v+\eta_0)\int_0^t\int\int\f{|\na_v\na_x^2\pa^{\alpha}G|^2
+(1+|v|^2)|\na_x^2\pa^{\alpha}G|^2}{\Mloc_*}dvdxd\tau+C\v\\ \nm
&\quad +C\eta_0\int_0^t\|\na^2\pa^{\alpha}(n,w,T)\|^2d\tau+C(\d+\v)\int_0^t(\|\na\pa^{\alpha}(n,w,T)\|^2
+\|\na(\tilde\phi,\tilde\psi,\tilde\z)\|^2)d\tau.
\end{align}
From \eqref{assum-1}, we have $\|\pa^{\alpha}u\|_{L^{\infty}}\leq C(\d+\v)$. By Sobolev inequality, it follows that
\begin{align}
K_1&\leq \nm C\int_0^t\|\pa^{\alpha}u\|_{L^{\infty}}\Big(\int\int\f{|\na_v\na_x^2G|^2}{\Mloc_*}dvdx\Big)^{\f12}
\cdot\Big[\|\na^2\pa^{\alpha}(n,w,T)\|\\ \nm
&\quad+C(\d+\v)(\|\na\pa^{\alpha}(n,w,T)\|+\|\na(n,w,T)\|^2_{L^4})
+\Big(\int\int\f{|\na_x^2\pa^{\alpha}G|^2}{\Mloc_*}dvdx\Big)^{\f12}\Big]d\tau\\ \nm
& \leq C(\d+\v)\int_0^t(\|\na^2\pa^{\alpha}(n,w,T)\|^2+\|\na\pa^{\alpha}(n,w,T)\|^2+\|\na(\tilde\phi,\tilde\psi,\tilde\z)\|^2)d\tau+C\v
\\ \nm
&\quad
+C(\d+\v)\int_0^t\int\int\f{|(\na^2_x\pa^{\alpha}G,\na_v\na^2_xG)|^2}{\Mloc_*}dvdxd\tau.
\end{align}
Likewise,
\begin{align} \nm
K_2& \leq C(\d+\v)\int_0^t(\|\na^2\pa^{\alpha}(n,w,T)\|^2+\|\na\pa^{\alpha}(n,w,T)\|^2+\|\na(\tilde\phi,\tilde\psi,\tilde\z)\|^2)d\tau+C\v
\\ \nm
&\quad
+C(\d+\v)\int_0^t\int\int\f{|(\na^2_x\pa^{\alpha}G,\na_v\na_x\pa^{\alpha}G)|^2}{\Mloc_*}dvdxd\tau.
\end{align}
From \eqref{assumption}, we have $\|\na^2u\|\leq C(\d+\v)$. By Cauchy inequality, we have
\begin{align}
K_3&\leq  \nm C\int_0^t\|\na^2u\|\Big\|\int\f{|\na_v\pa^{\alpha}G|^2}{\Mloc_*}dv\Big\|^{\f12}_{L_x^{\infty}}
\cdot\Big[\|\na^2\pa^{\alpha}(n,w,T)\|\\  \nm
&\quad+C(\d+\v)(\|\na\pa^{\alpha}(n,w,T)\|+\|\na(n,w,T)\|^2_{L^4})
+\Big(\int\int\f{|\na_x^2\pa^{\alpha}G|^2}{\Mloc_*}dvdx\Big)^{\f12}\Big]d\tau\\ \nm
&\leq C(\d+\v)\int_0^t\Big(\int\int\f{|(\na_v\pa^{\alpha}G,\na_v\na_x\pa^{\alpha}G,\na_v\na^2_x\pa^{\alpha}G)|^2}{\Mloc_*}dvdx\Big)^{\f12}
\cdot\Big[\|\na^2\pa^{\alpha}(n,w,T)\|\\ \nm
&\quad+C(\d+\v)(\|\na\pa^{\alpha}(n,w,T)\|+\|\na(n,w,T)\|^2_{L^4})
+\Big(\int\int\f{|\na_x^2\pa^{\alpha}G|^2}{\Mloc_*}dvdx\Big)^{\f12}\Big]d\tau\\ \nm
&\leq C(\d+\v)\int_0^t(\|\na^2\pa^{\alpha}(n,w,T)\|^2+\|\na\pa^{\alpha}(n,w,T)\|^2+\|\na(\tilde\phi,\tilde\psi,\tilde\z)\|^2)d\tau+C\v\\ \nm
&\quad+C(\d+\v)\int_0^t\int\int\f{|(\na^2_x\pa^{\alpha}G,\na_v\na^2_x\pa^{\alpha}G,\na_v\na_x\pa^{\alpha}G,\na_v\pa^{\alpha}G)|^2}{\Mloc_*}dvdxd\tau.
\end{align}
From \eqref{1-ns-x}$_2$ and the assumption \eqref{assumption}, we deduce that $\|\na\pa^{\alpha} u\|\leq C(\d+\v)$. Similarly to $K_3$, we have
\begin{align}
K_4&\leq  \nm C\int_0^t\|\na\pa^{\alpha}u\|\Big\|\int\f{|\na_v\na_xG|^2}{\Mloc_*}dv\Big\|^{\f12}_{L_x^{\infty}}
\cdot\Big[\|\na^2\pa^{\alpha}(n,w,T)\|\\  \nm
&\quad+C(\d+\v)(\|\na\pa^{\alpha}(n,w,T)\|+\|\na(n,w,T)\|^2_{L^4})
+\Big(\int\int\f{|\na_x^2\pa^{\alpha}G|^2}{\Mloc_*}dvdx\Big)^{\f12}\Big]d\tau\\ \nm
&\leq C(\d+\v)\int_0^t(\|\na^2\pa^{\alpha}(n,w,T)\|^2+\|\na\pa^{\alpha}(n,w,T)\|^2+\|\na(\tilde\phi,\tilde\psi,\tilde\z)\|^2)d\tau+C\v\\ \nm
&\quad+C(\d+\v)\int_0^t\int\int\f{|(\na^2_x\pa^{\alpha}G,\na_v\na^3_xG,\na_v\na^2_xG,\na_v\na_xG)|^2}{\Mloc_*}dvdxd\tau,
\end{align}
where in the last inequality we have used
$$
\Big\|\int\f{|\na_v\na_xG|^2}{\Mloc_*}dv\Big\|_{L_x^{\infty}}
\leq C\int\int\f{|(\na_v\na_xG,\na_v\na_x^2G,\na_v\na^3_xG)|^2}{\Mloc_*}dvdx.
$$
From \eqref{assum-4}, we have $\|\na^2\pa^{\alpha}u\|\leq C(\d+\v)$, and
\begin{align} \nm
K_5& \leq C\int_0^t\|\na^2\pa^{\alpha}u\|\Big\|\int\f{|\na_vG|^2}{\Mloc_*}dv\Big\|^{\f12}_{L_x^{\infty}}
\cdot\Big[\|\na^2\pa^{\alpha}(n,w,T)\|\\ \nm
&\quad+C(\d+\v)(\|\na\pa^{\alpha}(n,w,T)\|+\|\na(n,w,T)\|^2_{L^4})
+\Big(\int\int\f{|\na_x^2\pa^{\alpha}G|^2}{\Mloc_*}dvdx\Big)^{\f12}\Big]d\tau\\ \nm
&\leq C(\d+\v^{\f12})\int_0^t(\|\na^2\pa^{\alpha}(n,w,T)\|^2+\|\na^2\pa^{\alpha}(w-u)\|^2
+\|\na\pa^{\alpha}(n,w,T)\|^2)d\tau\\ \nm
&\quad+C(\d+\v^{\f12})\int_0^t\|\na(\tilde\phi,\tilde\psi,\tilde\z)\|^2d\tau
+C(\d+\v^{\f12})\int_0^t\int\int\f{|\na^2_x\pa^{\alpha}G|^2}{\Mloc_*}dvdxd\tau\\ \nm
&\quad+C(\d+\v)\int_0^t\int\int\f{|(\na_v\na_x^2G,\na_v\na_xG,\na_v\wt{G})|^2}{\Mloc_*}dvdxd\tau+C\v,
\end{align}
where in the last inequality we have used
$$
\ba
\Big\|\int\f{|\na_vG|^2}{\Mloc_*}dv\Big\|_{L_x^{\infty}}
&\leq C\int\int\f{|(\na_vG,\na_v\na_xG,\na_v\na^2_xG)|^2}{\Mloc_*}dvdx\\
&\leq C\int\int\f{|(\na_v\wt{G},\na_v\na_xG,\na_v\na^2_xG)|^2}{\Mloc_*}dvdx
+C\|\bar u_{1x_1}\|^2_{L^2_{x_1}},
\ea
$$
and
$$
\|\na^2\pa^{\alpha}u\|\leq \|\na^2\pa^{\alpha}w\|+\|\na^2\pa^{\alpha}(w-u)\|.
$$
By \eqref{assum-1},
$$
\|\pa^{\alpha}(w-u)\|_{L^{\infty}}=\|\pa^{\alpha}(\tilde\psi-\psi)\|_{L^{\infty}}
\leq \|\pa^{\alpha}(\psi,\tilde\psi)\|\leq C(\d+\v).
$$
Similarly to $K_1$, we have
\begin{align}
K_6&\leq \nm C\int_0^t\|\pa^{\alpha}(w-u)\|_{L^{\infty}}(\|\na^2(n,w,T)\|+\|\na(n,w,T)\|_{L^4}^2)
\cdot\Big[\|\na^2\pa^{\alpha}(n,w,T)\|\\ \nm
&\quad+C(\d+\v)(\|\na\pa^{\alpha}(n,w,T)\|+\|\na(n,w,T)\|^2_{L^4})
+\Big(\int\int\f{|\na_x^2\pa^{\alpha}G|^2}{\Mloc_*}dvdx\Big)^{\f12}\Big]d\tau\\ \nm
& \leq C(\d+\v)\int_0^t(\|\na^2\pa^{\alpha}(n,w,T)\|^2+\|\na\pa^{\alpha}(n,w,T)\|^2+\|\na(\tilde\phi,\tilde\psi,\tilde\z)\|^2)d\tau+C\v
\\ \nm
&\quad
+C(\d+\v)\int_0^t\int\int\f{|\na^2_x\pa^{\alpha}G|^2}{\Mloc_*}dvdxd\tau.
\end{align}
Likewise,
\begin{align} \nm
K_7&\leq \nm C\int_0^t\|\na(w-u)\|_{L^{\infty}}(\|\na\pa^{\alpha}(n,w,T)\|+\||\na(n,w,T)||\pa^{\alpha}(n,w,T)|\|)\\ \nm
&\quad\cdot\Big[\|\na^2\pa^{\alpha}(n,w,T)\|
+C(\d+\v)(\|\na\pa^{\alpha}(n,w,T)\|+\|\na(n,w,T)\|^2_{L^4})
+\Big(\int\int\f{|\na_x^2\pa^{\alpha}G|^2}{\Mloc_*}dvdx\Big)^{\f12}\Big]d\tau\\ \nm
& \leq C(\d+\v)\int_0^t(\|\na^2\pa^{\alpha}(n,w,T)\|^2+\|\na\pa^{\alpha}(n,w,T)\|^2
+\|\pa^{\alpha}(\tilde\phi,\tilde\psi,\tilde\z)\|^2)d\tau
\\ \nm
&\quad
+C(\d+\v)\int_0^t\int\int\f{|\na^2_x\pa^{\alpha}G|^2}{\Mloc_*}dvdxd\tau+C\v.
\end{align}
Similarly to $K_3$, we have
\begin{align} \nm
 K_8& \leq C\int_0^t\|\pa^{\alpha}(n,w,T)\|_{L^{\infty}}\|\na^2(w-u)\|\cdot
\Big[\|\na^2\pa^{\alpha}(n,w,T)\|\\ \nm
&\quad+C(\d+\v)(\|\na\pa^{\alpha}(n,w,T)\|+\|\na(n,w,T)\|_{L^4}^2)
+\Big(\int\int\f{|\na_x^2\pa^{\alpha}G|^2}{\Mloc_*}dvdx\Big)^{\f12}\Big]d\tau\\ \nm
& \leq C(\d+\v)\int_0^t(\|\na^2\pa^{\alpha}(n,w,T)\|^2+\|\na\pa^{\alpha}(n,w,T)\|^2
+\|\na(\tilde\phi,\tilde\psi,\tilde\z)\|^2)d\tau+C\v\\  \nm
&\quad +C(\d+\v)\int_0^t\|\na^2(w-u)\|^2d\tau
+C(\d+\v)\int_0^t\int\int\f{|\na^2_x\pa^{\alpha}G|^2}{\Mloc_*}dvdxd\tau.
\end{align}
Likewise,
\begin{align} \nm
 K_9& \leq C\int_0^t\|\na(n,w,T)\|_{L^{\infty}}\|\na\pa^{\alpha}(w-u)\|\cdot
\Big[\|\na^2\pa^{\alpha}(n,w,T)\|\\ \nm
&\quad+C(\d+\v)(\|\na\pa^{\alpha}(n,w,T)\|+\|\na(n,w,T)\|_{L^4}^2)
+\Big(\int\int\f{|\na_x^2\pa^{\alpha}G|^2}{\Mloc_*}dvdx\Big)^{\f12}\Big]d\tau\\ \nm
& \leq C(\d+\v)\int_0^t(\|\na^2\pa^{\alpha}(n,w,T)\|^2+\|\na\pa^{\alpha}(n,w,T)\|^2
+\|\na(\tilde\phi,\tilde\psi,\tilde\z)\|^2)d\tau+C\v\\  \nm
&\quad +C(\d+\v)\int_0^t(\|\na^2(w-u)\|^2+\|\na u_t\|^2)d\tau
+C(\d+\v)\int_0^t\int\int\f{|\na^2_x\pa^{\alpha}G|^2}{\Mloc_*}dvdxd\tau.
\end{align}
To fully exploit the relaxation friction damping in the momentum exchange arising from the fluid-particle interactions, we need to handle  $K_{10}$ carefully. We decompose $K_{10}$ as
$$
\ba
K_{10}&=\int_0^t\int\int\r\pa_j\pa_i\pa^{\alpha}(w-u)\cdot\na_v\Mloc\f{\pa_j\pa_i\pa^{\alpha}\Mloc}{\Mloc_*}dvdxd\tau\\
&+\int_0^t\int\int\r\pa_j\pa_i\pa^{\alpha}(w-u)\cdot\na_v\Mloc\f{\pa_j\pa_i\pa^{\alpha}G}{\Mloc_*}dvdxd\tau
=:K_{10,1}+K_{10,2}.
\ea
$$
A direct computation yields
\be\label{J11-1}
\ba
K_{10,1}&
=-\int_0^t\int\r\pa_j\pa_i\pa^{\alpha}(w-u)\cdot\Big(\int\f{v-w}{T}\Big(\f{\Mloc}{\Mloc_*}-1\Big)\pa_j\pa_i\pa^{\alpha}\Mloc\,dv\Big)dxd\tau\\
&\quad-\int_0^t\int\r\pa_j\pa_i\pa^{\alpha}(w-u)\cdot\Big(\int\f{v-w}{T}\pa_j\pa_i\pa^{\alpha}\Mloc\,dv\Big)dxd\tau,
\ea
\ee
where
\begin{align}
\nm&-\int_0^t\int\r\pa_j\pa_i\pa^{\alpha}(w-u)\cdot\Big(\int\f{v-w}{T}\pa_j\pa_i\pa^{\alpha}\Mloc\,dv\Big)dxd\tau\\
\nm&\qquad\leq-\int_0^t\int  \f{\r n}{T}\pa_j\pa_i\pa^{\alpha}(w-u)\cdot\pa_j\pa_i\pa^{\alpha}w \,dxd\tau\\
\nm&\qquad
+C(\d+\v)\int_0^t(\|\na^2\pa^{\alpha}(w-u)\|^2+\|\na\pa^{\alpha}(n,w,T)\|^2+\|\na(\tilde\phi,\tilde\psi,\tilde\z)\|^2)d\tau+C\v,
\end{align}
which together with \eqref{J11-1} implies
\be\label{J11-M}
\ba
K_{10,1}
&\leq -\int_0^t\int  \f{\r n}{T}\pa_j\pa_i\pa^{\alpha}(w-u)\cdot\pa_j\pa_i\pa^{\alpha}w \,dxd\tau
+C(\d+\v+\eta_0)\int_0^t\|\na^2\pa^{\alpha}(w-u,n,w,T)\|^2d\tau\\
&+C(\d+\v)\int_0^t(\|\na\pa^{\alpha}(n,w,T)\|^2+\|\na(\tilde\phi,\tilde\psi,\tilde\z)\|^2)d\tau+C\v.
\ea
\ee
By Cauchy inequality, we have
\be\label{J11-G}
\ba
K_{10,2}&
=-\int_0^t\int \r\pa_j\pa_i\pa^{\alpha}(w-u)\cdot\Big(\int\f{v-w}{T}\Big(\f{\Mloc}{\Mloc_*}-1\Big)\pa_j\pa_i\pa^{\alpha}Gdv\Big)dxd\tau\\
&\quad-\int_0^t\int \r \pa_j\pa_i\pa^{\alpha}(w-u)\cdot\Big(\int\f{v-w}{T}\pa_j\pa_i\pa^{\alpha}Gdv\Big)dxd\tau\\
& \leq C\eta_0\int_0^t\|\na^2\pa^{\alpha}(w-u)\|^2d\tau
+C\eta_0\int_0^t\int\int\f{|\na^2_x\pa^{\alpha}G|^2}{\Mloc_*}dvdxd\tau,
\ea
\ee
where we have used
$$
-\int_0^t\int\r\pa_j\pa_i\pa^{\alpha}(w-u)\cdot\Big(\int\f{v-w}{T}\pa_j\pa_i\pa^{\alpha}Gdv\Big)dxd\tau=0.
$$
Combining \eqref{J11-M} and \eqref{J11-G} yields
$$
\ba
K_{10}&\leq -\int_0^t\int \f{\r n}{T}  \pa_j\pa_i\pa^{\alpha}(w-u)\cdot\pa_j\pa_i\pa^{\alpha}w \,dxd\tau
+C(\d+\v+\eta_0)\int_0^t\|\na^2\pa^{\alpha}(w-u,n,w,T)\|^2d\tau
\\
&+C(\d+\v)\int_0^t(\|\na\pa^{\alpha}(n,w,T)\|^2+\|\na(\tilde\phi,\tilde\psi,\tilde\z)\|^2)d\tau
+C\eta_0\int_0^t\int\int\f{|\na^2_x\pa^{\alpha}G|^2}{\Mloc_*}dvdxd\tau+C\v.
\ea
$$
Integrating by parts in $v$ over $\mathbb{R}^3$, we estimate $K_{11}$ as follows.
\begin{align} \nm
 K_{11}&=-\int_0^t\int\r (w-u)\cdot\Big(\int\na_v\Big(\f{1}{\Mloc_*}\Big)\f{(\pa_j\pa_i\pa^{\alpha}f)^2}{2}dv\Big)dxd\tau\\ \nm
&=-\int_0^t\int\r (w-u)\cdot\Big(\int\f{v-w_*}{T_*}\f{(\pa_j\pa_i\pa^{\alpha}f)^2}{2\Mloc_*}dv\Big)dxd\tau\\ \nm
&\leq C\d\int_0^t(\|\na^2\pa^{\alpha}(n,w,T)\|^2+\|\na\pa^{\alpha}(n,w,T)\|^2+\|\na(\tilde\phi,\tilde\psi,\tilde\z)\|^2)d\tau\\ \nm
&+C\d\int_0^t\int\int\f{1+|v|^2}{\Mloc_*}|\na^2_x\pa^{\alpha}G|^2dvdxd\tau+C\v.
\end{align}
The terms $L_i~(i=1,\cdots,9,11)$ and
$$
\int_0^t\int\int  (\cdots)\f{\pa_j\pa_i\pa^{\alpha}f}{\Mloc_*}\, dv dx d\tau
$$
can be estimated similarly to the corresponding $K_i~(i=1,\cdots,9,11)$. We therefore focus on $L_{10}$ in detail, in order to exploit the relaxation friction damping in the energy exchange. As for $K_{10}$, we decompose $L_{10}$ as
$$
\ba
L_{10}&=-\int_0^t\int\int\r\pa_j\pa_i\pa^{\alpha}(T-\t)\Delta_v\Mloc\f{\pa_j\pa_i\pa^{\alpha}\Mloc}{\Mloc_*}dvdxd\tau\\
&\quad-\int_0^t\int\int\r\pa_j\pa_i\pa^{\alpha}(T-\t)\Delta_v\Mloc\f{\pa_j\pa_i\pa^{\alpha}G}{\Mloc_*}dvdxd\tau
=:L_{10,1}+L_{10,2}.
\ea
$$
A direct computation yields
\be\label{L10-1}
\ba
L_{10,1}&=-\int_0^t\int\r\f{\pa_j\pa_i\pa^{\alpha}(T-\t)}{T}
\Big(\int\big(\f{|v-w|^2}{T}-3\big)\f{\Mloc}{\Mloc_*}\pa_j\pa_i\pa^{\alpha}\Mloc \,dv\Big)dxd\tau\\
&=-\int_0^t\int\r\f{\pa_j\pa_i\pa^{\alpha}(T-\t)}{T}\Big(\int\big(\f{|v-w|^2}{T}-3\big)
\Big(\f{\Mloc}{\Mloc_*}-1\Big)\pa_j\pa_i\pa^{\alpha}\Mloc\,dv\Big)dxd\tau\\
&\quad-\int_0^t\int\r\f{\pa_j\pa_i\pa^{\alpha}(T-\t)}{T}\Big(\int\big(\f{|v-w|^2}{T}-3\big)\pa_j\pa_i\pa^{\alpha}\Mloc\,dv\Big)dxd\tau,
\ea
\ee
where
$$
\ba
&-\int_0^t\int\r\f{\pa_j\pa_i\pa^{\alpha}(T-\t)}{T}\Big(\int\big(\f{|v-w|^2}{T}-3\big)\pa_j\pa_i\pa^{\alpha}\Mloc\,dv\Big)dxd\tau\\
&\qquad\leq-\int_0^t\int  \f{3\r n}{T^2}\pa_j\pa_i\pa^{\alpha}(T-\t)\pa_j\pa_i\pa^{\alpha}T dxd\tau\\
&\qquad+C(\d+\v)\int_0^t(\|\na^2\pa^{\alpha}(T-\t)\|^2+\|\na\pa^{\alpha}(n,w,T)\|^2+\|\na(\tilde\phi,\tilde\psi,\tilde\z)\|^2)d\tau+C\v,
\ea
$$
which together with \eqref{L10-1} implies
\begin{align}\label{L10-M} \nm
L_{10,1}&\leq -\int_0^t\int  \f{3\r n}{T^2}\pa_j\pa_i\pa^{\alpha}(T-\t)\pa_j\pa_i\pa^{\alpha}T \,dxd\tau
+C(\d+\v+\eta_0)\int_0^t\|\na^2\pa^{\alpha}(T-\t,n,w,T)\|^2d\tau\\
&\quad+C(\d+\v)\int_0^t(\|\na\pa^{\alpha}(n,w,T)\|^2+\|\na(\tilde\phi,\tilde\psi,\tilde\z)\|^2)d\tau+C\v.
\end{align}
Similarly to \eqref{J11-G}, we have
\be\label{L10-G}
\ba
L_{10,2}
& \leq C\eta_0\int_0^t\|\na^2\pa^{\alpha}(T-\t)\|^2d\tau
+C\eta_0\int_0^t\int\int\f{|\na^2_x\pa^{\alpha}G|^2}{\Mloc_*}dvdxd\tau.
\ea
\ee
Combining \eqref{L10-M} and \eqref{L10-G} yields
$$
\ba
L_{10}&\leq -\int_0^t\int  \f{3\r n}{T^2}\pa_j\pa_i\pa^{\alpha}(T-\t)\pa_j\pa_i\pa^{\alpha}T \,dxd\tau
+C(\d+\v+\eta_0)\int_0^t\|\na^2\pa^{\alpha}(T-\t,n,w,T)\|^2d\tau\\
&+C(\d+\v)\int_0^t(\|\na\pa^{\alpha}(n,w,T)\|^2+\|\na(\tilde\phi,\tilde\psi,\tilde\z)\|^2)d\tau
+C\eta_0\int_0^t\int\int\f{|\na^2_x\pa^{\alpha}G|^2}{\Mloc_*}dvdxd\tau+C\v.
\ea
$$
Substituting the above estimates for $K_{i}$ and $L_i~(i=1\cdots,11)$ into \eqref{3f}, using Lemma \ref{co-1}, and choosing $\d$, $\v$, and  $\eta_0$ suitably small, we obtain
\begin{align}\label{es-f}
\nm& \int\int\f{|\na_x^2\pa^{\alpha}f(t)|^2}{2\Mloc_*}dv dx
+\f{\bar\sigma}{2}\int_0^t\int\int\f{|\na_v\na_x^2\pa^{\alpha}G|^2+(1+|v|^2)|\na_x^2\pa^{\alpha}G|^2}{\Mloc_*}dvdxd\tau\\
&+\int_0^t\int \Big(\f{\r n}{T}  \na^2\pa^{\alpha}(w-u)\cdot\na^2\pa^{\alpha}w
+\f{3\r n}{T^2}\na^2\pa^{\alpha}(T-\t)\cdot\na^2\pa^{\alpha}T\Big) dxd\tau\\
\nm &\leq C\int\int\f{|\na_x^2\pa^{\alpha}f(0)|^2}{\Mloc_*}dv dx
 +C(\d+\v^{\f12}+\eta_0)\int_0^t\|\na^2\pa^{\alpha}(w-u,T-\t,n,w,T)\|^2d\tau\\
\nm &\quad +C(\d+\v^{\f12})\int_0^t(\|\na\pa^{\alpha}(n,w,T)\|^2+\|\na^2(w-u,T-\t)\|^2
+\|\na(u_t,\t_t)\|^2+\|\pa^{\alpha}(\tilde\phi,\tilde\psi,\tilde\z)\|^2)d\tau
\\
\nm &\quad +C(\d+\v)\int_0^t\int\int\f{|(\na_v\na_x\pa^{\alpha}G,\na_v\pa^{\alpha}G,\na_v\wt{G})|^2}{\Mloc_*}dvdxd\tau+C\v.
\end{align}
Adding \eqref{2xt-re} and \eqref{es-f}, and using
$$
\ba
&\quad\int_0^t\int \f{3\r n}{T^2}\na^2\pa^{\alpha}(T-\t)\cdot\na^2\pa^{\alpha}T
-\f{3\r n}{\t T}\na^2\pa^{\alpha}(T-\t)\cdot\na^2\pa^{\alpha} \t\, dxd\tau\\
& =\int_0^t\int \f{3\r n}{\t T}|\na^2\pa^{\alpha}(T-\t)|^2dxd\tau
+\int_0^t\int\big(\f{3\r n}{T^2}-\f{3\r n}{\t T}\big)\na^2\pa^{\alpha}(T-\t)\cdot\na^2\pa^{\alpha}T \,dxd\tau\\
&\geq \int_0^t\int \f{3\r n}{\t T}|\na^2\pa^{\alpha}(T-\t)|^2dxd\tau
-C(\d+\v)\int_0^t\|\na^2\pa^{\alpha}(T-\t,T)\|^2d\tau.
\ea
$$
This yields the desired inequality \eqref{es-all}. The proof of Lemma \ref{le3.7} is complete.

\hfill $\Box$

Proposition \ref{Prop3.3} follows from Lemmas \ref{le3.4}-\ref{le3.7}.

\noindent{{\emph{Proof of Proposition \ref{Prop3.3}}}}:
Multiplying \eqref{e-2x} by a sufficiently large constant $C$, adding the resulting inequality to \eqref{3-rho}, and choosing $\d$ and $\varepsilon$ suitably small, we obtain
\begin{align}\label{xx-2}
\nm &
\|\na^2(\r,u,\t,n,w,T)(t)\|^2+\int_0^t(\|\na^3(\r,n,w,T)\|^2+\|\na^2(w-u,T-\t)\|^2)d\tau
\\
& \leq C(\mathcal{N}(0)+\v)+C\|\na^3(\r,n)(t)\|^2
+C\int_0^t\|\na^3(w-u,T-\t)\|^2d\tau\\
\nm &\quad +C(\d+\v)\int_0^t(\|\na^2(\r,n,w,T)\|^2+\|\na(\phi,\tilde\phi,\tilde\psi,\tilde\z)\|^2+\|\na(w-u,T-\t)\|^2)d\tau
\\
\nm &
\quad +C(\d+\v)\int_0^t\int\int\f{|(\na_v\na^2_xG,\na_x\pa^{\alpha}G,
\na_v\na_xG,\pa^{\alpha}G,\na_v\wt{G},\wt{G})|^2}{\Mloc_*}dvdxd\tau\\
\nm &\quad +C\int_0^t\int\int\f{|\na_x^2\pa^{\alpha}G|^2}{\Mloc_*}dvdxd\tau.
\end{align}

Next, we estimate $\int_0^t\|\na^2(\r_t,u_t,\t_t,n_t,w_t,T_t)\|^2d\tau$. To this end, we rewrite \eqref{2-nw-xx}$_3$ with non-fluid diffusive effect as
\begin{align}\label{3.102} \nm
&\f32\pa_j\pa_i T_t+\f32w\cdot\na\pa_j\pa_i T+T\pa_j\pa_i\dv w+\f32\pa_jw\cdot\na\pa_i T
+\f32\pa_j\pa_i w\cdot\na T+\f32\pa_iw\cdot\na\pa_j T
\\ \nm
&+\pa_j T\pa_i\dv w+\pa_j\pa_i T\dv w+\pa_i T\pa_j\dv w
=-3\r\pa_j\pa_i(T-\t)-3\pa_j\r\pa_i(T-\t)-3\pa_j\big(\pa_i\r(T-\t)\big)\\
&+\pa_j\pa_i\Big(\f{1}{n}\big(-\int\f12|v|^2v\cdot\na_xG dv+w\cdot\int v\otimes v\cdot\na_x G dv\big)\Big).
\end{align}
Multiplying \eqref{2-e-xx}$_1$, \eqref{2-e-xx}$_2$, \eqref{2-e-xx}$_3$, \eqref{2-nw-xx}$_1$, \eqref{nw-2x-G}, and \eqref{3.102} by $\pa_j\pa_i\r_t$, $\pa_j\pa_i u_t$, $\pa_j\pa_i\t_t$, $\pa_j\pa_i n_t$, $\pa_j\pa_i w_t$, and $\pa_j\pa_i T_t$, respectively, then integrating the resulting equations by parts over $[0,t]\times\Omega$, we obtain
\begin{equation}\label{xxt}
\ba
&\int_0^t\|\na^2(\r_t, u_t, \t_t, n_t, w_t,T_t)\|^2 d\tau \leq C\int_0^t(\|\na^3(\r,n,w,T)\|^2+\|\na^2(w-u,T-\t)\|_{H^1}^2)d\tau
 \\
&\quad+C(\d+\v)\int_0^t(\|\na^2(\r,n,w,T)\|^2+\|\na(w-u,T-\t)\|^2+\|\na(\phi,\tilde\phi)\|^2)d\tau+C\v\\
&\quad+C\int_0^t\int\int\frac{|\na^3_xG|^2}{\Mloc_*} dv dxd\tau
+C(\d+\v)\int_0^t\int\int\f{|\na^2_xG|^2}{\Mloc_*}dvdxd\tau.
\ea
\end{equation}
Multiplying \eqref{xx-2} by a sufficiently large constant $C$, adding the resulting inequality to \eqref{xxt}, and choosing $\d$ and $\varepsilon$ suitably small, we obtain
\begin{align}\label{full-3x} \nm
&
\|\na^2(\r,u,\t,n,w,T)(t)\|^2
+\int_0^t(\|\na^2\pa^{\alpha}(\r,n,w,T)\|^2+\|\na^2(u_t,\t_t)\|^2+\|\na^2(w-u,T-\t)\|^2)d\tau\\ \nm
&\leq C(\mathcal{N}(0)+\v)+C\|\na^3(\r,n)(t)\|^2
+C\int_0^t\|\na^3(w-u,T-\t)\|^2d\tau+C\int_0^t\int\int\f{|\na_x^2\pa^{\alpha}G|^2}{\Mloc_*}dvdxd\tau\\
&\quad +C(\d+\v)\int_0^t(\|\na^2(\r,n,w,T)\|^2+\|\na(\phi,\tilde\phi,\tilde\psi,\tilde\z)\|^2+\|\na(w-u,T-\t)\|^2)d\tau\\ \nm
&\quad +C(\d+\v)\int_0^t\int\int\f{|(\na_v\na^2_xG,\na_x\pa^{\alpha}G,
\na_v\na_xG,\pa^{\alpha}G,\na_v\wt{G},\wt{G})|^2}{\Mloc_*}dvdxd\tau.
\end{align}
Multiplying \eqref{full-3x} by a sufficiently large constant $C$, adding the resulting inequality to \eqref{eG-2x}, and choosing $\d$ and $\varepsilon$ suitably small, we have
\begin{align}\label{full+G} \nm
&\|\na^2(\r,u,\t,n,w,T)(t)\|^2
+\int_0^t(\|\na^2\pa^{\alpha}(\r,n,w,T)\|^2+\|\na^2(u_t,\t_t)\|^2+\|\na^2(w-u,T-\t)\|^2)d\tau
\\ \nm
&+\int\int \f{|\na_x\pa^{\alpha}G(t)|^2}{\Mloc_*}dv dx
+\int_0^t\int\int\f{|\na_v\na_x\pa^{\alpha}G|^2+(1+|v|^2)|\na_x\pa^{\alpha}G|^2}{\Mloc_*}dvdxd\tau
\\
&\leq C(\mathcal{N}(0)+\v) +C\|\na^3(\r,n)(t)\|^2+C\int_0^t\|\na^3(w-u,T-\t)\|^2d\tau\\ \nm
&
\quad +C(\d+\v^{\f12})\int_0^t(\|\na\pa^{\alpha}(\r,n,w,T)\|^2+\|\na(u_t,\t_t)\|^2
+\|\na(w-u,T-\t)\|^2+\|\pa^{\alpha}(\phi,\tilde\phi,\tilde\psi,\tilde\z)\|^2
)d\tau\\ \nm
&\quad +C\int_0^t\int\int\f{|\na_x^2\pa^{\alpha}G|^2}{\Mloc_*}dvdxd\tau
+C(\d+\v)\int_0^t\int\int\f{|(\na_v\pa^{\alpha}G,\pa^{\alpha}G,\na_v\wt{G},\wt{G})|^2}{\Mloc_*}dvdxd\tau.
\end{align}
Multiplying \eqref{es-all} by a sufficiently large constant $C$, adding the resulting inequality to \eqref{full+G}, and choosing $\d$, $\v$, and $\eta_0$ suitably small, we obtain the desired inequality \eqref{high-2}. The proof of Proposition \ref{Prop3.3} is complete.

\hfill $\Box$

%%%%%%%%%%%%%%%%%%%%%%%%%%%%%%%%%%%%%%%%%%%%%%%%%%%%%%%%%%%%%%%%%%%%%%%%%%%%%%

\begin{appendix}
\section{Proof of Lemma \ref{coefficient}}
\label{Appen-1}
\noindent {\it Proof.} By the definition in  \eqref{basis} and \eqref{p0}, a direct computation yields
\begin{equation}\label{P0P1}
\mb{P}_0^{\Mloc_{[n,w,T]}}=\mb{P}_0^{\Mloc_{[1,w,T]}},\quad \mb{P}_1^{\Mloc_{[n,w,T]}}=\mb{P}_1^{\Mloc_{[1,w,T]}}.
\end{equation}
Using \eqref{P0P1}, a direct calculation yields
$$
\ba
\quad\mb{P}_1(v\cdot\na_x\Mloc)&=\mb{P}_1^{\Mloc_{[1,w,T]}}(v\cdot\na_x\Mloc_{[n,w,T]})\\
&=n\mb{P}_1^{\Mloc_{[1,w,T]}}\Big(\big(\f{|v-w|^2}{2T}-\f32\big)(v-w)\cdot\f{\na_xT}{T}\Mloc_{[1,w,T]}\Big)\\
&+n\mb{P}_1^{\Mloc_{[1,w,T]}}\Big(\f{(v-w)\cdot\na_xw\cdot(v-w)}{T}\Mloc_{[1,w,T]}\Big),
\ea
$$
where
$$
\mb{P}_1^{\Mloc_{[1,w,T]}}\Big(\big(\f{|v-w|^2}{2T}-\f32\big)(v-w)\cdot\f{\na_xT}{T}\Mloc_{[1,w,T]}\Big)
=\big(\f{|v-w|^2}{2T}-\f52\big)(v-w)\cdot\f{\na_xT}{T}\Mloc_{[1,w,T]},
$$
and
$$
\mb{P}_1^{\Mloc_{[1,w,T]}}\Big(\f{(v-w)\cdot\na_xw\cdot(v-w)}{T}\Mloc_{[1,w,T]}\Big)
=\Big(\f{(v-w)\cdot\na_xw\cdot(v-w)}{T}-\f13\dv_x w\f{|v-w|^2}{T}\Big)\Mloc_{[1,w,T]}.
$$
Combining the above identities gives
\be\label{P1-vm}
\ba
\mb{P}_1(v\cdot\na_x\Mloc)=:nh\Mloc_{[1,w,T]}+ng\Mloc_{[1,w,T]},
\ea
\ee
where
$$
h:=\big(\f{|v-w|^2}{2T}-\f52\big)(v-w)\cdot\f{\na_x T}{T},\quad
g:=\f{(v-w)\cdot\na_xw\cdot(v-w)}{T}-\f13\dv_x w\f{|v-w|^2}{T}.
$$
Recalling the definition of the linearized Fokker-Planck operator $\Lop_{[w,T]}$ in \eqref{L}, we have
\begin{equation}\label{Lh}
\ba
\Lop_{[w,T]}(h\Mloc_{[1,w,T]})&=T(\Delta_v h-\f{v-w}{T}\cdot\na_v h)\Mloc_{[1,w,T]}\\
&=T(-\f{3h}{T})\Mloc_{[1,w,T]}=-3h\Mloc_{[1,w,T]},
\ea
\end{equation}
where we have used  $\Delta_v h-\f{v-w}{T}\cdot\na_v h=-\f{3h}{T}$. Likewise,
\begin{equation}\label{Lg}
\ba
\Lop_{[w,T]}(g\Mloc_{[1,w,T]})&
=T(\Delta_v g-\f{v-w}{T}\cdot\na_v g)\Mloc_{[1,w,T]}=T(-\f{v-w}{T}\cdot\na_v g)\Mloc_{[1,w,T]}\\
&=T(-\f{2g}{T})\Mloc_{[1,w,T]}=-2g\Mloc_{[1,w,T]},
\ea
\end{equation}
where we have used  $\Delta_v g-\f{v-w}{T}\cdot\na_v g=-\f{v-w}{T}\cdot\na_v g=-\f{2g}{T}$.
From \eqref{Lh} and \eqref{Lg}, We deduce  that
\begin{equation}\label{L-ni}
\Lop^{-1}_{[w,T]}(h\Mloc_{[1,w,T]})=-\f13h\Mloc_{[1,w,T]},\quad
\Lop^{-1}_{[w,T]}(g\Mloc_{[1,w,T]})=-\f12g\Mloc_{[1,w,T]}.
\end{equation}
We now derive the explicit formula for the viscosity term in $\eqref{ns-2}_2$. From \eqref{P1-vm} and \eqref{L-ni}, the $i$-th component $(i=1,2,3)$ of the viscosity coefficient is given by
$$
\begin{array}{l}
\di\quad -\int v_i v\cdot\na_x\big(\r^{-1}\Lop^{-1}_{[w,T]}\big[\mb{P}_1(v\cdot\nabla_x\Mloc)\big]\big) dv\\
\di =\sum_{j=1}^3\pa_{x_j}\big[\f{n}{3\r}\int v_iv_j h \Mloc_{[1,w,T]}dv+\f{n}{2\r}\int v_iv_j g \Mloc_{[1,w,T]}dv\big]\\
\di =\sum_{j=1}^3\pa_{x_j}\big[\f{n}{2\r}\int (v_i-w_i)(v_j-w_j) g \Mloc_{[1,w,T]}dv\big]\\
\di =\sum_{j=1}^3\pa_{x_j}\big[\f{n T}{2\r}(\pa_{x_i}w_j+\pa_{x_j}w_i-\f23\d_{ij}\dv_x w)\big],
\end{array}
$$
and the explicit formula for the heat conductivity coefficient is given by
$$
\begin{array}{l}
\di\quad -\int \f{1}{2}|v|^2 v\cdot\na_x\big(\r^{-1}\Lop^{-1}_{[w,T]}\big[\mb{P}_1(v\cdot\nabla_x\Mloc)\big]\big) dv\\
\di =\sum_{j=1}^3\pa_{x_j}\big(\f{n}{3\r}\int \f{1}{2}|v|^2v_j h \Mloc_{[1,w,T]}dv+\f{n}{2\r}\int \f{1}{2}|v|^2v_j g \Mloc_{[1,w,T]}dv\big)\\
\di =\sum_{j=1}^3\pa_{x_j}\big(\f{n}{3\r}\int \f12|v-w|^2(v_j-w_j) h \Mloc_{[1,w,T]}dv
+\f{n}{2\r} w\cdot\int(v-w)(v_j-w_j) g \Mloc_{[1,w,T]}dv\big)\\
\di =\sum_{j=1}^3\pa_{x_j}\big(\f{n}{3\r}\f52T\pa_{x_j}T+\f{n}{2\r}T\sum_{i=1}^3w_i(\pa_{x_j}w_i+\pa_{x_i}w_j-\f23\d_{ij}\dv_x w)\big).
\end{array}
$$
This yields \eqref{viscosity} and \eqref{heat-conduct}, and completes the proof of Lemma \ref{coefficient}.

\hfill $\Box$

%%%%%%%%%%%%%%%%%%%%%%%%%%%%%%%%%%%%%%%%%%%%%%%%%%%%%%%%%%%%%%%%%%%%%%%%%%%%%%
\section{Proof of Lemma \ref{le-diss}}
\label{Appen-2}
\noindent {\it Proof.} By \eqref{L}, a direct computation yields
\be\label{diss-0}
\ba
-\int\f{f\Lop_{[w,T]} f }{\Mloc}dv&=\int\big((v-w)f+T\na_vf\big)\cdot\na_v(\f{f}{\Mloc})dv\\
&=\int\f{|T\na_v f+(v-w)f|^2}{T\Mloc}dv\geq 0.
\ea
\ee
From the right-hand side of \eqref{diss-0}, we deduce that
\be\label{d-1}
\ba
\int\f{|T\na_v f+(v-w)f|^2}{T\Mloc}dv
&=\int\f{T^2|\na_vf|^2+|v-w|^2f^2}{T\Mloc}dv+\int\f{1}{\Mloc}(v-w)\cdot\na_v f^2dv\\
&=\int\f{T^2|\na_vf|^2+|v-w|^2f^2}{T\Mloc}dv-\int\dv_v(\f{v-w}{\Mloc})f^2dv\\
&=\int\f{T|\na_v f|^2-3f^2}{\Mloc}dv.
\ea
\ee
Moreover,
$$
\ba
\int\f{|T\na_v f+(v-w)f|^2}{T\Mloc}dv&=\int\f{1}{T\Mloc}(T\na_v f+\f{v-w}{2}f+\f{v-w}{2}f)^2dv\\
&=\int\f{1}{T\Mloc}|T\na_vf+\f{v-w}{2}f|^2dv+\int\f{|v-w|^2f^2}{4T\Mloc}dv\\
&+\int\f{f}{T\Mloc}(T\na_v f+\f{v-w}{2}f)\cdot(v-w)dv,
\ea
$$
where the last term can be estimated as
$$
\ba
\int\f{f}{T\Mloc}(T\na_v f+\f{v-w}{2}f)\cdot(v-w) dv
&=-\int \dv_v(\f{v-w}{2\Mloc})f^2dv+\int\f{|v-w|^2f^2}{2T\Mloc}dv
\\
&=-\int \f{f^2}{2T\Mloc}(3T+|v-w|^2)dv+\int\f{|v-w|^2f^2}{2T\Mloc}dv\\
&=-\int\f{3f^2}{2\Mloc}dv,
\ea
$$
which implies
\be\label{d-3}
\ba
\int\f{|T\na_v f+(v-w)f|^2}{T\Mloc}dv
&=\int\f{1}{T\Mloc}|T\na_vf+\f{v-w}{2}f|^2dv+\int\f{1}{\Mloc}(\f{|v-w|^2}{4T}-\f32)f^2dv\\
&\geq \int\f{1}{\Mloc}(\f{|v-w|^2}{4T}-\f32)f^2dv\\
&\geq \int\f{1}{\Mloc}\big(\f{|v|^2}{8T}-\f{|w|^2}{4T}-\f32\big)f^2dv.
\ea
\ee
Replacing $f$ by $\mb{P}_1f$ in \eqref{d-1} and \eqref{d-3}, and adding the resulting inequalities, we have
\be\label{d-6}
 2\int\f{|T\na_v(\mb{P}_1f)+(v-w)\mb{P}_1f|^2}{T\Mloc}dv\\
\geq
T_\star\int\f{|\na_v(\mb{P}_1f)|^2}{\Mloc}dv
+\int\f{1}{\Mloc}\big(\f{|v|^2}{8T^\star}-\f{|w|^2}{4T_\star}-\f92\big)|\mb{P}_1f|^2dv.
\ee

 {Next, since $\mb{P}_1f\in\mathfrak{N}^\bot$, the spectral bound \eqref{spec} of
Lemma \ref{Borel} together with \eqref{diss-0} gives}
\be\label{d-4}
\int\f{|T\na_v(\mb{P}_1f)+(v-w)\mb{P}_1f|^2}{T\Mloc}dv
\geq  {2}\int\f{|\mb{P}_1f|^2}{\Mloc}dv.
\ee
 {We emphasize that \eqref{d-4} is sharp: equality holds for $f=(v_1-w_1)(v_2-w_2)\Mloc$, which lies in $\mathfrak{N}^\bot$ and satisfies
$\Lop_{[w,T]}f=-2f$ by \eqref{conj}. In particular, no bound of the form $\geq T\int|\mb{P}_1f|^2\Mloc^{-1}dv$ is available when $T>2$, so the spectral estimate, and not a Poincar\'e inequality with constant one, is what must be used here.}

Multiplying \eqref{d-4} by {$(|w|^2+20T_\star)/8T_\star$}, adding the resulting inequality to \eqref{d-6}, we obtain
\be\label{d-5}
\ba
& {\f{|w|^2+36T_\star}{8T_\star}}\int\f{|T\na_v(\mb{P}_1f)+(v-w)\mb{P}_1f|^2}{T\Mloc}dv\\
& \geq T_\star\int\f{|\na_v(\mb{P}_1f)|^2}{\Mloc}dv
+\int\big(\f{|v|^2}{8T^\star}+\f12\big)\f{|\mb{P}_1f|^2}{\Mloc}dv\\
&\geq \f{1}{8}\min\{T_\star,\f{1}{T^\star},1\}      
\int\f{|\na_v(\mb{P}_1f)|^2+(1+|v|^2)|\mb{P}_1f|^2}{\Mloc}dv.
\ea
\ee
By \eqref{L} and \eqref{p0}, we have
$$
\ba
\Lop_{[w,T]}(\mb{P}_0f)&=\dv_v\big((v-w)\mb{P}_0f+T\na_v(\mb{P}_0f)\big)\\
&=\sum_{j=1}^3\pa_{v_j}\big((v_j-w_j)\mb{P}_0 f+T(-\f{v_j-w_j}{T}\mb{P}_0 f+\langle f,\chi_j\rangle\f{\chi_0}{\sqrt{T}}
+\langle f,\chi_4\rangle\f{2}{\sqrt{6T}}\chi_j)\big)\\
&=\sum_{j=1}^3\pa_{v_j}\big(\sqrt{T}\langle f,\chi_j\rangle\chi_0+\f{2\sqrt{T}}{\sqrt{6}}\langle f,\chi_4\rangle\chi_j\big)
\\
&=-\sum_{j=1}^3\langle f,\chi_j\rangle\chi_j-2\langle f,\chi_4\rangle\chi_4.
\ea
$$
Combining \eqref{p0}, \eqref{d-1} and \eqref{d-5}, and integrating by parts with respect to $v$, we obtain
\begin{align}\label{d-7}
-\int\f{f\Lop_{[w,T]} f}{\Mloc} dv &=-\int\f{\mb{P}_0f\Lop_{[w,T]}(\mb{P}_0f)}{\Mloc}dv
-\int\f{\mb{P}_1f\Lop_{[w,T]}(\mb{P}_1f)}{\Mloc}dv\\
\nm & =\sum_{j=1}^3|\langle f, \chi_j \rangle|^2+2|\langle f, \chi_4 \rangle|^2
+\int\f{|T\na_v(\mb{P}_1f)+(v-w)\mb{P}_1f|^2}{T\Mloc}dv\\
\nm & \geq {\min\{T_\star,\f{1}{T^\star},1\}\cdot\f{T_\star}{|w|^2+36T_\star}}
\int\f{|\na_v(\mb{P}_1f)|^2+(1+|v|^2)|\mb{P}_1f|^2}{\Mloc}dv,
\end{align}
where we used the  orthogonality relations $\mb{P}_0f\perp\Lop_{[w,T]}(\mb{P}_1f)$ and $\mb{P}_1f\perp\Lop_{[w,T]}(\mb{P}_0f)$. Choosing  {$\sigma(w,T_\star,T^\star)=\min\{T_\star,\f{1}{T^\star},1\}\cdot\f{T_\star}{|w|^2+36T_\star}$} in \eqref{d-7} yields the desired inequality \eqref{diss-1}. The proof of Lemma \ref{le-diss} is complete.

\hfill $\Box$

\end{appendix}

%%%%%%%%%%%%%%%%%%%%%%%%%%%%%%%%%%%%%%%%%%%%%%%%%%%%%%%%%%%%%%%%%%%%%%%%%%%%%%

\section*{Acknowledgment}
The research of H. Li was partially supported by the National Natural Science Foundation of China (No.11931010, 11871047 and 11225102).
T. Wang was partially supported by the Beijing Natural Science Foundation No. 1252001, the National Key R\&D Program of China No. 2022YFA1007700, the National Natural Science Foundation of China No. 12371215.

%%%%%%%%%%%%%%%%%%%%%%%%%%%%%%%%%%%%%%%%%%%%%%%%%%%%%%%%%%%%%%%%%%%%%%%%%%%%%%

\end{document}